\documentclass[11pt]{article}
\usepackage{amssymb}
\usepackage{amsmath}
\usepackage{amsfonts}
\usepackage{graphicx}
\usepackage{hyperref}
\usepackage{amsthm}
\usepackage{color}
\usepackage{natbib}
\bibpunct{(}{)}{;}{a}{,}{,} 
\theoremstyle{plain}
\newtheorem{theorem}{Theorem}
\newtheorem{lemma}{Lemma}
\newtheorem{corollary}{Corollary}
\newtheorem{proposition}{Proposition}

\theoremstyle{remark}
\newtheorem{remark}{Remark}

\newcommand{\II}{\mathcal{I}}
\newcommand{\trace}{\text{tr}}
\newcommand{\Ai}{\text{Ai}}
\renewcommand{\aa}{\alpha}
\newcommand{\zz}{\zeta}
\newcommand{\ww}{\omega}
\newcommand{\kk}{\kappa}
\newcommand{\ee}{\varepsilon}
\newcommand{\tee}{\tilde{\varepsilon}}
\newcommand{\tmu}{\tilde{\mu}}
\newcommand{\tsigma}{\tilde{\sigma}}
\newcommand{\taa}{\tilde{\alpha}}

\newcommand{\MM}{\mathbf{M}}
\newcommand{\EE}{\mathbf{E}}
\newcommand{\NN}{\mathbf{N}}
\newcommand{\labs}{\left|}
\newcommand{\rabs}{\right|}
\renewcommand{\to}{\rightarrow}

\newcommand{\Oh}[1]{\ensuremath{\mathrm{O}\left({#1}\right)}}

\newcommand{\convd}{\stackrel{\mathcal D}{\longrightarrow}}

\newcommand{\tauexpr}[1]{\tmu_{n,N} + {#1}\tsigma_{n,N}}
\newcommand{\xnexpr}[1]{\mu_{n,N} + {#1}\sigma_{n,N}}

\newcommand{\rmt}{Random Matrix Theory}

\newcommand{\tw}{Tracy-Widom law}

\title{Accuracy of the Tracy-Widom limit for the largest
  eigenvalue in white Wishart matrices}
\author{Zongming Ma\\
\emph{Stanford University}}
\begin{document}
\maketitle

\begin{abstract}
  Let $A$ be a $p$-variate real Wishart matrix on $n$ degrees of
  freedom with identity covariance. The distribution of the largest
  eigenvalue in $A$ has important applications in multivariate
  statistics. Consider the asymptotics when $p$ grows in proportion to
  $n$, it is known from \citet{johnstone01} that after centering and
  scaling, these distributions approach the orthogonal \tw{} for
  real-valued data, which can be numerically evaluated and tabulated in
  software.

  Under the same assumption, we show that more carefully chosen
  centering and scaling constants improve the accuracy of the
  distributional approximation by the Tracy-Widom limit to second order:
  \Oh{(n\wedge p)^{-2/3}}.
  Together with the numerical simulation, it implies that the \tw{} is
  an attractive approximation to the distributions of these largest
  eigenvalues, which is important for using the
  asymptotic result in practice. We also provide a parallel accuracy
  result for the smallest eigenvalue of $A$ when $n > p$.
\end{abstract}


\section{Introduction}
\label{sec:introduction}

The central object of multivariate statistical analysis is an
$n\times p$ data matrix $X$, where each of the $n$ rows corresponds
to an observation of a random vector in a $p$-dimensional space. If
we assume that the row vectors are i.i.d. samples from a
multivariate Gaussian distribution $N_{p}(\mu,\Sigma)$, much of the
classical theory in multivariate statistical analysis is reduced to
study of the eigen-decomposition of a random matrix following a
Wishart distribution. Typical examples include but are not limited
to principal component analysis (PCA), factor analysis and
multidimensional scaling (MDS). The fundamental setting is the
determinantal equation
\begin{equation*}
  \det(A-\lambda I) = 0\, ,
\end{equation*}
where $A$ follows a central Wishart distribution with covariance
matrix $\Sigma$.

In this setting, a common null hypothesis is $H_0: \Sigma = I$. For
instance, in PCA, this is the hypothesis of isotropic variation over
all the principal components; see, for example, \citet[Section
8.4.3]{mkb}. If $H_0$ is true, we say that we are in the null case and
call $A$ a (real) white Wishart matrix.  For testing this particular
hypothesis, as for many others in multivariate statistics, there are
two different systematic strategies: one is the likelihood ratio test
(LRT), which uses all the eigenvalues of $A$; the other is the union
intersection test (UIT) initiated by \citet{roy53}, which utilizes
only the largest (or smallest) eigenvalue of $A$ for the current
problem.


An inconvenience of using UIT is that the exact evaluation of the
marginal distribution of the extreme sample eigenvalues is not simply
tractable, even in the null case considered here. Interested readers
are referred to \citet[Section 9.7]{muirhead} for the expressions of
the marginal distributions in terms of hypergeometric function of
matrix argument; see, in particular, Corollary 9.7.2 and 9.7.4 there.
We remark that recent work of \citet{koev} has developed efficient
evaluations of hypergeometric functions of matrix argument and made
the computation of the exact marginal distributions possible when both
$n$ and $p$ are small.

An alternative approach is to approximate these exact finite sample
distributions of the extreme eigenvalues by some other well-understood
asymptotic distribution. This kind of approximation is ubiquitous in
statistics: the normal approximation to the distribution of the Wald
and score statistics, the Chi-square approximation to the Pearson
statistic in fitting contingency tables, etc. For the problem studied
here, \citet[Chapter 13]{anderson} provides a complete summary of the
established results in the conventional regime of asymptotics:
\begin{equation*}
\text{$p$ is fixed and $n\to \infty$.}
\end{equation*}
However, many modern data (microarray data, stock prices, weather
forecasting, etc.) we are now dealing with typically have the number
of features $p$ very large while the number of observations $n$ much
smaller than or just comparable to $p$. For these situations, the
classical asymptotics is no longer always appropriate and new
asymptotic results that could handle this type of data are desirable.


An advance in this direction was made in
\citet{johnstone01}, where the asymptotic regime was switched to
\begin{equation}
\label{eq:asymptotics}
p \to \infty, n = n(p)\to \infty \text{ and } n/p\to \gamma\in
(0,\infty).
\end{equation}
To state his result, let $X$ be an $n\times p$ data matrix with the
$n$ rows i.i.d. following $N_p(0,I)$. The $p\times p$ matrix $A =
X'X$ has a standard Wishart distribution: $A\sim W_p(I,n).$ We
denote the ordered eigenvalues of $A$ by $\lambda_1 \geq \cdots \geq
\lambda_p$. Borrowing tools from the field of \rmt{} (RMT),
especially those established in \citet{tw94, tw96, tw98}, Johnstone
showed that if we define centering and scaling constants as
\begin{equation}
\label{eq:orig-center-scale}
 \mu_{p} = \left(\sqrt{n-1} +
   \sqrt{p}\right)^2\, , \quad
 \sigma_{p} = \left(\sqrt{n-1} + \sqrt{p}\right)
   \left(\frac{1}{\sqrt{n-1}} +
     \frac{1}{\sqrt{p}}\right)^{1/3}\, ,
\end{equation}
then under condition \eqref{eq:asymptotics},
\begin{equation}
\label{eq:loe-lim}
\frac{\lambda_1 - \mu_{p}}{\sigma_{p}}\, \convd\,
W_1 \sim F_1\, ,
\end{equation}
where $F_1$ is the orthogonal \tw{}, which was originally found by
\citet{tw96} as the limiting law of the largest eigenvalue of a
$p\times p$ real Gaussian symmetric matrix. We remark that, prior to
\citet{johnstone01}, as a byproduct of his analysis on random growth
model, \citet{johansson00} established the scaling limit for the
largest eigenvalue in complex white Wishart matrix, which turns out
to be the unitary Tracy-Widom law $F_2$. We'd also like to mention
that for the weak limit \eqref{eq:loe-lim} to hold, \citet{nek-inf}
extended the asymptotic regime \eqref{eq:asymptotics} to include the
cases where $n/p\to 0$ or $\infty$.

This type of asymptotic result, albeit emerging only recently in
the statistics literature, has already found its relevance to
applications with modern data. For instance, based on the weak limit
\eqref{eq:loe-lim}, \citet{patterson} developed a formal test for
the presence of population heterogeneity in a biallelic dataset and
suggested a systematic way for assigning statistical significance to
successive eigenvectors, which in turn has been used to correct
population stratification \citep{price} and to perform genetic
matching \citep{luca} in genome-wide association studies.

From a statistical point of view, to inform the use of any
asymptotic result in practice, we need to have an understanding of
the accuracy of the approximation to finite distributions by the
limit, which usually appears in the form of a rate of convergence
result. In the complex domain, \citet{nek06} established such a
result for Johansson's theorem with carefully chosen centering and
scaling constants. With his choice, the error term in the
Tracy-Widom approximation could be controlled at the order
\Oh{(n\wedge p)^{-2/3}}, as opposed to \Oh{(n\wedge p)^{-1/3}} by
using the original centering and scaling constants in
\citet{johansson00}. For an up-to-date survey of higher order
accuracy results of this fashion, we refer to \citet[Section
3]{johnstone06}.


In statistics, we are typically more interested in real-valued data.
However, for technical reasons, results for complex-valued data are
usually easier to derive under the asymptotic regime
\eqref{eq:asymptotics} than in the real case. Recently, in analyzing
the parallel problem for the greatest root statistics for pairs of
Wishart matrices \citep[Definition 3.7.2]{mkb}, \citet{johnstone07}
figured out a way to connect the central object of study in the real
case to that in the complex case. To be more specific, in both real
and complex cases, the problem reduces to the study of operator
convergence in appropriate metrics by using standard techniques from
\rmt{}. The key observation there is that the crucial element of the
operator kernel in the real case could be represented in closed form
as a rank one perturbation of the complex kernel; see
\citet[Eq.(50)]{johnstone07}, which is a consequence of
\citet[Proposition 4.2]{adler00}.

Inspired by \citet{johnstone07}, we investigate in this paper the
rate of convergence for the distributions of properly rescaled
largest eigenvalues in real white Wishart matrices to the orthogonal
Tracy-Widom law. We remark that, instead of using \citet[Proposition
4.2]{adler00}, the central formula \eqref{eq:central} for the
``complex to real'' connnection in our paper is derived from a
slightly earlier result given in \citet[Section 4]{widom99} which is
specific to white Wishart matrices. This new approach not only helps
to avoid introducing a further nonlinear transformation after
rescaling the largest eigenvalues as in \citet{johnstone07} but also
enables us to make direct use of the analysis done in \citet{nek06}
for complex white Wishart matrices.

\vskip 1em
\textit{Statement of the theoretical result.} It was
suggested in \citet{johnstone06} that if we modify the centering and
scaling constants from \eqref{eq:orig-center-scale} to
\begin{equation}
  \label{eq:new-center-scale}
  \begin{split}
    \tmu_{np} & = \left(\sqrt{n-\tfrac{1}{2}} + \sqrt{p -
      \tfrac{1}{2}}\right)^2,\\
  \tsigma_{np} & = \left(\sqrt{n-\tfrac{1}{2}} + \sqrt{p -
      \tfrac{1}{2}}\right) \left(\frac{1}{\sqrt{n-\tfrac{1}{2}}} +
    \frac{1}{\sqrt{p - \tfrac{1}{2}}}\right)^{1/3},
  \end{split}
\end{equation}
we might obtain second order accuracy in the Tracy-Widom
approximation.

Indeed, the main theoretical result of the paper can be formulated as the
following theorem, which establishes the above conjecture.
\begin{theorem}
\label{thm:main}
Let $A\sim W_p(I,n)$ and $\lambda_1$ be its largest eigenvalue. Define
centering and scaling constants $(\tmu_{np}, \tsigma_{np})$ as in
\eqref{eq:new-center-scale}, then under condition
\eqref{eq:asymptotics}, there exists a continuous and nonincreasing
function $C(\cdot)$, such that for all real $s_0$, there exists an
integer $N_0(s_0,\gamma)$ for which we have that for any $s\geq s_0$
and $n\wedge p\geq N_0(s_0,\gamma)$,
\begin{equation*}
  \bigl|P\{\lambda_1\leq \tmu_{np} + \tsigma_{np}s\} - F_1(s)\bigr| \leq C(s_0)(n\wedge
  p)^{-2/3}\exp(-s/2)\, .
\end{equation*}
\end{theorem}


The theorem provides theoretical support for using the
Tracy-Widom law $F_1$ as approximate largest eigenvalue distribution
in the null case. In addition, the numerical investigation pursued
in Section \ref{subsec:numer-perf} shows that the approximation
yields reasonable accuracy even when $n$ and $p$ are as small as
$2$. Therefore, both theoretical and numerical results provide us
with the confidence in using the Tracy-Widom approximation for
nearly all finite $n\times p$ distributions, at least under the
Wishart assumption.

\begin{remark}
  In fact, Theorem \ref{thm:main} will be proved only when $p$ is even
  and $n \neq p$ since our method relies on a determinant formula of
  \citet{db55} which was only established for $p$ even and the
  Laguerre polynomials which are essential for building the
  convergence rate are not well-defined when $n = p$.  It would be of
  interest to have some theoretical support for the $p$ odd and the
  square cases. However, numerical experiments suggest that the
  Tracy-Widom approximation works just as well for $p$ odd as for $p$
  even and for $n= p$ as for $n \neq p$.
\end{remark}

\vskip 1em

\textit{Organization of the paper.} In Section \ref{sec:stat-impl},
we first investigate the numerical quality of the Tracy-Widom
approximation for finite $n\times p$ distributions, then review some
important statistical settings to which our result is relevant and
finally discuss several interesting issues involved in this study,
including a parallel result for the smallest eigenvalue. The rest of
the paper is dedicated to the proof of Theorem \ref{thm:main},
mainly with tools from Random Matrix Theory. In Section
\ref{sec:prelim}, we start with the formulation of Theorem
\ref{thm:main} in RMT terminology. After that, we derive the central
formula \eqref{eq:central} in this paper and reduce our problem to
the study of operator convergence in some appropriate metric. We
sketch our proof of the main result in Section \ref{sec:proof} by
assembling operator theoretic tools and asymptotic bounds on
transformed Laguerre polynomials. Finally, Section
\ref{sec:lagu-polyn-asympt} gives details of the Laguerre
asymptotics required in the proof. Appendix \ref{appendix} collects
various necessary technical details not spelled out fully in the
main text. Appendix \ref{appendix-b} discusses the issues mentioned
in Section \ref{subsec:other} in a more concrete manner.

\section{Statistical Implications and Discussion}
\label{sec:stat-impl}

\subsection{Quality of the approximation}
\label{subsec:numer-perf}

An important motivation for the current study is to promote
practical use of the Tracy-Widom approximation. For example, one
could tabulate the $F_1$ table and use it to compute $p$-values.
With such motivation, we investigate the quality of the
approximation with numerical experiments.

\paragraph{Distributional approximation.}

First of all, we study the numerical accuracy of the approximation
using our centering and scaling constants
\eqref{eq:new-center-scale} and compare it with that of the original
proposal \eqref{eq:orig-center-scale} in \citet{johnstone01}, with
results summarized in Table \ref{table:prob-tb}. We first look at
the square cases with $n = p = 2$, $5$, $20$ and $100$ and then the
cases with the same $p$'s but with the ratio $n/p$ fixed at $4:1$,
and finally the cases where $p = 5$ and $10$ with $n/p$ raised to as
high as $100:1$ and $1000:1$, which, in some sense, fall into the
situation $n/p\to \infty$ as discussed in \citet{nek-inf}. Finally,
in all these cases, we use $R = 40,000$ replications.

In terms of accuracy, from the last three columns of Table
\ref{table:prob-tb}, the approximation seems reasonable at
conventional significance levels of $10\%$, $5\%$ and $1\%$
(corresponding to right-hand tails of the distributions) even when
$p$ is as small as $2$ or $5$, and keeps improving as $p$ grows
large, regardless of the $n/p$ ratio. When $p$ is large, for
instance, in the $100\times 100$ and $400\times 100$ cases, the
Tracy-Widom law yields reasonable approximation over the whole range
of interest and matches the finite distributions almost exactly on
the right-hand tail.


In terms of the comparison with the original centering and scaling
constants, we could see from the first block of Table
\ref{table:prob-tb} that in the square cases, neither method seems
superior to the other. However, when the ratio $n/p$ is changed to
$4:1$ or larger (see the second and the third blocks of Table
\ref{table:prob-tb}), the improvement by using the new constants is
substantial. The new constants not only provide better absolute
accuracy in most of the cases, but also seem to result in a faster
convergence to the limiting distribution $F_1$.

Last but not least, the good performance on the right tail and the
faster convergence by using the new constants, as reflected in Table
\ref{table:prob-tb}, support our theoretical bound in Theorem
\ref{thm:main}.

\begin{table}[!tb]
\begin{center}
\begin{small}
\begin{tabular}{c|ccccccccc}
  \hline
  \hline
  Percentiles & $-3.8954$ & $-3.1804$ & $-2.7824$ & $-1.9104$ &  $-1.2686$ & $-0.5923$ & $0.4501$ & $0.9793$ &
  $2.0234$\\
  \hline
  TW  & .01 & .05 & .10 & .30 & .50 & .70 & \textbf{.90} & \textbf{.95} & \textbf{.99} \\
  \hline\hline
      $2\times 2$ & .000 & .000 & .000 & .034 & .379 & .690 & \textbf{.908} & \textbf{.953} & \textbf{.988} \\
                  & (.000) & (.000) & (.000) & (.015) & (.345) & (.669) & (.902) & (.950) & (.987) \\
      $5\times 5$ & .000 & .002 & .021 & .218 & .465 & .702 & \textbf{.908} & \textbf{.954} & \textbf{.989} \\
                  & (.000) & (.002) & (.020) & (.213) & (.460) & (.698) & (.907) & (.953) & (.989) \\
    $20\times 20$ & .003 & .029 & .071 & .275 & .490 & .700 & \textbf{.902} & \textbf{.952} & \textbf{.990} \\
                  & (.003) & (.029) & (.071) & (.274) & (.489) & (.699) & (.901) & (.952) & (.990) \\
  $100\times 100$ & .008 & .044 & .091 & .291 & .495 & .699 & \textbf{.901} & \textbf{.951} & \textbf{.990} \\
                  & (.008) & (.043) & (.091) & (.291) & (.495) & (.699) & (.901) & (.951) & (.990) \\
  \hline
      $8\times 2$ & .000 & .000 & .013 & .200 & .458 & .704 & \textbf{.913} & \textbf{.956} & \textbf{.990} \\
                  & (.000) & (.004) & (.031) & (.274) & (.534) & (.755) & (.931) & (.966) & (.992) \\
     $20\times 5$ & .001 & .018 & .057 & .262 & .486 & .703 & \textbf{.905} & \textbf{.952} & \textbf{.990} \\
                  & (.002) & (.028) & (.077) & (.305) & (.533) & (.739) & (.919) & (.960) & (.992) \\
    $80\times 20$ & .005 & .035 & .082 & .287 & .497 & .700 & \textbf{.902} & \textbf{.951} & \textbf{.990} \\
                  & (.006) & (.043) & (.096) & (.312) & (.524) & (.723) & (.911) & (.956) & (.992) \\
  $400\times 100$ & .009 & .047 & .095 & .298 & .499 & .700 & \textbf{.899} & \textbf{.949} & \textbf{.989} \\
                  & (.010) & (.052) & (.103) & (.312) & (.514) & (.712) & (.905) & (.952) & (.990) \\
  \hline
    $500\times 5$ & .010 & .050 & .100 & .303 & .502 & .705 & \textbf{.905} & \textbf{.953} & \textbf{.992} \\
                  & (.022) & (.084) & (.154) & (.387) & (.589) & (.770) & (.932) & (.968) & (.995) \\
   $5000\times 5$ & .012 & .056 & .108 & .311 & .511 & .711 & \textbf{.910} & \textbf{.957} & \textbf{.993} \\
                  & (.027) & (.098) & (.169) & (.406) & (.606) & (.783) & (.938) & (.971) & (.995) \\
  $1000\times 10$ & .010 & .049 & .099 & .296 & .500 & .701 & \textbf{.902} & \textbf{.952} & \textbf{.991} \\
                  & (.017) & (.073) & (.136) & (.363) & (.567) & (.754) & (.925) & (.964) & (.993) \\
 $10000\times 10$ & .012 & .054 & .104 & .306 & .506 & .707 & \textbf{.903} & \textbf{.950} & \textbf{.991} \\
                  & (.022) & (.084) & (.148) & (.381) & (.579) & (.764) & (.927) & (.964) & (.994) \\
   \hline\hline
$2\times\text{SE}$ & .001 & .002 & .003 & .005 & .005 & .005 &
.003 & .002 & .001\\
   \hline\hline
\end{tabular}
\end{small}
\end{center}
\caption{\label{table:prob-tb} Simulations for finite $n\times p$ vs. Tracy-Widom
  approximation: accuracy comparison of the new centering and scaling
  constants \eqref{eq:new-center-scale} with that in
  \citet{johnstone01}. For each combination of $n$ and $p$, we show in
  the first line the estimated cumulative probabilities for
  $\lambda_1$, rescaled using \eqref{eq:new-center-scale}; and in the
  second line with parentheses, rescaled using \citet[Eq.(1.3) and
  (1.4)]{johnstone01}, both computed from $R = 40,000$
  repeated draws using the Beta-ensemble sampling technique proposed
  by \citet{de02}. The conventional significance levels are
  highlighted in bold font and the last line gives approximate
  standard errors based on binomial sampling. The orthogonal Tracy-Widom
  distribution $F_1$ was computed using the method proposed in
  \citet{persson05} with percentiles obtained by inverse
  interpolation. }
\end{table}

\paragraph{Approximate percentiles.}

Except for computing $p$-values, $F_1$ could also be used to compute
approximate percentiles of finite $n\times p$ distributions. To
measure the accuracy of this approximation, we consider the relative
error $r_{\aa} = \theta^{TW}_{\aa}/\theta_{\aa}-1$, where
$\theta_{\aa}$ is the exact $100\aa$-th percentile of the largest
eigenvalue in finite $n\times p$ model and $\theta^{TW}_{\aa}$ is
its counterpart obtained from the Tracy-Widom law.

\begin{figure}[!tb]
    \centering
    \includegraphics[width = 3.05in]{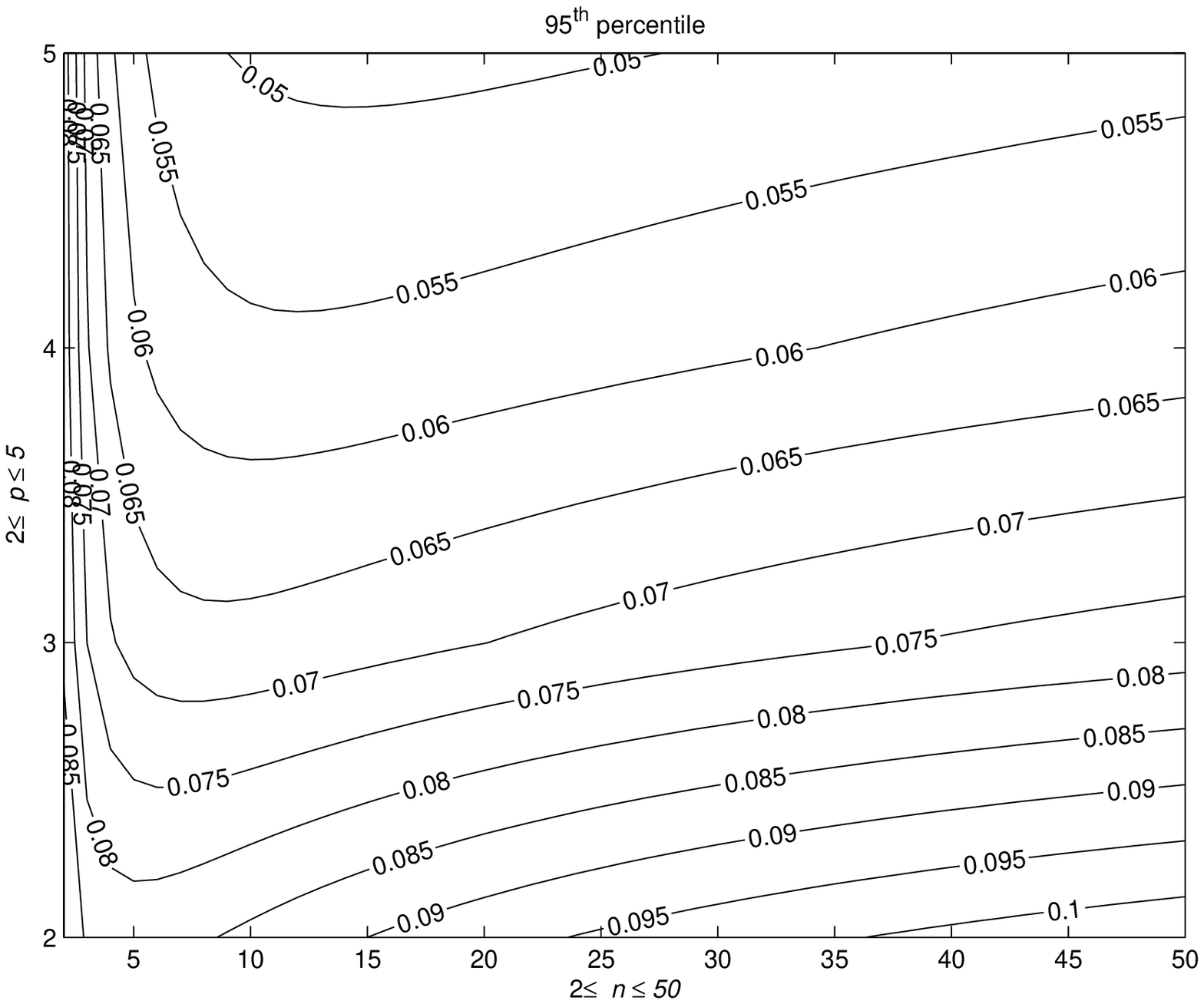}
    \includegraphics[width = 3.05in]{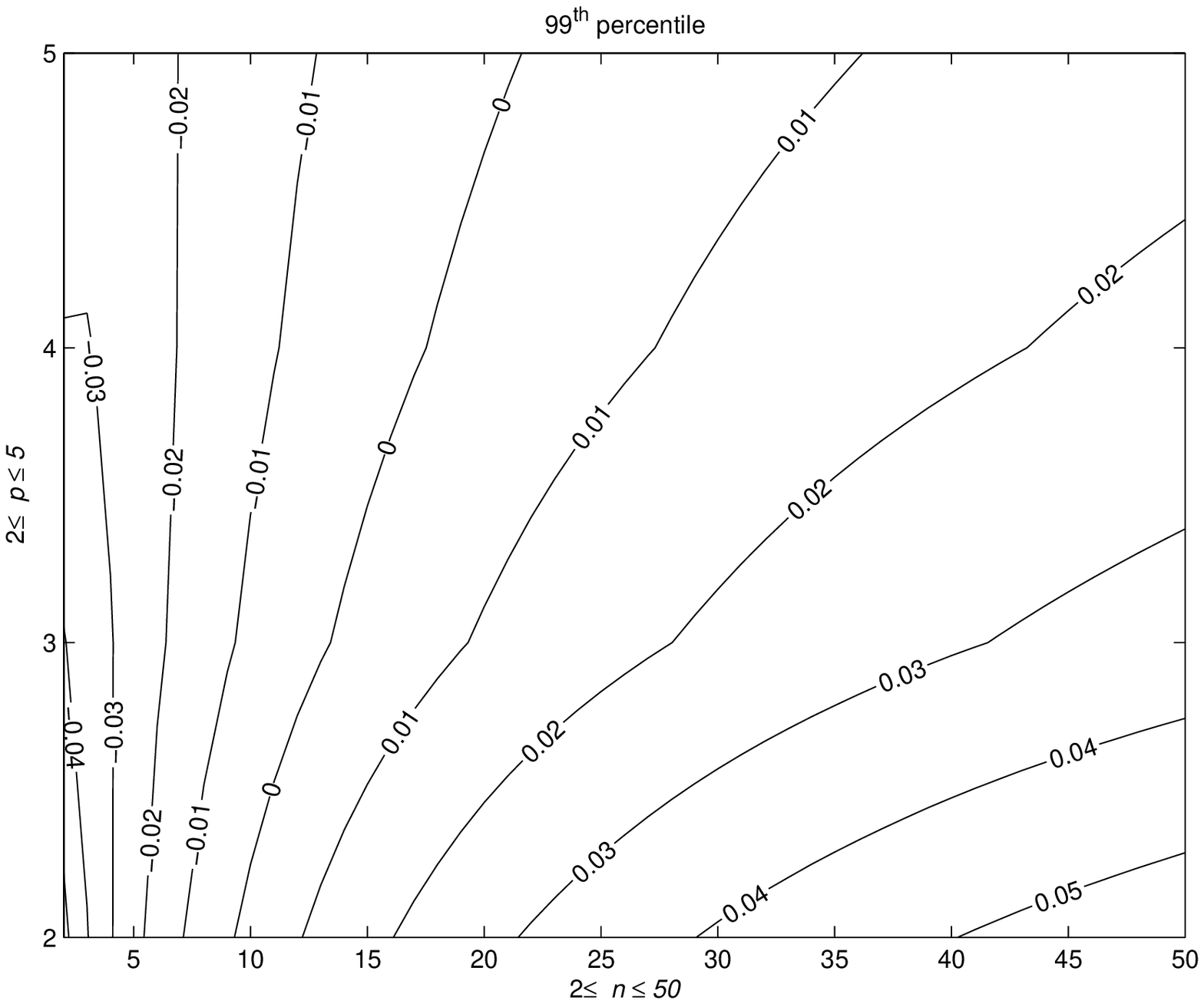}
    \caption{\label{fig:rr} Plots of the relative error
    $r_\aa = \theta^{TW}_\aa/\theta_\aa - 1$ for the approximate percentiles
    computed from $F_1$: (a) $95$-th percentile; (b) $99$-th percentile.
    The exact finite $n\times p$ largest
      eigenvalue distributions are computed using Plamen Koev's
      implementation in \textsc{matlab} of the recursive method
      proposed in \citet{koev} and the orthogonal Tracy-Widom law
      $F_1$ is computed using the method proposed in
      \citet{persson05}. The percentiles are always obtained
      from inverse interpolation.}
\end{figure}

In Figure \ref{fig:rr}, we plot the relative error $r_\aa$ for $\aa
= 0.95$ and $0.99$, with $p$ ranging from $2$ to $5$ and $n$ from
$2$ to $50$. Although the minimum of $n$ and $p$ is no larger than
$5$, the numerical accuracy is reasonably satisfactory. For the
$95$-th percentile, the relative error ranges from $5\%$ to $10\%$
for most of the cases and slightly exceeds $10\%$ only for the cases
where $p = 2$ and the $n/p$ ratios are high. The approximation to
the $99$-th percentile is even better, with the absolute relative
error $\labs r_{.99}\rabs \leq 5\%$ for most of the cases. Due to
the computational limitation \citep{koev}, we could not compute the
exact percentiles when $n$ and $p$ are large. However, we expect the
approximate percentiles to become more accurate as the consequence
of better distributional approximation.

\subsection{Related statistical settings}
\label{subsec:relat-stat-sett}

In this part, we review several common settings in multivariate
statistics to which our result is relevant.

\paragraph{Principal component analysis.} Suppose that $X =
[{x}_1,\cdots, {x}_n]'$ is a Gaussian data matrix.  Write the sample
covariance matrix ${S} = n^{-1}{X}'{H}{X},$ where ${H} =
{I}-n^{-1}\mathbf{1}\mathbf{1}'$ is the centering matrix, principal
component analysis looks for a sequence of standardized vectors
${a}_1, \cdots, {a}_p$ in $\mathbb{R}^p$, such that for $i =
1,\cdots, p$, where ${a}_i$ successively solves the following
optimization problem:
\begin{equation*}
  \max\{{a}'{S}{a}:\, {a}'{a}_j =
  0,\, j\leq i \}\, ,
\end{equation*}
where ${a}_0$ can be taken as the zero vector. The successive sample
principal component eigenvalues $\ell_1\geq \cdots \geq \ell_p$ then
satisfy $\ell_i = {a}'_i{S}{a}_i$. From a different perspective,
these $\ell_i$'s may also be found as the roots of the determinantal
equation
\begin{equation*}
  \det({S} - \lambda{I}) = 0\,.
\end{equation*}

One basic question in the application of PCA is testing the
hypothesis of isotropic variation, i.e., the hypothesis that all the
population principal component eigenvalues are equal. Under this
null hypothesis, the population covariance matrix of the row vectors
in ${X}$ is $\sigma^2{I}.$ For simplicity, let us suppose that
$\sigma^2 = 1$ (if $\sigma^2$ is an unknown value, we can estimate
it by some $\hat{\sigma}^2$ first and divide ${S}$ by
$\hat{\sigma}^2$). Then the sample covariance matrix ${S}$ satisfies
\begin{equation*}
  n{S}\sim W_p({I}, n-1).
\end{equation*}

The largest principal component eigenvalue $\ell_1$ of ${S}$ is a
natural test statistic for a union intersection test. Our result
applies for $n\ell_1$.

\paragraph{Multidimensional scaling.} Let ${X}$ be an
$n\times p$ data matrix. Consider the centered inner product matrix
${B} = {HXX'H}$, i.e. ${B}_{ij} = ({x}_i - \bar{{x}})'({x}_j -
\bar{{x}})$. In a typical setting of multidimensional scaling, we
are usually only given the matrix ${B}$ instead of the original
observations ${X}$. Let $\lambda_1 \geq \cdots \geq \lambda_p$ be
the ordered eigenvalues of ${B}$ and ${v}_{i}$ be the corresponding
eigenvector. As defined in \citet[Section 14.3]{mkb}: for fixed $k$
($1\leq k\leq p$), the rows of ${V}_k = ({v}_{1},\cdots, {v}_{k})$
are called the principal coordinates of ${X}$ in $k$ dimensions,
which constitute the classical $k$-dimensional solution to the
multidimensional scaling problem.

We observe that the matrix ${B}$ shares its non-zero eigenvalues
with $n{S} = {X'HX}$. For the principal coordinate method to make
sense, it is important that non-zero eigenvalues of ${B}$ and hence
all the eigenvalues of $n{S}$ do not equal a common value.
Translated to the population level, the population covariance matrix
${\Sigma} \neq \sigma^2 {I}.$ Assuming $\sigma^2 = 1$ (or dividing
${B}$ by $\sigma^2$ or its estimate $\hat{\sigma}^2$), the null
hypothesis can be written as $H_0: {\Sigma} = {I}$. As in the
situation of PCA, our result is useful for the test statistic
$\ell_1$, where $\ell_1$ is the largest eigenvalue of ${B}$.

\paragraph{Testing that a covariance matrix equals a specified
matrix.} Suppose that we have the Gaussian data matrix $X$ with the
rows $x_1, \cdots, x_n$ be independent $N_p(\mu,\Sigma)$ random
vectors and consider the null hypothesis $H: \Sigma = \Sigma_0$,
where $\Sigma_0$ is a specified positive definite matrix.

If the mean vector $\mathbf{\mu}$ is unknown, let ${S} =
n^{-1}{X}'{H}{X}$ be the sample covariance matrix. The union
intersection test uses the largest eigenvalue of the matrix
${\Sigma}_0^{-1}{S}$, denoted by $\lambda_1({\Sigma}_0^{-1}{S})$, as
the test statistic \citep[see][p.130]{mkb}.

We observe that $\lambda_1({\Sigma}_0^{-1}{S}) =
\lambda_1({\Sigma}_0^{-1/2}{S}{\Sigma}_0^{-1/2})$, where under the
null hypothesis,
$$n{\Sigma}_0^{-1/2}{S}{\Sigma}_0^{-1/2} \sim
W_p({I}, n-1).$$
Hence, our result is available for
$n\lambda_1({\Sigma}_0^{-1}{S})$.

\paragraph{Singular value decomposition.} For $X$ a real $n\times p$
matrix, there exists orthogonal matrices $U(n\times n)$ and
$V(p\times p)$, such that
\begin{equation*}
  X = UDV^T,
\end{equation*}
where $D = \text{diag}(\sigma_1,\cdots,\sigma_{\min(n,p)}) \in
\mathbb{R}^{n\times p}$, and $\sigma_1\geq \cdots\geq
\sigma_{\min(n,p)}\geq 0$. This representation is called the
singular value decomposition of $X$ [See \citet[Theorem
2.5.2]{golub}]. For $1\leq i\leq \min(n,p)$, $\sigma_i$ is called
the $i$-th singular value of $X$. Theorem \ref{thm:main} then
provides an accurate distributional approximation for $\sigma_1^2$
when the entries of $X$ are independent standard normal.

\subsection{Other issues}
\label{subsec:other}

For here, we provide brief remarks on several interesting issues
that we come across during the development of this work. More
details about them could be found in Appendix \ref{appendix-b}.

\paragraph{Transformation.} In the analysis of the greatest root
statistic, \citet{johnstone07} suggested that a nonlinear
transformation [the logit transformation: $\tau(x)=\log[x/(1-x)]$ in
his case] helps improve the distributional approximation by the
Tracy-Widom law, see Theorem 1, Table 1 and Fig. 1 there. In
addition to its numerical effect, the transformation has an
geometric explanation and yields a very natural integral
representation for the correlation kernel which later appears in the
central formula Eq.(50) there; see \citet[Section 2.2, also Eq.'s
(16) and (46)]{johnstone07}, \citet[Proposition 4.11]{pjfbook} and
\citet[Proposition 4.2]{adler00}.

Following \citet[Proposition 4.11]{pjfbook}, if we wanted to employ
a comparable transformation for our white Wishart case, it would be
the logarithmic transformation: $\tau(x) = \log x$. In fact, in our
study, we first looked into some depth along this direction and
could conclude the following second order accuracy result: under the
condition of Theorem \ref{thm:main}, let $\nu_{np} =
\log{\tmu_{np}}$ and $\tau_{np} =\tsigma_{np} / \tmu_{np}$, there
exists a continuous and nonincreasing function $C(\cdot)$, such that
for all real $s_0$, there exists an integer $N_0(s_0,\gamma)$ for
which we have that for any $s\geq s_0$ and $n\wedge p\geq
N_0(s_0,\gamma)$,
\begin{equation}\label{eq:res-log}
  \bigl|P\{\log\lambda_1\leq \nu_{np} + \tau_{np}s\} - F_1(s)\bigr| \leq C(s_0)(n\wedge
  p)^{-2/3}\exp(-s/2)\, .
\end{equation}

Some comments on how this result could be derived are included in
\ref{sec:log-trans}.

Although the rates of convergence are the same, numerical
experiments suggest that using the nonlinear transformation does not
yield as good numerical results in distributional approximation for
small to moderate $n$ and $p$ as simply rescaling $\lambda_1$ using
\eqref{eq:new-center-scale}, especially on the right-hand tail which
is of the most statistical interest. When $n$ and $p$ grow large,
using the transformation or not does not have as much influence,
as they approach the same limit.

In consideration of the actual quality of approximation, especially
for small to moderate $n$ and $p$, we suggest not using the
logarithmic transformation for the largest eigenvalues.  However, it
is of theoretical interest to know why such natural transformation
works for the greatest root statistic in \citet{johnstone07} but not
for the largest eigenvalue in white Wishart matrices here.

\paragraph{The smallest eigenvalue.} Following the principle of union
intersection tests, the smallest eigenvalue could also serve as the
test statistic in some cases, see, for instance, \citet[Section
5.2.2c]{mkb}. Hence, what we have established for the largest
eigenvalue is also worth investigation for the smallest one.
Moreover, understanding the deviation of the smallest eigenvalue
from its almost sure limit is also of independent interest. For
example, it plays an important role in the theory of sparse signal
recovery from large underdetermined linear system. See, for example,
\citet{donoho04} and \citet{candestao06}. In fact, as we studied the
accuracy result for the largest eigenvalue using the logarithmic
transformation, we obtained a parallel result for smallest
eigenvalues as a pleasant byproduct. We state without proof the
result here.

Suppose that $n - 1\geq p$ and $n/p \to \gamma \in (1,\infty)$ and
introduce the reflect Tracy-Widom law \citep{paul06} as
\begin{equation*}
G_1(s) = 1 - F_1(-s).
\end{equation*}
Let
\begin{equation*}
\mu_{np}^- = \left(\sqrt{n-\tfrac{1}{2}} -
\sqrt{p-\tfrac{1}{2}}\right)^2, \sigma_{np}^- =
\left(\sqrt{n-\tfrac{1}{2}} -
\sqrt{p-\tfrac{1}{2}}\right)\left(\frac{1}{\sqrt{p-\tfrac{1}{2}}} -
\frac{1}{\sqrt{n-\tfrac{1}{2}}}\right)^{1/3},
\end{equation*}
and then define
\begin{equation}
\label{eq:small-center-scale}
    \tau^-_{np} = \sigma^-_{np}/\mu^-_{np}, \quad \text{and} \quad \nu^-_{np} =
    \log\mu^-_{np} + \frac{1}{8}\left(\tau^-_{np}\right)^2.
\end{equation}
We then have that for the smallest eigenvalue $\lambda_p$ of a
$p\times p$ white Wishart matrix with $n$ degrees of freedom, there
exists a continuous and nondecreasing function $C(\cdot)$, such that
for all real $s_0$, there exists an integer $N_0(s_0,\gamma)$ for
which we have that for any $s\leq s_0$ and $p\geq N_0(s_0,\gamma)$,
\begin{equation}\label{eq:res-log-smallest}
  \bigl|P\{\log\lambda_p\leq \nu^-_{np} + \tau^-_{np}s\} - G_1(s)\bigr| \leq C(s_0)
  p^{-2/3}\exp(s/2)\, .
\end{equation}

See \ref{sec:smallest} for remarks on how to prove this result.

Unlike the case for $\lambda_1$, the logarithmic transformation
improves the numerical accuracy of the distributional approximation
for $\lambda_p$ significantly, especially when $p$ is small and
$n/p$ is close to $1$. We feel that an intuitive explanation to this
phenomenon could be the following: for $\lambda_p$, the lower bound
at $0$ strongly affects the approximation on the original scale,
especially when both $p$ and $n/p$ are small. However, by
transforming $\lambda_p$ to $\log\lambda_p$, one maps the lower
bound to $-\infty$ and hence avoids this `hard edge' effect. The
largest eigenvalue does not enjoy such a benefit for it does not
have an algebraic upper bound.

As a numerical illustration, in Table \ref{table:prob-tb-small}, we
present some simulation results on the Tracy-Widom approximation to
smallest eigenvalues transformed as above for two $n/p$ ratios:
$2:1$ and $4:1$, both with $p = 5$, $10$ and $100$. Again, for each
combination of $n$ and $p$, we run $R = 40,000$ replications. The
approximation seems good on the left-hand tail (where traditional
significance levels locate) even for $p$ as small as $5$, regardless
of the $n/p$ ratio. Moreover, for both $n/p$ ratios, when $p$ grows
to $100$, the approximation becomes reasonably accurate over the
entire range under investigation and is almost perfect on the
left-hand tail. Therefore, the numerical results agree well with the
theory for the smallest eigenvalues, too.

\begin{table}[!htb]
\begin{center}
\begin{small}
\begin{tabular}{c|ccccccccc}
  \hline
  \hline
  Percentiles & 3.8954 & 3.1804 & 2.7824 & 1.9104 & 1.2686 & 0.5923 & -0.4501 & -0.9793 &
  -2.0234\\
  \hline
  RTW  & .99 & .95 & .90 & .70 & .50 & .30 & \textbf{.10} & \textbf{.05} & \textbf{.01} \\
  \hline
  \hline
  $10\times 5$    &1.000 & .995 & .976 & .796 & .553 & .306 & \textbf{.093} & \textbf{.045} & \textbf{.012} \\
  $20\times 10$   & .999 & .984 & .952 & .760 & .536 & .305 & \textbf{.098} & \textbf{.049} & \textbf{.011} \\
  $200\times 100$ & .993 & .958 & .910 & .708 & .504 & .301 & \textbf{.099} & \textbf{.050} & \textbf{.010} \\
  \hline
  $20\times 5$    & .998 & .977 & .939 & .745 & .527 & .306 & \textbf{.097} & \textbf{.049} & \textbf{.010} \\
  $40\times 10$   & .996 & .969 & .926 & .726 & .511 & .300 & \textbf{.098} & \textbf{.048} & \textbf{.010} \\
  $400\times 100$ & .993 & .955 & .905 & .703 & .501 & .301 & \textbf{.100} & \textbf{.050} & \textbf{.010} \\
  \hline\hline
$2\times\text{SE}$ & .001 & .002 & .003 & .005 & .005 & .005 &
.003 & .002 & .001\\
   \hline\hline
\end{tabular}
\end{small}
\end{center}
\caption{\label{table:prob-tb-small} Simulations for finite $n\times p$ vs. Tracy-Widom
  approximation: the smallest eigenvalue. For each combination of $n$
  and $p$, the estimated cumulative probabilities are computed for
  $(\log\lambda_p - \nu^-_{np})/\tau^-_{np}$ with $R = 40,000$ draws
  from $W_p(I, n)$. The methods of sampling, computing $F_1$ and
  obtaining percentiles are the same as in Table
  \ref{table:prob-tb}. The conventional significance levels are
  highlighted in bold font and the last line gives approximate
  standard errors based on binomial sampling. }
\end{table}

\section{Random Matrix Theory}
\label{sec:prelim}

The establishment of Theorem \ref{thm:main} relies heavily on results
and methods from Random Matrix Theory (RMT) literature. In particular,
those about unitary and orthogonal Laguerre matrix ensembles play an
important role. In this section, we first restate our main result
using RMT terminology.  With a Lipschitz-type bound, we transform the
problem into the study of convergence rate of operators with matrix
kernels and derive the closed form representation \eqref{eq:central}
of the top-left entry in the kernel for Laguerre orthogonal
ensemble. Finally, we study the effect of scaling on our kernel
representation and carefully formulate the analysis problem to be
solved in later sections.

\subsection{Restatement of Theorem \ref{thm:main} in Random Matrix Theory}
\label{subsec:form-rand-matr}

Suppose $A$ is an $N\times N$ matrix following a $W_N(I,n)$
distribution with $n> N$. [Here and after, following the RMT
notational convention, we use $N$ rather than $p$ to denote the
number of features.] The celebrated joint probability density
function of the eigenvalues $x_1\geq \cdots \geq x_N\geq 0$ is given
by \citep{muirhead}:
\begin{equation*}
  p_N(x_1, \cdots, x_N) = d_{n,N}^{-1}\prod_{1\leq j<k\leq N}(x_j -
  x_k)\prod_{j=1}^N x_j^{n-N-1}e^{-x_j/2},
\end{equation*}
where $d_{n,N}$ is a normalizing constant depending only on $n$ and
$N$.

On the other hand, RMT people have investigated Laguerre Orthogonal
Ensembles (LOE), where `ensemble' stands for distribution of
matrices and `orthogonal' refers to the invariance of the
distribution under orthogonal transformations. The LOE($N,\taa$)
model ($\taa> -1$) has the matrix eigenvalue density as
\begin{equation}
  \label{eq:loe-pdf}
  \tilde{p}_N(x_1,\cdots, x_N) = d^{-1}_{\taa,N}\prod_{1\leq j<k\leq N}(x_j -
  x_k)\prod_{j=1}^N x_j^{\taa}e^{-x_j/2},
\end{equation}
where $x_1\geq \cdots\geq x_N\geq 0$ and $d_{\taa,N}$ is a normalizing
constant depending only on $N$ and $\taa$.

If we define $\aa_N = n - N$, the joint eigenvalue density of white
Wishart matrix $A$ is exactly the eigenvalue density of the LOE($N,
\aa_N - 1$) model. By this observation, we can formulate Theorem
\ref{thm:main} in terms of RMT as the following.

\begin{theorem}
  \label{thm:rmt}
  Let $x_1$ be the largest eigenvalue in the LOE($N,\aa_N-1$)
  model and $F_1$ be the orthogonal Tracy-Widom law. Define
  $\tmu_{n,N}$ and $\tsigma_{n,N}$ as
\begin{equation}
  \label{eq:tmu-tsigma}
    \tmu_{n,N} = \left(\sqrt{n - \tfrac{1}{2}} + \sqrt{N -
        \tfrac{1}{2}}\right)^2,\,
    \tsigma_{n,N} = \left(\sqrt{n - \tfrac{1}{2}} + \sqrt{N -
        \tfrac{1}{2}}\right)\left(\frac{1}{\sqrt{n - \tfrac{1}{2}}} +
      \frac{1}{\sqrt{N - \tfrac{1}{2}}}\right)^{1/3}.
\end{equation}
If $n > N$, $N\to \infty, n = n(N)\to \infty$ and $n/N\to \gamma\in
[1,\infty)$, there exists a continuous and nonincreasing
function $C(\cdot)$, such that for all real $s_0$, there exists an
integer $N_0(s_0,\gamma)$ for which we have that for any $s\geq s_0$
and $N\geq N_0(s_0,\gamma)$,
\begin{equation*}
  \bigl|P\{x_1\leq \tmu_{n,N} + \tsigma_{n,N}s\} - F_1(s)\bigr| \leq C(s_0)
  N^{-2/3}\exp(-s/2)\, .
\end{equation*}
\end{theorem}

\begin{remark}
  The theorem is stated only for situations where $n > N$. It
  works equally well when $n < N$ by switching $n$ and $N$. This results
  from the following observations: (a) constants in \eqref{eq:tmu-tsigma}
  are symmetric in $n$
  and $N$ and (b) switching $n$ and $N$ does not change the
  distribution of $x_1$.
\end{remark}

\subsection{Operator determinant and kernel representation}
\label{sec:determ-form-orth}

We focus on the LOE($N,\taa$) model in \eqref{eq:loe-pdf} for the
moment. For general orthogonal ensembles, \citet[Section 9]{tw98}
showed that when $N$ is even, for $\chi = I_{x>x'}$:
\begin{equation}
\label{eq:F-N-1}
  F_{N,1}(x') \equiv P\{x_1\leq x'\} = \sqrt{\det(I-K_N\chi)},
\end{equation}
with $K_N$ an operator with a $2\times 2$ matrix kernel:
\begin{equation}
\label{eq:K-N-kernel}
  K_{N}(x,y) = \begin{pmatrix} I & -\partial_2 \\ \ee_1 & T
\end{pmatrix} S_{N,1}(x,y) - \begin{pmatrix} 0 & 0 \\ \ee(x-y) & 0
\end{pmatrix}\ ,
\end{equation}
where $\partial_2$ is the differential operator with respect to the
second argument, $\ee_1$ is the convolution operator acting on the
first argument with the kernel $\ee(x-y)=\frac{1}{2}\text{sgn}(x-y)$
and $TS(x,y) = S(y,x)$ for any kernel $S$. However, no explicit
representation of $S_{N,1}$ was given there.

In a follow-up paper, \citet{widom99} derived explicit expression of
the kernel $S_{N,1}$ for Gaussian and Laguerre orthogonal ensembles,
which is summarized in \citet[Eq.(4.3)]{adler00} in a more friendly
form. In particular, for the LOE($N,\taa$) model of our interest, we
have [Warning: we need to switch $x$ and $y$ in
\citet[Eq.(4.3)]{adler00}.]:
\begin{equation}
  \label{eq:widom-S1}
  \begin{split}
    S_{N,1}(x,y)  = S_{N,2}(x,y) + & \frac{N!}{4\Gamma(N+\tilde{\aa})}
    x^{\tilde{\aa}/2} e^{-x/2} \left(\frac{d}{dx} L_N^{\tilde{\aa}}(x)
    \right)\\
    & \times \int_0^\infty
    \text{sgn}(y-z) z^{\tilde{\aa}/2-1} e^{-z/2} [L_{N}^{\tilde{\aa}}(z) - L_{N-1}^{\tilde{\aa}}(z)]dz,
  \end{split}
\end{equation}
where $L_k^{\taa}$ $(k=N-1, N)$ are Laguerre polynomials defined in
\citet[Chapter V]{szego} and $S_{N,2}(x,y)$ is the kernel related to
the Laguerre unitary ensemble (LUE) with parameter ($N,\taa$), which
has the following eigenvalue density:
\begin{equation*}
  p_{N}(x_1,\cdots, x_N) = c_{n,N}^{-1}\prod_{1\leq j<k\leq N}(x_j -
  x_k)^2\prod_{j=1}^N x_j^{\tilde{\aa}}e^{-x_j},\quad
  x_1\geq\cdots\geq x_N\geq 0.
\end{equation*}

With \eqref{eq:widom-S1}, we start to derive an closed form
representation for $S_{N,1}$ after some necessary definitions. As in
\citet{johnstone01}, we define a basis $\{\phi_k\}_{k=0}^\infty$ on
$L^2([0,\infty))$ with transformed Laguerre polynomials
\begin{equation}
  \label{eq:phi-k}
  \phi_k(x;\taa) = \sqrt{\frac{k!}{(k+\taa)!}}\,x^{\taa/2}e^{-x/2}L_k^{\taa}(x).
\end{equation}
Then calling $a_N = \sqrt{N(N+\taa)}$, we follow \citet[Section
2]{nek06} to introduce for $x\geq 0$,
\begin{equation}
  \label{eq:phi-psi}
  \phi(x;\taa) = (-1)^N\sqrt{\frac{a_N}{2}}\phi_N(x;\taa - 1)x^{-1/2};\quad
  \psi(x;\taa) = (-1)^{N-1}\sqrt{\frac{a_N}{2}}\phi_{N-1}(x;\taa +
  1)x^{-1/2}.
\end{equation}

With the definition in \eqref{eq:phi-psi}, for the first term in
\eqref{eq:widom-S1}, \citet[Eq.(3.6)]{johnstone01} and \citet[Appendix
A.5]{nek06} gave the following integral representation
\begin{equation*}
  S_{N,2}(x,y) = \int_0^\infty \phi(x+z)\psi(y+z) + \psi(x+z)\phi(y+z) dz.
\end{equation*}
For the second term, we could apply \citet[Eq.(5.1.13),
(5.1.14)]{szego} to obtain that it equals
\begin{align*}
  -\frac{N!}{4\Gamma(N+\tilde{\aa})}x^{\tilde{\aa}/2}e^{-x/2}L_{N-1}^{\tilde{\aa}+1}(x)
  \int_0^\infty \text{sgn}(y-z)z^{\tilde{\aa}/2 - 1}e^{-z/2}L_N^{\tilde{\aa}-1}(z)dz
  =  \psi(x)\int_0^\infty \ee(y-z)\phi(z)dz.
\end{align*}
Hence, we obtain
\begin{equation*}
  S_{N,1}(x,y) = S_{N,2}(x,y) + \psi(x)\int_0^\infty \ee(y-z)\phi(z)dz.
\end{equation*}

Recall that for the white Wishart matrix $A\sim W_N(I,n)$, setting $\aa_N = n-N$,
it is connected to the LOE($N,\taa$) model by the identity $\taa =
\aa_N - 1$. Thus, if we use the parameters $N$ and $\aa_N$, then the
above calculation gives the following representation for $S_{N,1}$:
\begin{equation}
  \label{eq:central}
\boxed{
  S_{N,1}(x,y;\aa_N - 1) = S_{N,2}(x,y;\aa_N-1) +
  \psi(x;\aa_N-1)\int_0^\infty \ee(y-z)\phi(z;\aa_N-1)dz.
}
\end{equation}

\subsubsection{Framework for deriving the determinant formula}
\label{sec:det}

The determinant formula \eqref{eq:F-N-1} introduced at the beginning
of this subsection provides the foundation for the convergence
arguments. However, it is worth clarification under which framework it
is derived.

\citet{tw05} described with care the operator convergence of $K_N\chi$
to the limit $K_{GOE}$ for the Hermite finite $N$ ensemble. We adapt
and extend their approach to the Laguerre finite $N$
ensemble. Therefore, we paraphrase their remarks on the weighted
Hilbert spaces and regularized $2$-determinants under the current
setting.

In the kernel $K_N$ given in \eqref{eq:K-N-kernel}, the first term on
the right hand side has each of its entries finite rank operators and
hence a trace class operator. However, this is not true for
$\ee(x-y)$. According to \citet[Theorem VI.23]{reed&simon80}, it is
even not Hilbert-Schmidt on $L^2([x',\infty))$. One way to take care
of this problem is to introduce the weighted $L^2$ space and to
generalize the operator determinant as in \citet{tw05}.

To this end, let $\rho$ be any weight function which satisfies the
following two conditions:
\begin{enumerate}
\item[(1)] its reciprocal $\rho^{-1}\in L^1([0,\infty))$; and
\item[(2)] each operator that constitutes elements in the first term
  on the right hand side of \eqref{eq:K-N-kernel} is in
  $L^2([x',\infty);\rho)\cap L^2([x',\infty);\rho^{-1})$.
\end{enumerate}
Then, as remarked in \citet{tw05}, $\ee$: $L^2([x',\infty);\rho)\to
L^2([x',\infty);\rho^{-1})$ is Hilbert-Schmidt. Moreover, $K_N$
could now be regarded as a $2\times 2$ matrix kernel on the space
$L^2([x',\infty);\rho)\oplus L^2([x',\infty);\rho^{-1})$.

We have thus made clear on which space the kernel $K_N$ acts. In
order for the determinant formula \eqref{eq:F-N-1} to hold, we need
a generalization of the usual Fredholm determinant for trace class
operators to determinant for Hilbert-Schmidt operators.

By our condition on $\rho$, for $K_N = [K_{ij}]_{1\leq i,j\leq 2}$, we
regard $K_{11}$ and $K_{22}$ as trace class operators on
$L^2([x',\infty);\rho)$ and $L^2([x',\infty);\rho^{-1})$
respectively and off-diagonal elements as Hilbert-Schmidt operators:
\begin{equation*}
  K_{12}: L^2([x',\infty);\rho^{-1})\to
  L^2([x',\infty);\rho)\quad \text{and}\quad K_{21}: L^2([x',\infty);\rho)\to
  L^2([x',\infty);\rho^{-1}).
\end{equation*}

Hence, $\trace(K_N) = \trace(K_{11}) + \trace(K_{22})$ is well
defined. The regularized 2-determinant of Hilbert-Schmidt operator $T$
with eigenvalues $\mu_k$ is defined by
\begin{equation*}
  {\det}_2(I-T) = \prod_{k}(1-\mu_k)e^{\mu_k}.
\end{equation*}
Then one naturally extends the operator definition of determinants to
Hilbert-Schmidt operator matrix $T$ with trace class diagonal entries
by setting
\begin{equation}
  \label{eq:det}
  \det(I-T) = {\det}_2(I-T)\exp(-\trace T).
\end{equation}

Finally, as remarked in \citet{tw05}, the resulting notion of
$\det(I-K_N)$ is independent of the choice of $\rho$ and allows the
derivation in \citet{tw98} that yields \eqref{eq:F-N-1},
\eqref{eq:K-N-kernel} and eventually \eqref{eq:central}.

Later in Section \ref{sec:rho-choice}, we will make a specific
choice of $\rho$, which not only makes our arguments more explicit but
also eases the derivation of the right tail exponential decay in our
desired bound.

\subsection{Scaling the kernel}
\label{subsec:repr-orth-kern}

Fixing any real number $s_0$ and introducing the linear
transformation $\tau(s) = \tauexpr{s}$, we are interested in the
convergence rate of $F_{N,1}(\tau(s'))$ to $F_1(s')$ for all $s'\geq
s_0$.

Define the rescaled kernel $K_\tau$ as the following:
\begin{equation*}
K_\tau(s,t) = \sqrt{\tau'(s)\tau'(t)}\,K_N(\tau(s),\tau(t))
=\tsigma_{n,N}K_N(\tauexpr{s},\tauexpr{t}).
\end{equation*}
We have $\det(I-K_N) = \det(I-K_\tau)$ by noticing that $K_N$ and
$K_\tau$ share the spectrum. We give below an explicit
representation of $K_\tau$ for later use.

Before we proceed, we apply the $\tau$-scaling to $\phi$,
$\psi$ and $S_{N,2}$ and thus define
\begin{equation}
\label{eq:phi-psi-tau}
  \phi_\tau(s) = \tsigma_{n,N}\phi(\tauexpr{s}),\quad
  \psi_\tau(s) = \tsigma_{n,N}\psi(\tauexpr{s})
\end{equation}
and
\begin{equation}
\begin{split}
\label{eq:s-tau}
  S_\tau(s,t)  &= \tsigma_{n,N}S_{N,2}(\tauexpr{s},\tauexpr{t})\\
   &= \int_0^\infty \phi_\tau(s+z)\psi_\tau(t+z) + \psi_\tau(s+z)\phi_\tau(t+z)dz.
\end{split}
\end{equation}
For later convenience, $\phi_\tau(s)$ and $\psi_\tau(s)$ are assumed
to be $0$ when $\tau(s) = \tauexpr{s} < 0$, and hence they are
well-defined on the entire real line.

 Finally, we introduce the short notation
\begin{equation}
  \label{eq:S-tau-R}
  S_\tau^R(s,t) = S_\tau(s,t) +
  \psi_\tau(s)\int_{-\infty}^\infty\ee(t-z)\phi_\tau(z)dz = S_\tau(s,t) +
  \psi_\tau(s)\left(\ee\phi_\tau\right)(t).
\end{equation}
[We remind the reader that in the above discussion, we have dropped
the explicit dependence on $\taa$ or $\aa_N-1$ to avoid notation
nightmare. Henceforth, we mention the explicit dependence only for
eliminating ambiguity.]

We further observe that the determinant formula does not change if
we modify $K_\tau$ as
\begin{equation*}
  K_\tau(s,t) = \begin{pmatrix} K_{\tau,11}(s,t) &
    \tsigma_{n,N} K_{\tau,12}(s,t) \\
    \tsigma_{n,N}^{-1}K_{\tau,21}(s,t) & K_{\tau,22}(s,t)
  \end{pmatrix}\ ,
\end{equation*}
for the spectrum does not change. Based on this observation and our
detailed calculation in \ref{subsec:a-k-tau}, we could represent the
entries of $K_\tau$ as
\begin{equation}
  \label{eq:k-tau-entries}
  \begin{split}
    K_{\tau,11}(s,t) = S^R_\tau(s,t), \qquad &
    K_{\tau.12}(s,t) = -\partial_t S^R_\tau(s,t), \\
    K_{\tau,21}(s,t) = (\ee_1 S^R_\tau)(s,t) - \ee(s-t),\qquad &
    K_{\tau,22}(s,t) = K_{\tau,11}(t,s).
  \end{split}
\end{equation}

For the desired limit $F_1(s')$ of the sequence $F_{N,1}(s')$,
\citet{tw05} showed that $F_1(s') = \sqrt{\det(I-K_{GOE})}$, where
the operator $K_{GOE}$ has the matrix kernel
\begin{equation*}
    K_{GOE}(s,t) = \begin{pmatrix}
        S(s,t) & SD(s,t)\\
        IS(s,t)-\ee(s-t) & S(t,s)
    \end{pmatrix}
\end{equation*}
with the entries given by
\begin{equation*}
    \begin{split}
        S(s,t) = S_A(s,t) + \frac{1}{2}\Ai(s)\left( 1-\int_t^\infty \Ai(u)du
        \right),\quad
        SD(s,t) = -\partial_t S_A(s,t) -
        \frac{1}{2}\Ai(s)\Ai(t)\\
        \text{and}\quad
        IS(s,t) = -\int_s^\infty S_A(u,t)du +
        \frac{1}{2}\left( \int_t^s \Ai(u)du + \int_s^\infty \Ai(u)du \int_t^\infty \Ai(u)du
        \right).
    \end{split}
\end{equation*}
Here $S_A(s,t) = \int_0^\infty \Ai(s+u)\Ai(t+u)du$ represents the
Airy kernel with $\Ai(\cdot)$ the Airy function defined in
\citet[p.53, Eq.(8.01)]{olver74}.

By our discussion in Section \ref{sec:det}, it is necessary that both
$K_\tau$ and $K_{GOE}$ belong to the
following class $\mathcal{A}$ of operators
\begin{equation}
\label{eq:op-class}
  \begin{split}
    \mathcal{A} \equiv \{& \text{$2\times 2$ Hilbert-Schmidt operator
      matrices $A$ on}  \\
    & \quad \text{$L^2([s',\infty);\rho\circ\tau)\oplus
    L^2([s',\infty);\rho^{-1}\circ\tau)$ with trace class diagonal entries} \}.
  \end{split}
\end{equation}
This fact will be verified after we choose a specific $\rho$ in
Section \ref{sec:rho-choice}. For the convenience of argument, let
us assume it for the moment.

\subsection{Lipschitz bound and kernel difference}
\label{sec:repr-kern-diff}

Let $p_N = F_{N,1}(s')$ and $p_\infty = F_1(s')$, we note that $
|p_N - p_\infty| \leq |p_N^2 - p_\infty^2|/p_\infty \leq
  C(s_0)|p_N^2 - p_\infty^2|$,
where $C(s_0) = 1/F_1(s_0)$ which is continuous and non-increasing in
$s_0$. Thus, we are led to the difference of the determinants
\begin{equation}
\label{eq:diff-det}
  \labs F_{N,1}(s') - F_1(s') \rabs \leq C(s_0)\labs \det(I-K_\tau) -
  \det(I - K_{GOE}) \rabs.
\end{equation}

\begin{remark}
Here and after, we use $C(s_0)$ to denote in general any continuous
and non-increasing function of $s_0$ and $C$ any universal constant,
where the actual function and constant might be different from
display to display.
\end{remark}

To study the quantity on the right hand side of \eqref{eq:diff-det},
our basic tool is the following Lipschitz-type bound on the matrix
operator determinant for operators in $\mathcal{A}$.
\begin{proposition}
\label{prop:det-bd} For operators $A$ and $B$ in class $\mathcal{A}$
and determinants of $I-A$ and $I-B$ defined as in \eqref{eq:det}, if
$\sum_i\|A_{ii} - B_{ii}\|_1 + \sum_{i\neq j}\|A_{ij} -
B_{ij}\|_2\leq 1/2$, then
\begin{equation}
\label{eq:det-bd}
  \labs \det(I-A) - \det(I-B) \rabs \leq
  M(B)\left(\sum_{i}\|A_{ii}-B_{ii}\|_1 + \sum_{i\neq j}\|A_{ij} -
  B_{ij}\|_2 \right),
\end{equation}
where $M(B) = 2\,\labs\det(I-B)\rabs + 2\exp\left[2\left(1+
\|B\|_2\right)^2 + \sum_i\|B_{ii}\|_1 \right].$
\end{proposition}

Note that the leading term on the right hand side of
\eqref{eq:det-bd} depends only on $B$. In this sense, Proposition
\ref{prop:det-bd} is a refinement of Proposition 3 in
\citet{johnstone07}. Its proof could be found in
\ref{subsec:a-proof}.

By Proposition \ref{prop:det-bd}, if we could control the entry-wise
convergence rate of $K_\tau$ to $K_{GOE}$, we will be able to bound
the right hand side of \eqref{eq:diff-det} and hence prove our
theorem. To this end, a convenient expression of the kernel
difference $K_\tau - K_{GOE}$ is helpful. we derive such an
expression below by essentially adapting the arguments in
\citet[Section 8.3]{johnstone07} to the current context.

According to \citet[Eq.(4.2)]{nagao95}, we could calculate [see
\ref{subsec:a-constants} for detail] that when $N$ is even,
\begin{equation*}
\begin{split}
\int_{-\infty}^\infty \psi_\tau(s;\aa_N-1)ds & = 0, \quad
    \text{and}\\
  \int_{-\infty}^\infty \phi_\tau(s;\aa_N-1)ds & =
  \frac{N^{1/4}(n-1)^{1/4} }
    {2^{(\aa_N-3)/2}(N+1)
    }\frac{\Gamma\left(\frac{N+3}{2}\right)}{\Gamma\left(\frac{n+1}{2}\right)}
    \left[\frac{\Gamma(n)}{\Gamma(N+1)}\right]^{1/2}.
\end{split}
\end{equation*}
For later use, we define $\beta_N = \frac{1}{2}\int_{-\infty}^\infty
\phi_\tau(s)ds$.

Bring in the right-tail integration operator $\tee$ introduced in
\citet[Section 8.3]{johnstone07} as
\begin{equation}\label{eq:tee}
  (\tee{g})(s) \equiv \int_s^\infty g(u)du,
\end{equation}
we have the identity $(\ee{g})(s) = \frac{1}{2}\int_{-\infty}^\infty
g(u)du -(\tee{g})(s)$ and hence obtain
\begin{equation*}
    \ee\phi_\tau = \beta_N - \tilde{\ee}\phi_\tau,\quad
    \text{and}\quad \ee\psi_\tau = -\tilde{\ee}\psi_\tau.
\end{equation*}

For $S_\tau$ defined in \eqref{eq:s-tau}, by Fubini's theorem
[justified by Lemma \ref{lemma:laguerre}]
\begin{equation*}
  \int_{-\infty}^\infty S_\tau(u,t)du = 2\beta_N \int_0^\infty
  \psi_\tau(t+z)dz = 2\beta_N(\tilde{\ee}\psi_\tau)(t).
\end{equation*}
Observing that for any kernel $A(s,t)$, $ (\ee_1 A)(s,t) =
\frac{1}{2}\int_{-\infty}^\infty A(u,t)du -
  \int_s^\infty A(u,t)du$, and introducing the abbreviation $a\otimes
  b$ for rank one operator with kernel $a(s)b(t)$, we have
$\ee_1 S_\tau = \beta_N \otimes \tilde{\ee}\psi_\tau -
\tilde{\ee}_1S_\tau$, and for $S_\tau^R$ in \eqref{eq:S-tau-R}, we
have $ S^R_\tau = S_\tau + \psi_\tau\otimes\beta_N -
\psi_\tau\otimes \tilde{\ee}\phi_\tau, $ which finally gives
\begin{equation*}
  \ee_1 S^R_\tau = -\tilde{\ee}_1\left(S_\tau - \psi_\tau \otimes
  \tilde{\ee}\phi_\tau\right) + \beta_N\left(1\otimes \tilde{\ee}\psi_\tau -
  \tilde{\ee}\psi_\tau\otimes 1\right).
\end{equation*}

By the explicit expressions for $K_\tau$ entries in
\eqref{eq:k-tau-entries},
\begin{equation*}
  K_\tau(s,t) =
  \begin{pmatrix}
    S^R_\tau(s,t) & -\partial_t S^R_\tau(s,t) \\
    \ee_1 S^R_\tau(s,t) & S^R_\tau(t,s)
  \end{pmatrix} +
  \begin{pmatrix}
    0 & 0 \\
    -\ee(s-t) & 0
  \end{pmatrix}
  = LS^R_\tau + K^\ee,
\end{equation*}
where
\begin{equation*}
  L =
  \begin{pmatrix}
    I & -\partial_2 \\
    \ee_1 & T
  \end{pmatrix}
  \quad \text{and} \quad
  K^\ee =
  \begin{pmatrix}
    0 & 0 \\
    -\ee & 0
  \end{pmatrix}.
\end{equation*}

We then decompose $K_\tau$ and $K_{GOE}$ as follows:
\begin{equation}
\label{eq:kern-decomp}
  K_\tau = K^R_\tau + K^F_{\tau,1} + K^F_{\tau,2} + K^\ee \quad\text{and}\quad
  K_{GOE} = K^R + K^F_1 + K^F_2 + K^\ee,
\end{equation}
where by defining $G = \Ai/\sqrt{2}$ and the matrix kernels
$\tilde{L} =
  \begin{pmatrix}
    I & -\partial_2 \\
    -\tee_1 & T
  \end{pmatrix}$,
$  L_1 =
  \begin{pmatrix}
    I & 0 \\
    -\tee_1 & 0
  \end{pmatrix}$, and
$  L_2 =
  \begin{pmatrix}
    0 & 0 \\
    \tee_2 & I
  \end{pmatrix}$, we
could write down the unspecified components in
\eqref{eq:kern-decomp} explicitly as
\begin{equation}
\label{eq:kern-comp}
  \begin{aligned}
    K^R_\tau = \tilde{L}[S_\tau -
  \psi_\tau\otimes\tilde{\ee}\phi_\tau], &\quad &
  K^F_{\tau,1} = \beta_N L_1[\psi_\tau(s)], &\quad &
  K^F_{\tau,2} = \beta_N L_2[\psi_\tau(t)], \\
  \text{and}\quad K^R = \tilde{L}[S_A - G\otimes \tilde{\ee}G], &\quad &
  K^F_1 = \frac{1}{\sqrt{2}}L_1[G(s)], &\quad &
  K^F_2 = \frac{1}{\sqrt{2}}L_2[G(t)].
  \end{aligned}
\end{equation}

For $\Delta_N$ to be defined in \eqref{eq:Delta-N}, we will
establish in Lemma \ref{lemma:laguerre} that $\phi_\tau = G +
\Delta_N G' + \Oh{N^{-2/3}}$, so set $G_N = G + \Delta_N G'$, we
write the difference
\begin{equation}
\label{eq:delta-r-delta-f0}
  \begin{split}
    K^R_\tau - K^R &= \tilde{L}[S_\tau -
  \psi_\tau\otimes\tilde{\ee}\phi_\tau - S_A +
  G\otimes\tilde{\ee}(G_N-\Delta_N G')] \\
  &= \tilde{L}[S_\tau - S_A + \Delta_N G\otimes G] -
  \tilde{L}[\psi_\tau\otimes \tilde{\ee}\phi_\tau - G\otimes
  \tilde{\ee}G_N] = \delta^{R} + \delta^F_0.
  \end{split}
\end{equation}
Set
\begin{equation*}
    S_{A_N}(s,t) = \int_0^\infty G(s+z)G_N(t+z) + G_N(s+z) G(t+z)dz,
\end{equation*}
since
$
  \int_0^\infty \Ai(s+z)\Ai'(t+z) + \Ai'(s+z)\Ai(t+z)dz = \int_0^\infty
  \frac{d}{dz}\left[\Ai(s+z)\Ai(t+z)\right]dz = -\Ai(s)\Ai(t),
$
we obtain
\begin{equation}
\label{eq:delta-r}
  \delta^R = \tilde{L}[S_\tau - S_{A_N}].
\end{equation}

Finally, we organize the components of $K_\tau - K_{GOE}$ as
\begin{equation}
  \label{eq:kern-diff}
  K_\tau - K_{GOE} = \delta^R + \delta^F_0 + \delta^F_1 +
    \delta^F_2
\end{equation}
where except for $\delta^F_0$ and $\delta^R$ given in
\eqref{eq:delta-r-delta-f0} and \eqref{eq:delta-r}, we further
define $\delta^F_i = K^F_{\tau,i} - K^F_i$ for $i = 1,2$.

By the bounds \eqref{eq:diff-det} and \eqref{eq:det-bd}, we need
entrywise bounds on $K_\tau - K_{GOE}$ to get our final convergence
rate. By the decomposition in \eqref{eq:kern-diff}, the problem
reduces to entrywise bounds for each of the $\delta$-terms. Since
all these entries have explicit representations, this becomes an
analysis problem which is to be solved in the next two sections.

\section{Proof} \label{sec:proof}

In this section, we prove Theorem \ref{thm:rmt} [and hence Theorem
\ref{thm:main}] by focusing on the entries of the $\delta$-terms in
\eqref{eq:kern-diff}. Besides the RMT analysis performed in Section
\ref{sec:prelim}, the proof needs two additional toolkits: a)
asymptotics of transformed Laguerre polynomials, and b) several
operator theoretic bounds of Hilbert-Schmidt and trace class norms.

\subsection{Preliminaries}
\label{subsec:preliminaries}

Here, we introduce some basic results for later repeated use in the
proof. Moreover, we make a specific choice of the weight function
$\rho$.

We start with Laguerre polynomial asymptotics. Recall that with
constants $\tmu_{n,N}, \tsigma_{n,N}$ in \eqref{eq:tmu-tsigma} and
functions $\phi, \psi$ defined in \eqref{eq:phi-psi}, we have
defined transformed Laguerre polynomials $\phi_\tau$ and $\psi_\tau$
in \eqref{eq:phi-psi-tau}. Moreover, for the Airy function, we
define
\begin{equation}
  \label{eq:G}
  G(s) = \frac{1}{\sqrt{2}}\,\Ai(s).
\end{equation}

By \eqref{eq:kern-decomp}, \eqref{eq:kern-comp} and
\eqref{eq:kern-diff}, the kernels $K_\tau$ and $K_{GOE}$ and hence
their difference are essentially expressed in terms of $\phi_\tau,
\psi_\tau, G$ and their variants. Therefore, we will find the
following set of asymptotic bounds helpful to the analysis of their
behavior.
\begin{lemma}
  \label{lemma:laguerre}
  Let $\phi_\tau$, $\psi_\tau$ and $G$ be defined as in
  \eqref{eq:phi-psi-tau} and \eqref{eq:G} and $\Delta_N$ to be
  defined in \eqref{eq:Delta-N}. If $n>N$, $N\to \infty$, $n = n(N)\to \infty$
  and $n/N\to \gamma\in [1,\infty)$, there exists a continuous and
  nonincreasing function $C(\cdot)$, such that for any real number
  $s_0$, there exists an integer $N_0(s_0,\gamma)$ for which we have
  that for all $s\geq s_0$ and $N\geq N_0(s_0,\gamma)$,
\begin{align}
  \labs \psi_\tau(s)\rabs, \labs \psi_\tau'(s) \rabs & \leq
  C(s_0)\exp(-s); \\
  \labs \phi_\tau(s)\rabs, \labs \phi_\tau'(s) \rabs & \leq
  C(s_0)\exp(-s); \\
  \labs \psi_\tau(s) - G(s) \rabs, \labs \psi_\tau'(s) - G'(s) \rabs &
  \leq C(s_0)N^{-2/3}\exp(-s); \label{eq:psi'-bd} \\
  \labs \psi_\tau(s) - G(s) - \Delta_N G'(s)\rabs, \labs \psi_\tau'(s)
  - G'(s) - \Delta_N G''(s) \rabs &
  \leq C(s_0)N^{-2/3}\exp(-s). \label{eq:phi'-bd}
\end{align}
\end{lemma}

In order not to distract us from the cause of proving Theorem
\ref{thm:rmt}, we defer the proof of Lemma \ref{lemma:laguerre} to
Section \ref{sec:lagu-polyn-asympt}. For the rest of Section
\ref{sec:proof}, let us assume temporarily that Lemma
\ref{lemma:laguerre} is already established.

In addition to the Laguerre polynomial asymptotics, we need some
operator theoretic bounds of Hilbert-Schmidt and trace class norms.
This set of tools has been previously established in \citet[Section
8.4.1]{johnstone07}. For the sake of completeness, we state them
here with some corrections and modifications that are helpful to our
context.

From now on, we fix a real number $s_0$ and consider any $s'\in
[s_0,\infty)$. In general, let an operator $T: L^2([s',\infty),
\rho_1) \to L^2([s',\infty), \rho_2)$ defined by
\begin{equation}\label{eq:operator-def}
f \mapsto Tf: (Tf)(u) = \int_{s'}^\infty T(u,v)f(v)dv
\end{equation}
for some kernel $T(u,v)$. We obtain that the Hilbert-Schmidt norm
$\|T\|_2$ of $T$ satisfies
\begin{equation*}
\|T\|^2_2  = \iint_{[s',\infty)^2} \labs T(u,v)\rabs^2
\rho_1^{-1}(v)\rho_2(u)dudv.
\end{equation*}

Following the notation in \citet{johnstone07}, we introduce the symbol
$\diamond$ for the following convolution type operator:
\begin{equation*}
(a\diamond b)(u,v) \equiv \int_0^\infty a(u+z)b(v+z)dz.
\end{equation*}
Among all the operators defined by \eqref{eq:operator-def}, we are
interested in those with kernels $D$ of the form $
D(u,v)=\aa(u)\beta(v)$, or $D(u,v) = \aa(u)\beta(v)(a\diamond
b)(u,v)$. We use the following notation for a Laplace-type transform:
\begin{equation*}
\mathcal{L}(\rho)[t] \equiv \int_{s'}^\infty e^{-tz}\rho(z)dz.
\end{equation*}

For an operator with kernel of the form $D(u,v) =
\aa(u)\beta(v)(a\diamond b)(u,v)$, we have the following bound on its
operator norm:
\begin{lemma}\label{lemma:hs-tr-bd}
Let $D$ be an operator taking $L^2([s',\infty), \rho_1)$ to
$L^2([s',\infty), \rho_2)$ and having kernel
$D(u,v) = \aa(u)\beta(v)(a\diamond b)(u,v)$,
where we assume, for $u\geq s'$, that
\begin{equation}
  \label{eq:hs-tr-bd1}
  |\aa(u)|\leq \aa_0e^{\aa_1 u}, \quad |\beta(u)|\leq \beta_0e^{\beta_1 u},\quad
  |a(u)|\leq a_0e^{-a_1 u}, \quad |b(u)|\leq b_0e^{-b_1 u}.
\end{equation}
If both $\mathcal{L}(\rho_1^{-1})$ and $\mathcal{L}(\rho_2)$
converge for $t> c$, and $2(a_1-\alpha_1), 2(b_1-\beta_1)
> c$, the Hilbert-Schmidt norm satisfies
\begin{equation*}
\|D\|_2 \leq
\frac{\aa_0\beta_0a_0b_0}{a_1+b_1}\Bigl\{\mathcal{L}(\rho_2)[2(a_1-\alpha_1)]
\mathcal{L}(\rho_1^{-1})[2(b_1-\beta_1)]\Bigr\}^{1/2}.
\end{equation*}
If $\rho_1=\rho_2,$ then the trace norm $\|D\|_1$ satisfies the same
bound.
\end{lemma}

Next, we investigate rank one operators with kernels of the form
$D(u,v) = \aa(u)\beta(v)$. First, a remark taken verbatim from
\citet{tw05}: the norm of an operator $D = \aa\otimes \beta$ taking
$L^2(\rho_1)$ to $L^2(\rho_2)$ with kernel $D(u,v)=\aa(u)\beta(v)$
is given by $\|D\| = \|\aa\|_{2,\rho_2}\|\beta\|_{2,\rho_1^{-1}}$.
Here, the norm can be trace class (if $\rho_1=\rho_2$) or
Hilbert-Schmidt since they agree for rank one operators. Moreover,
if $\aa$ and $\beta$ satisfies the bound \eqref{eq:hs-tr-bd1},
similar derivation to that for proving Lemma \ref{lemma:hs-tr-bd}
will give the following lemma specific for rank one operators.
\begin{lemma}
\label{lemma:rankone}
 Let $D=\aa\otimes\beta$ be a rank one operator taking
$L^2([s',\infty), \rho_1)$ to $L^2([s',\infty), \rho_2)$ and having
kernel $D(u,v) = \aa(u)\beta(v)$,
where we assume, for $u\geq s'$, that
$|\aa(u)| \leq \aa_0e^{\aa_1 u}$ and $|\beta(u)|\leq
\beta_0e^{\beta_1 u}$.
If both $\mathcal{L}(\rho_1^{-1})$ and $\mathcal{L}(\rho_2)$
converge for $t> c$ that $-2\alpha_1, -2\beta_1
> c$, the Hilbert-Schmidt norm satisfies
\begin{equation*}
  \|D\|_2 \leq
  {\aa_0\beta_0}\Bigl\{\mathcal{L}(\rho_2)[-2\alpha_1]\mathcal{L} (\rho_1^{-1})
  [-2\beta_1]\Bigr\}^{1/2}.
\end{equation*}
If $\rho_1=\rho_2,$ then the trace norm $\|D\|_1$ satisfies the same
bound.
\end{lemma}

\subsubsection{Choice of the weight function $\rho$}
\label{sec:rho-choice}
In order to make our arguments explicit and to obtain the
exponential decay of the right tail in our bound, we feel it
convenient to make a specific choice of the weight function $\rho$.

In particular, for $\nu \in (0,1]$ and to be specified later in
\eqref{eq:nu}, on the $s$-scale, let
\begin{equation}
  \label{eq:rho}
  \rho\circ \tau(s) = 1 + \exp\left(\nu|s|\right).
\end{equation}
The above definition implies that on the $x$-scale, we specify the
weight function as
\begin{equation*}
  \rho(x) = 1 + \exp\left(\frac{\nu\labs x -
  \tmu_{n,N}\rabs}{\tsigma_{n,N}}\right).
\end{equation*}
We remark that on the $x$-scale, our choice of $\rho$ depends on
$N$.

First of all, we check that our choice of $\rho$ [on the $x$-scale]
satisfies the two required conditions spelled out in Section
\ref{sec:det}. Condition (1) holds for $\rho^{-1}(x) \asymp
\exp(-\nu x/\tsigma_{n,N})$ as $x\to \infty$. Condition (2) holds if
$\phi, \psi, \phi', \psi', \tee\phi$ and $\tee\psi$ belong to
$L^2([x',\infty);\rho)\cap L^2([x',\infty);\rho^{-1})$. We take
$\phi$ and $\psi$ as examples, while the argument for the rest is
essentially the same. By the definition of $\phi_k$ in
\eqref{eq:phi-k}, the right tails of both $\phi$ and $\psi$ are
bounded by $\exp(-x/3)$. On the other hand, as $x\to \infty$,
$\rho^{\pm 1}(x) \asymp \exp(\pm \nu x/\tsigma_{n,N})$ with
$\nu/\tsigma_{n,N} \leq 1/\tsigma _{n,N} = \Oh{N^{-2/3}}$. These two
facts suffice to show that both $|\phi|^2\rho^{\pm 1}$ and
$|\psi|^2\rho^{\pm 1}$ are integrable over the region $[x',\infty)$,
at least when $N$ is large. Condition (2) is hence satisfied.

By \eqref{eq:rho}, the operator class $\mathcal{A}$ in
\eqref{eq:op-class} is now concrete. We now make valid all the
formal derivation in Section \ref{sec:prelim} by verifying that
$K_\tau, K_{GOE}\in \mathcal{A}$. Observing that $\tau$ is linear,
by \citet[Theorem VI.22(h) and Theorem VI.23]{reed&simon80},
condition (2) on $\rho$ implies that $K_\tau - K^\ee \in
\mathcal{A}$. The super exponential decay \eqref{eq:airy-asym} of
the Airy functions, together with the same theorems as above,
guarantees that $K_{GOE} - K^\ee\in \mathcal{A}$. Hence, we need
only to verify that $\ee: L^2([s',\infty);\rho\circ\tau) \to
L^2([s',\infty);\rho^{-1}\circ\tau)$ is Hilbert-Schmidt, which is an
immediate consequence of condition (1) on $\rho$.

From now on, we use $\rho$ to denote $\rho\circ\tau$ in
\eqref{eq:rho} with no ambiguity, for all the remaining discussion
in this paper focuses on the $s$-scale.

For the operator-theoretic bounds, by our choice of $\rho$ in
\eqref{eq:rho}, we could adapt Lemma \ref{lemma:hs-tr-bd} and Lemma
\ref{lemma:rankone} into a more convenient form as follows.
\begin{corollary}
  \label{cor:bd}
With $\rho$ as specified in \eqref{eq:rho}, for $\nu \leq \eta/2$, we have
\begin{equation}
  \label{eq:L-bd}
  \mathcal{L}\left(\rho^{\pm 1}\right)[\eta] \leq \frac{4}{\eta -
    \nu}\exp\left(-\eta s'\pm \nu |s'| \right) \leq
  \frac{8}{\eta}\exp\left(-\eta s'\pm \nu |s'| \right).
\end{equation}

In particular, under the assumption of Lemma \ref{lemma:hs-tr-bd}, if
$\{\rho_1,\rho_2\}\subset \{\rho, \rho^{-1}\}$ and $a_1 - \aa_1,
b_1 - \beta_1 \geq \nu$, then
\begin{equation}
\label{eq:D-bd-1}
  \|D\|_2, \|D\|_1 \leq C\ \frac{\aa_0 \beta_0 a_0
    b_0}{a_1+b_1}\ \exp\left[-(a_1+b_1-\aa_1-\beta_1)s' + \nu|s'|
    \right],
\end{equation}
where $C = C(a_1,\aa_1, b_1, \beta_1)$.

Under the assumptions of Lemma \ref{lemma:rankone}, if
$\{\rho_1,\rho_2\}\subset \{\rho, \rho^{-1}\}$ and $-\aa_1,
-\beta_1 \geq \nu$, then
\begin{equation}
\label{eq:D-bd-2}
  \|D\|_2, \|D\|_1 \leq C\aa_0\beta_0\exp\left[(\aa_1+\beta_1)s' +
    \nu|s'|\right],
\end{equation}
where $C = C(\aa_1, \beta_1)$.
\end{corollary}

The proof of \eqref{eq:L-bd} follows directly from the derivation in
\citet[p.50]{johnstone07}; see, in particular, Eq.(205), (206)
there. Then the operator bounds \eqref{eq:D-bd-1} and
\eqref{eq:D-bd-2} are obtained by plugging \eqref{eq:L-bd} into the
bounds in Lemma \ref{lemma:hs-tr-bd} and Lemma \ref{lemma:rankone}.

\subsection{Operator convergence}
\label{subsec:operator-convergence}

With the tools from the previous subsection, we work out here
entrywise bounds for each $\delta$ term given in the decomposition
\eqref{eq:kern-diff}.

\paragraph{$\delta^R$ term.} Using the $\diamond$ operator, we have
$\delta^R = \tilde{L}[S_\tau - S_{A_N}]$ with $S_\tau =
\phi_\tau\diamond \psi_\tau + \psi_\tau\diamond \phi_\tau$ and
$S_{A_N} = G_N\diamond G + G \diamond G_N$. We shall use the
abbreviation $D^{(k)}f$, $k = -1, 0$ and $1$ to denote $\tee f$, $f$
and $f'$ respectively. Regardless of the signs, we have the
following unified expression for the entries of $\delta^R$:
\begin{equation}
\label{eq:delta-R-entries}
  \begin{split}
    \delta^R_{ij} \,=\, & D^{(k)}(\phi_\tau - G_N)\diamond
    D^{(l)}\psi_\tau + D^{(k)}G_N\diamond D^{(l)}(\psi_\tau - G)\\
    & + D^{(k)}(\psi_\tau - G)\diamond D^{(l)}\phi_\tau +
    D^{(k)}G\diamond D^{(l)}(\phi_\tau - G_N),
  \end{split}
\end{equation}
for $i,j\in \{1, 2\}$, $k\in \{-1, 0\}$ and $l\in \{0, 1\}$. By
Lemma \ref{lemma:laguerre} and asymptotics of the Airy function [see
\eqref{eq:airy-asym}], we find that for any of the four terms in
\eqref{eq:delta-R-entries}, the condition \eqref{eq:hs-tr-bd1} is
satisfied with $\aa_0 = \beta_0 = 1$, $\aa_1 = \beta_1 = 0$, $a_1 =
b_1 = 1$ and $a_0, b_0$ as shown in the following table.
\begin{table}[!h]
  \centering
  \begin{tabular}{|l|r|r|}
\hline
    & $a_0$  & $b_0$  \\
\hline
$D^{(k)}(\phi_\tau - G_N)\diamond D^{(l)}\psi_\tau$ &
$C(s_0)N^{-2/3}$
& $C(s_0)$ \\
\hline $D^{(k)}G_N\diamond D^{(l)}(\psi_\tau - G)$ & $C(s_0)$ &
$C(s_0)N^{-2/3}$  \\
\hline
$D^{(k)}(\psi_\tau - G)\diamond D^{(l)}\phi_\tau$ & $C(s_0)N^{-2/3}$
 & $C(s_0)$  \\
\hline $D^{(k)}G\diamond D^{(l)}(\phi_\tau - G_N)$ & $C(s_0)$  &
$C(s_0)N^{-2/3}$  \\
\hline
  \end{tabular}
\end{table}

We apply Corollary \ref{cor:bd} and obtain that for $\nu \leq 1$,
\begin{equation}
  \label{eq:delta-R-bd}
  \|\delta^R_{ij}\| \,\leq\, C(s_0)N^{-2/3}\exp\left( -2s' + \nu|s'|
  \right).
\end{equation}
Here and after, the unspecified norm $\|\cdot\|$ denotes
Hilbert-Schmidt norm $\|\cdot\|_2$ if $i\neq j$ and trace class norm
$\|\cdot\|_1$ otherwise. We remark that by a simple triangular
inequality, we could choose the $C(s_0)$ function in the last
display as the sum of products of continuous and non-increasing
functions, which could be seen from the term
$(\aa_0\beta_0a_0b_0)/(a_1+b_1)$ in \eqref{eq:D-bd-1}. Moreover, the
term $C$ in \eqref{eq:D-bd-1} is a universal constant for fixed
$a_1, \aa_1, b_1$ and $\beta_1$ here. Hence, the final $C(s_0)$
function remains continuous and non-increasing. For the other
$\delta$ terms, we will have the same result by the same arguments
and hence will be omitted.

\paragraph{$\delta^F_0$ term.} We reorganize $\delta^F_0$ as
\begin{equation*}
  \delta^F_0 = -\tilde{L}[\psi_\tau\otimes \tee \phi_\tau - G\otimes
  \tee G_N] = -\tilde{L}[\psi_\tau \otimes \tee(\phi_\tau - G_N) +
  (\psi_\tau - G)\otimes \tee G_N] = \delta^{F,1}_0 + \delta^{F,2}_0.
\end{equation*}
The entries of $\delta^{F,i}_0$, $i=1, 2$ are all of the form
$\aa(s)\beta(t)$ with the multipliers chosen from
$D^{(k)}\psi_\tau$, $D^{(k)}(\phi_\tau - G_N)$, $D^{(k)}(\psi_\tau -
G)$ and $D^{(k)}G_N$ for $k\in \{-1, 0, 1\}$. For these multipliers,
the condition for Lemma \ref{lemma:rankone} holds with the constants
$\aa_1 = \beta_1 = -1$ and $\aa_0$ (or $\beta_0$) specified below.
\begin{table}[!htb]
  \centering
  \begin{tabular}{|l|r|}
    \hline
      & $\aa_0$ (or $\beta_0$)  \\
    \hline
      $D^{(k)}\psi_\tau$ & $C(s_0)$  \\
      \hline
      $D^{(k)}(\phi_\tau - G_N)$ & $C(s_0)N^{-2/3}$  \\
      \hline
      $D^{(k)}(\psi_\tau - G)$ & $C(s_0)N^{-2/3}$  \\
      \hline
      $D^{(k)}G_N$ & $C(s_0)$  \\
      \hline
  \end{tabular}
\end{table}

We apply Corollary \ref{cor:bd} for these rank one terms and obtain
that for $\nu \leq 1$,
\begin{equation}
  \label{eq:delta-0-bd}
  \|\delta^F_{0,ij}\| \leq \|\delta^{F,1}_{0,ij}\| +
  \|\delta^{F,2}_{0,ij}\| \leq C(s_0)N^{-2/3}\exp\left(-2s' +
    \nu|s'|\right).
\end{equation}

\paragraph{$\delta^F_1$ and $\delta^F_2$ terms.} For these two terms, we
have
\begin{equation*}
  \delta^F_1 = L_1\left[\psi_\tau\otimes \beta_N - G \otimes
  \tfrac{1}{\sqrt{2}}\right] \quad \text{and} \quad
\delta^F_2 = L_2\left[\beta_N\otimes\psi_\tau -
  \tfrac{1}{\sqrt{2}}\otimes G\right].
\end{equation*}
By their similarity, we take $\delta^F_1$ as example and the same
analysis applies to $\delta^F_2$ with obvious modification. For
$\delta^F_1$, we reorganize it as
\begin{equation*}
  \delta^F_1 = L_1\left[ (\psi_\tau - G)\otimes \beta_N + G\otimes
    \left(\beta_N - \tfrac{1}{\sqrt{2}}\right) \right] =
  \delta^{F,1}_1 + \delta^{F,2}_1.
\end{equation*}
For analysis of the terms here, Corollary \ref{cor:bd} no longer
works and we give an alternative bound which was derived in full
detail in \citet{johnstone07}. In particular, consider matrices of
rank one operators on $L^2([s',\infty); \rho)\otimes
L^2([s',\infty); \rho^{-1})$, we denote, here and after, the
$L^2$-norm on $L^2([s',\infty);\rho)$ and
$L^2([s',\infty);\rho^{-1})$ by $\|\cdot\|_+$ and $\|\cdot\|_-$
respectively. \citet[Eq.(214)]{johnstone07} gives the following
bound
\begin{equation*}
  \begin{pmatrix}
    \|a_{11}\otimes b_{11} \|_1 & \|a_{12}\otimes b_{12} \|_2 \\
    \|a_{21}\otimes b_{21} \|_2 & \|a_{22}\otimes b_{22} \|_1
  \end{pmatrix}
  \leq
  \begin{pmatrix}
    \|a_{11}\|_+\|b_{11}\|_- & \|a_{12}\|_+\|b_{12}\|_+ \\
    \|a_{21}\|_-\|b_{21}\|_- & \|a_{22}\|_-\|b_{22}\|_+
  \end{pmatrix}.
\end{equation*}

By the inequality above and our reorganization of $\delta^F_1$, we
will see that the essential elements we need to bound are
$\|D^{(k)}(\psi_\tau - G)\|_\pm$, $\|D^{(k)}G\|_\pm$ and $\|1\|_-$ for
$k = -1$ and $0$.

For $\|D^{(k)}(\psi_\tau - G)\|_\pm$, we obtain from Lemma
\ref{lemma:laguerre} and \eqref{eq:L-bd}
that for $\nu \leq 1$:
\begin{equation*}
  \|D^{(k)}(\psi_\tau - G)\|_\pm^2 \leq
  C^2(s_0)N^{-4/3}\mathcal{L}(\rho^{\pm 1})[2] \leq
  C^2(s_0)N^{-4/3}\exp(-2s'+\nu|s'|).
\end{equation*}
For $\|D^{(k)}G\|_\pm$, asymptotics of the Airy function and
\eqref{eq:L-bd} give that for $\nu\leq 1$:
\begin{equation*}
  \|D^{(k)}G\|_\pm^2 \leq C^2(s_0)\mathcal{L}(\rho^{\pm 1})[2] \leq
  C^2(s_0)\exp(-2s' + \nu|s'|).
\end{equation*}
Finally, for $\|1\|_-$, we derive directly that
\begin{equation*}
  \begin{split}
    \|1\|_-^2 & = \int_{s'}^\infty \left[1 +
      \exp(\nu|s|)\right]^{-1}ds \leq \int_{s'}^\infty \exp(-\nu|s|)ds
    \\
    & \leq \int_0^\infty \exp(-\nu s)ds + \int_{-|s'|}^0 \exp(\nu s)ds
    = \frac{1}{\nu} + \frac{1}{\nu} - \frac{1}{\nu}\exp(-\nu |s'|)
    \leq \frac{2}{\nu}.
  \end{split}
\end{equation*}

By definition of the operator $L_1$ and our reorganization, we have
the first column of $\delta^F_1$ as following while the second
column of it are zeros:
\begin{equation*}
  \begin{pmatrix}
    \delta^F_{1,11} \\ \delta^F_{1,21}
  \end{pmatrix} =
  \begin{pmatrix}
    (\psi_\tau - G)\otimes \beta_N + G\otimes (\beta_N - 1/\sqrt{2})\\
    -\tee (\psi_\tau - G)\otimes \beta_N - \tee G\otimes (\beta_N -
    1/\sqrt{2})
  \end{pmatrix}.
\end{equation*}
Assuming $\beta_N - 1/\sqrt{2} = \Oh{N^{-1}}$ [for a proof, see
\ref{subsec:a-constants}], we have
\begin{equation*}
  \begin{split}
    \|\delta^F_{1,11}\|_1 \leq &\|(\psi_\tau - G)\otimes \beta_N\|_1
    +\| G\otimes(\beta_N - 1/\sqrt{2})\|_1 \\
  \leq & \|(\psi_\tau - G)\|_+ \|\beta_N\|_- + \|G\|_+ \|\beta_N -
  1/\sqrt{N}\|_- \\
  \leq & C(s_0)N^{-2/3}\nu^{-1/2}\exp\left(-s' + \nu |s'|/2\right) +
  C(s_0)N^{-1}\nu^{-1/2}\exp\left(-s' + \nu |s'|/2\right) \\
  \leq & C(s_0)N^{-2/3}\nu^{-1/2}\exp\left(-s' + \nu |s'|/2\right)
  \leq C(s_0)N^{-2/3}\exp(-s'/2).
  \end{split}
\end{equation*}
The last inequality holds by fixing $\nu$, for example, at $1$. By
the same calculation, this bound also holds for
$\|\delta^{F}_{1,12}\|_2$ and those entries of $\delta^F_2$.
Finally, we conclude our analysis with the following bound on
entries of $\delta^F_1$ and $\delta^F_2$: for $\nu = 1$, we have
\begin{equation}
  \label{eq:delta-1-2-bd}
  \|\delta^F_{1,ij}\|, \|\delta^F_{2,ij}\| \leq
  C(s_0)N^{-2/3}\exp\left(-s'/2\right).
\end{equation}

\subsection{Proof of Theorem \ref{thm:rmt}}
\label{subsec:proof}

Throughout the proof, we fix
\begin{equation}\label{eq:nu}
    \nu = 1
\end{equation}
in the weight function $\rho$ specified in \eqref{eq:rho}.

By \eqref{eq:kern-diff} and the bounds \eqref{eq:delta-R-bd},
\eqref{eq:delta-0-bd} and \eqref{eq:delta-1-2-bd}, we bound the
entries of $K_\tau - K_{GOE}$ using a simple triangular inequality
\begin{equation*}
    \|K_{\tau,ij} -K_{GOE,ij} \|\leq C(s_0)N^{-2/3}\exp(-s'/2).
\end{equation*}
Apply Proposition \ref{prop:det-bd} with $A = K_\tau$ and $B =
K_{GOE}$,
\begin{equation}\label{eq:det-diff-bd}
    \labs \det(I-K_\tau) - \det(I-K_{GOE}) \rabs \leq
    M(K_{GOE})C(s_0)N^{-2/3}\exp(-s'/2),
\end{equation}
where
\begin{equation*}
    M(K_{GOE}) = 2\,\det(I-K_{GOE}) + 2\exp\left[2\left( 1 + \|K_{GOE}\|_2\right)^2 + \sum_i
    \|K_{GOE,ii}\|_1\right] .
\end{equation*}

For the first term in $M(K_{GOE})$, we have $\det(I-K_{GOE}) =
F_1^2(s') \leq 1$. On the other hand, we have
\begin{equation*}
    \|K_{GOE}\|_2 \leq \sum_{i,j}\|K_{GOE,ij}\|_2 \leq
    \sum_i\|K_{GOE,ii}\|_1 + \sum_{i\neq j}\|K_{GOE,ij}\|_2.
\end{equation*}
In principle, one could show for each $i$ and $j$
\begin{equation*}
    \|K_{GOE,ij}\| \leq C(s_0),
\end{equation*}
with $C(s_0)$ continuous and non-increasing. Here, we only take
$\|K_{GOE,11}\|_1$ as an example for the proof of the others is
essentially the same. Let $H_\tau$ and $G_\tau$ be Hilbert-Schmidt
operators with kernels $\phi_\tau(x+y)$ and $\psi_\tau(x+y)$
respectively, then as operator
\begin{equation*}
    K_{GOE,11} = H_\tau G_\tau + G_\tau H_\tau + G\otimes
    \frac{1}{\sqrt{2}} - G\otimes \tee G.
\end{equation*}
By the relation $\|AB\|_1 \leq \|A\|_2\|B\|_2$,
\begin{equation*}
    \|K_{GOE,11}\|_1 \leq 2\,\|H_\tau\|_2\|G_\tau\|_2 +
    \frac{1}{\sqrt{2}}\|G\|_{2,\rho}\left\|1\right\|_{2,\rho^{-1}} +
    \|G\|_{2,\rho}\|\tee G\|_{2,\rho^{-1}}.
\end{equation*}
Each norm on the right hand side of the above inequality is the
square root of an integral of a positive function on $[s',\infty)$
or $[s',\infty)^2$ that is bounded by the corresponding integral
over $[s_0,\infty)$ or $[s_0,\infty)^2$, which in turn is continuous
and non-increasing in $s_0$. Hence, $\|K_{GOE,11}\|_1 \leq C(s_0)$.

By the above discussion, we could control the second term of
$M(K_{GOE})$ and hence itself by a continuous and non-increasing
$C(s_0)$. Finally, we complete the proof by combining this fact with
the initial bounds \eqref{eq:det-diff-bd} and \eqref{eq:diff-det}.

\section{Laguerre Polynomial Asymptotics}
\label{sec:lagu-polyn-asympt}

In this section, our goal is to establish Lemma
\ref{lemma:laguerre}. To this end, we exploit the Liouville-Green
approach to study the related asymptotics for Laguerre polynomials
of both large order and large degree. This approach has been
successfully used in \citet{johnstone01}, \citet{nek06} and more
recently, \citet{johnstone07} in deriving similar type of results.
The novelty here is the establishment of the bounds
\eqref{eq:psi'-bd} and \eqref{eq:phi'-bd} for the derivatives of
these polynomials.

To start with, let us consider the ``intermediate'' function
$F_{n,N}$ introduced in \citet[Section 2.2.2]{nek06} as
\begin{equation}\label{eq:fnN-def}
    F_{n,N}(x) \equiv (-1)^N \sigma_{n,N}^{-1/2} \sqrt{N!/n!}x^{(\aa_N +
    1)/2}e^{-x/2}L_N^{\aa_N}(x)
\end{equation}
with $\aa_N = n - N$. We could then relate $F_{n,N}$ to $\phi_N,
\phi$ and $\phi_\tau$ as
\begin{gather}
\phi_N(x;\aa_N) = (-1)^N\sigma_{n,N}^{1/2}x^{-1/2}F_{n,N}(x), \tag*{}\\
\phi(x;\aa_N - 1) =
\frac{N^{1/4}(n-1)^{1/4}}{\sqrt{2}}\sigma_{n-2,N}^{1/2}F_{n-2,N}(x)/x,\tag*{}\\
\text{and}\quad \phi_\tau(s) = \frac{1}{\sqrt{2}} \left(
\frac{N^{1/4}(n-1)^{1/4}\sigma_{n-2,N}^{1/2}\tsigma_{n,N}}{\mu_{n-2,N}}
\right) F_{n-2,N}(\tauexpr{s}) \left(
\frac{\mu_{n-2,N}}{\tauexpr{s}} \right), \tag*{}
\end{gather}
with $\mu_{n,N}$ and $\sigma_{n,N}$ defined as
\begin{equation*}
\mu_{n,N} = \left(\sqrt{n_+} + \sqrt{N_+}\right)^2\quad
\text{and}\quad \sigma_{n,N} = \left(\sqrt{n_+} +
\sqrt{N_+}\right)\left( \frac{1}{\sqrt{n_+}} + \frac{1}{\sqrt{N_+}}
\right)^{1/3},
\end{equation*}
using the abbreviations $n_+ = n + \tfrac{1}{2}$ and $N_+ = N +
\tfrac{1}{2}$. If we replace the subscripts $(n-2, N)$ in
$\mu_{n-2,N}, \sigma_{n-2,N}$ and $F_{n-2,N}$ by $(n-1, N-1)$ on the
right hand sides of the expressions for $\phi(x;\aa_N - 1)$ and
$\phi_\tau(s)$, we obtain the identities for $\psi(x;\aa_N-1)$ and
$\psi_\tau(s)$. Due to this close connection of $\phi_\tau$ and
$\psi_\tau$ to $F_{n,N}$, the essential element for proving the
desired asymptotic bounds reduces to the understanding of the
behavior of $F_{n,N}$ and its derivative, for which the
Liouville-Green approach is instrumental.

In the rest of this section, we first study in detail the
Liouville-Green approximation to the $F_{n,N}$ function and its
derivative. Then the result is used to facilitate the derivation of
the global bounds and the local as well as global Airy approximation
to $\phi_\tau$, $\psi_\tau$ and their derivatives.

\subsection{Liouville-Green approach}
\label{subsec:louville-green}

Many of the arguments in this part have been spelled out in some
detail in \citet{johnstone01} and \citet{nek06}. A more complete
account of the theory could be found in \citet[Chapter 11]{olver74}.
However, for completeness, we state them here briefly with notation
similar to that in \citet{nek06}.

Consider $w_N(x) = x^{(\aa_N + 1)/2}e^{-x/2}L_N^{\aa_N}(x)$ as a
multiple of $F_{n,N}$, we have
\begin{equation}\label{eq:diff-eqn}
\frac{d^2
w_N}{dx^2}=\left\{\frac{1}{4}-\frac{\kk_N}{x}+\frac{\lambda_N^2-1/4}
{x^2}\right\} w_N
\end{equation}
with $\kk_N = N + \tfrac{\aa_N + 1}{2} = \tfrac{n+N+1}{2}$ and
$\lambda_N = \tfrac{\aa_N}{2} = \tfrac{n-N}{2}$.

By a change of variable $\xi = x/\kk_N$, we obtain
\begin{equation*}
\frac{d^2w_N}{d\xi^2}=\left\{\kappa^2 f(\xi)+g(\xi)\right\}w_N,
\end{equation*}
where
\begin{equation*}
f(\xi)=\frac{(\xi-\xi_-)(\xi-\xi_+)}{4\xi^2} \quad\text{and}\quad
g(\xi)=\frac{1}{4\xi^2}\ ,
\end{equation*}
with $\xi_{\pm}=2\pm \sqrt{4-\ww_N^2}$ and
$\ww_N=2\lambda_N/\kk_N=\tfrac{2(n-N)}{n+N+1}$. The Liouville-Green method
introduces the change of independent variable as
\begin{equation*}
\frac{2}{3}\zz^{3/2} = \int_{\xi_+}^{\xi} \sqrt{f(t)}dt\quad
(\xi\geq \xi_+)\quad \text{and}\quad \frac{2}{3}(-\zz)^{3/2} =
\int_\xi^{\xi_+}\sqrt{-f(t)}dt\quad (\xi\leq \xi_+),
\end{equation*}
and defines a new dependent variable $W = (d\zz/d\xi)^{1/2}w_N$. For
the new pair $(W, \zz)$, we have the new differential equation as
\begin{equation*}
\frac{d^2W}{d\zeta^2}=\left\{\kk_N^2\zz+v(\ww_N, \zz)\right\}W.
\end{equation*}
Let $\hat{f} = f/\zz$, the recessive solution of \eqref{eq:diff-eqn}
satisfies \citep[p.399, Theorem 3.1]{olver74}
\begin{equation*}
  w_N(\kk_N\xi) \propto \hat{f}^{-1/4}(\xi)\{\Ai(\kk_N^{2/3}\zz) + \ee_2(\kk_N,
  \xi)\},
\end{equation*}
with the following estimates for the error term $\ee_2$ and its
derivative with $\xi\in [2,\infty)$:
\begin{align*}
\labs \ee_2(\kk_N, \xi) \rabs & \leq
\MM(\kk_N^{2/3}\zz)\EE^{-1}(\kk_N^{2/3}\zz)\left[ \exp\left(
\frac{\lambda_0}{\kk_N} F(\ww_N)\right) - 1 \right], \\
\labs \partial_\xi \ee_2(\kk_N, \xi) \rabs & \leq
\kk_N^{2/3}\hat{f}^{1/2}(\xi)\NN(\kk_N^{2/3}\zz)\EE^{-1}(\kk_N^{2/3}\zz)
\left[ \exp\left( \frac{\lambda_0}{\kk_N} F(\ww_N)\right) -
1\right].
\end{align*}
In the above bounds, $\MM, \EE$ are the modulus and weight functions
for the Airy function, and $\NN$ the phase function for its
derivative \citep[pp.394-396]{olver74}. Moreover, $\lambda_0 \doteq
1.04$ and $F(\ww_N)$ has been well studied in \citet[A.3]{nek06}.

For the function $F_{n,N}$ of our interest, we have from
\citet[Eq.(5) and A.1]{nek06} that
\begin{equation*}
F_{n,N}(x) = r_N\left( \frac{\kk_N}{\sigma_{n,N}^3}
\right)^{1/6}\hat{f}^{-1/4}(\xi)\{\Ai(\kk_N^{2/3}\zz) + \ee_2(\kk_N,
\xi)\},
\end{equation*}
with
\begin{equation}
\label{eq:r-N}
    r^2_N = \frac{2\pi \exp[-(n_+ + N_+)]n_+^{n_+}N_+^{N_+}}{N!n!} =
    1 + \Oh{n^{-1}, N^{-1}}.
\end{equation}

For the convenience of argument, we define an auxiliary function
$R_N(\xi) = (\dot{\zz}(\xi)/\dot{\zz}_N)^{-1/2}$ with $\dot{\zz}_N =
\dot{\zz}(\xi_+)$. We remark that by our definition, we have
$\sigma_{n,N} = (\kk_N^{-1/3}\dot{\zz}_N)^{-1}$ and $\hat{f} =
\dot{\zz}(\xi)^2$. Hence, $F_{n,N}$ could be rewritten as
\begin{equation}
\label{eq:fnN}
F_{n,N} = r_NR_N(\xi)\{\Ai(\kk_N^{2/3}\zz) +
\ee_2(\kk_N, \xi)\}.
\end{equation}

Finally, we conclude this part with some useful bounds and
asymptotics of $\MM, \EE, \NN$ and the Airy function
\citep[pp.392-397]{olver74}. As $x\to \infty$, we have
\begin{equation}
\label{eq:EMN}
    \EE(x)\sim \sqrt{2}e^{ \frac{2}{3}x^{3/2} },\quad
    \MM(x)\sim \pi^{-1/2}x^{-1/4},\quad \NN(x)\sim
    \pi^{-1/2}x^{1/4}.
\end{equation}
For all $x > 0$, the Airy function and its derivative are bounded as
\begin{equation}
\label{eq:airy-asym}
    0\leq \Ai(x)\leq
    \frac{e^{-\frac{2}{3}x^{3/2}}}{2\pi^{1/2}x^{1/4}},\quad \labs
    \Ai'(x) \rabs \leq \left(1 +
    \frac{7}{48x^{3/2}}\right)\frac{x^{1/4}e^{-\frac{2}{3}x^{3/2}}}{2\pi^{1/2}}.
\end{equation}
Finally, for all $x$, we have the following bounds
\begin{equation}
\label{eq:airy-bd}
    \labs \Ai(x)\rabs \leq \MM(x)\EE^{-1}(x),\quad \labs
    \Ai'(x)\rabs \leq \NN(x)\EE^{-1}(x),\quad \MM(x)\leq 1, \quad
    \EE(x)\geq 1,
\end{equation}
and finally, $\EE(x)$ is monotone increasing in $x$
\citep[p.395]{olver74}.

\subsection{Large $N$ asymptotics}
\label{subsec:asymptotics}

We now derive the large $N$ asymptotics of $\phi_\tau$, $\psi_\tau$
and related functions. First, we use the analysis done in
\citet{johnstone01} and \citet{nek06} to obtain bounds for
$|\psi_\tau|$ and $|\psi_\tau - G|$ without much extra effort. Then
we derive the bounds for $|\psi'_\tau|$ and $|\psi'_\tau - G|$,
which need some careful analysis to be detailed below and the bound
on $|\psi_\tau - G|$ is then further refined to match the claim in
Lemma \ref{lemma:laguerre}. Finally, corresponding results for
quantities related to $\phi_\tau$ could be obtained by understanding
the difference of the centering and scaling constants involved in
$\phi_\tau$ and $\psi_\tau$.

\subsubsection{Bounds for $|\psi_\tau(s)|$ and $|\psi_\tau(s) -
G(s)|$} \label{subsec:psi}

We define $x_{n,N}(s) = \xnexpr{s}$ and let
$$\theta_{n,N}(x_{n,N}(s)) \equiv
F_{n,N}(x_{n,N}(s))\left(\frac{\mu_{n,N}}{x_{n,N}(s)}\right).$$
\citet[A.8]{johnstone01} showed that under the condition of Lemma
\ref{lemma:laguerre},
\begin{equation*}
    \labs F_{n,N}(x_{n,N}(s))\sigma_{n,N}^{1/2}N^{-1/6} \rabs \leq
    C\exp(-s), \quad \text{for all } s\geq 0.
\end{equation*}
Simple manipulation gives $\sigma_{n,N}^{-1/2}N^{1/6}\leq
\sigma_{n,N}^{-1/2}N_+^{1/6}\leq \left( N_+/n_+ \right)^{1/2} < 1$,
and hence for all $s\geq 0$,
\begin{equation*}
\labs F_{n,N}(x_{n,N}(s)) \rabs \leq C\exp(-s).
\end{equation*}
If $s_0<0$, by \eqref{eq:r-N}, \eqref{eq:EMN},
\eqref{eq:R-N-minus-1-bd} and \citet[A.3]{nek06}, we obtain that
when $N\geq N_0(s_0,\gamma)$,
\begin{equation*}
\labs F_{n,N}(x_{n,N}(s))\rabs \leq r_N \labs R_N(\xi)
\rabs\MM(\kk_N^{2/3}\zz)\EE^{-1}(\kk_N^{2/3}\zz) \leq 2\
\EE^{-1}(\kk_N^{2/3}\zz) \leq 2.
\end{equation*}
If we let $M(s_0) = \max_{s\in [s_0,0]}\{2e^s\}$, and define
\begin{equation*}
C(s_0) = \max\{C, M(s_0)I_{s_0 < 0}, M(0)\},
\end{equation*}
it is then continuous and non-increasing in $s_0$ as desired and we
have that when $N\geq N_0(s_0,\gamma)$, $|F_{n,N}(x_{n,N}(s))|\leq
C(s_0)\exp(-s)$ for all $s\geq s_0$. Moreover, by noting
$\sigma_{n,N}/\mu_{n,N} = \Oh{N^{-2/3}}$, when $N$ is larger than
some constant that depends only on $s_0$,
\begin{equation*}
    \frac{\mu_{n,N}}{x_{n,N}(s)} \leq \left(1 +
    s_0 \frac{\sigma_{n,N}}{\mu_{n,N}}\right)^{-1} \leq 2, \quad
    \text{for all } s\geq s_0.
\end{equation*}
Hence, under the condition of Lemma \ref{lemma:laguerre}, we have
that when $N\geq N_0(s_0,\gamma)$,
\begin{equation*}
    \labs \theta_{n,N}(x_{n,N}(s)) \rabs \leq C(s_0)\exp(-s), \quad
    \text{for all } s\geq s_0.
\end{equation*}

Later on, \citet[Section 3.2]{nek06} showed that for any constant
$\varrho_N = 1 + \Oh{N^{-1}}$, if we define
$\Delta_{n,N}(x_{n,N}(s)) = |\varrho_N\theta_{n,N}(x_{n,N}(s)) -
\Ai(s)|$, then under the condition of Lemma \ref{lemma:laguerre}, we
have
\begin{equation*}
    N^{2/3}\Delta_{n,N}(x_{n,N}(s))\leq C(s_0)\exp(-s/2),\quad \text{for all } s\geq
    s_0.
\end{equation*}

For $\psi_\tau(s)$, we have $\tmu_{n,N} = \mu_{n-1,N-1}$ and
$\tsigma_{n,N} = \sigma_{n-1,N-1}$ and hence it is of the form
$\frac{1}{\sqrt{2}} \rho_{N} \theta_{n-1,N-1}(x_{n-1,N-1}(s))$.
Noting that $\rho_N = 1 + \Oh{N^{-1}}$ [see \ref{subsec:a-constants}
for a proof], we apply the bounds for $\theta_{n,N}$ and
$\Delta_{n,N}$ directly and obtain that under the condition of Lemma
\ref{lemma:laguerre}, when $N\geq N_0(s_0,\gamma)$,
\begin{equation*}
    |\psi_\tau(s)|\leq C(s_0)\exp(-s),\quad \text{and}\quad
    \labs \psi_\tau(s) - G(s) \rabs \leq
    C(s_0)N^{-2/3}\exp(-s/2),\quad \text{for all } s\geq s_0.
\end{equation*}
Actually the bound on $|\psi_\tau(s) - G(s)|$ could be further
improved to be that claimed in Lemma \ref{lemma:laguerre}: see
\eqref{eq:psitau-g-improved} for the refinement. We also remark that
we could not apply the results directly to $\phi_\tau$ since the
centering and scaling constants $(\mu_{n-2,N},\sigma_{n-2,N})$
specific to $F_{n-2,N}$ does not agree with the global constants
$(\tmu_{n,N}, \tsigma_{n,N})$ which we use.

\subsubsection{Bounds for $|\psi'_\tau(s)|$ and $|\psi'_\tau(s) -
G'(s)|$}

As we have seen, the analysis of $\psi_\tau$ depends on our
understanding of the function $\theta_{n,N}(x_{n,N}(s))$. To
investigate the bounds for $\psi'_\tau$ and its approximation by
$G'$, we start with a detailed analysis of the quantity
$\partial_s\theta_{n,N}(x_{n,N}(s))$.

We split $\partial_s\theta_{n,N}(x_{n,N}(s))$ into two parts:
\begin{equation*}
    \begin{split}
        \labs \partial_s \theta_{n,N}(x_{n,N}(s)) \rabs & \leq \labs
        \sigma_{n,N}F'_{n,N}(x_{n,N}(s))\frac{\mu_{n,N}}{x_{n,N}(s)}
        \rabs + \labs
        \sigma_{n,N}F_{n,N}(x_{n,N}(s))\frac{\mu_{n,N}}{x^2_{n,N}(s)}
        \rabs \\
        & = T_{N,1}(s) + T_{N,2}(s).
    \end{split}
\end{equation*}

\paragraph{$T_{N,2}$ term.} This term is relatively easy to bound.
Note that $T_{N,2}(s) =
|\theta_{n,N}(x_{n,N}(s))\sigma_{n,N}/x_{n,N}(s)|$ and that
$\sigma_{n,N}/\mu_{n,N} = \Oh{N^{-2/3}}$. When $N\geq N_0(s_0)$, the
ratio
\begin{equation*}
    |\sigma_{n,N}/x_{n,N}(s)| = |s + \mu_{n,N}/\sigma_{n,N}|^{-1} \leq
    C(s_0)N^{-2/3},\quad \text{for all } s\geq s_0.
\end{equation*}
Hence, by our previous bound on $|\theta_{n,N}|$, we obtain that
under the condition of Lemma \ref{lemma:laguerre},
\begin{equation*}
    T_{N,2}(s)\leq C(s_0)N^{-2/3}\exp(-s),\quad \text{for all }s\geq s_0.
\end{equation*}

\paragraph{$T_{N,1}$ term.}
Recalling that $\mu_{n,N}/x_{n,N}(s)$ could be bounded by $2$, we
focus on $\sigma_{n,N}F'_{n,N}$. Thinking of $x = x_{n,N}(s)$, we
have from \eqref{eq:fnN} that
\begin{equation*}
    \begin{split}
        \sigma_{n,N}F'_{n,N}(x) \ =\  & r_N
        \left(\frac{\sigma_{n,N}}{\kk_N}\right)R'_N(\xi)\left[\Ai(\kk_N^{2/3}\zz)
        + \ee_2(\kk_N, \xi)\right]\\
    & + r_N R_N(\xi)\left[ \Ai'(\kk_N^{2/3}\zz)R_N^{-2}(\xi) + \left(\frac{\sigma_{n,N}}
    {\kk_N}\right) \frac{\partial}{\partial \xi} \ee_2(\kk_N, \xi)
    \right].
    \end{split}
\end{equation*}
To facilitate our analysis, on the $s$-scale, we divide the whole
region $[s_0,\infty)$ as $I_{1,N}\cup I_{2,N}$ with $I_{1,N} = [s_0,
s_1N^{1/6})$ and $I_{2,N} = [s_1N^{1/6}, \infty)$. The choice of
$s_1$ is made explicit in \ref{subsec:a-s-1}.

\emph{Case $s\in I_{1,N}$.} In this case, we first reorganize
$\sigma_{n,N}F'_{n,N}(x)$ as $\sigma_{n,N}F'_{n,N}(x) = \sum_{i=1}^4
D^i_{n,N}$, with
\begin{align*}
D^1_{n,N} & = r_N
\left(\frac{\sigma_{n,N}}{\kk_N}\right)R'_N(\xi)\{\Ai(\kk_N^{2/3}\xi)
+ \ee_2(\kk_N, \xi)\}, &
D^2_{n,N} & = r_N[R_N^{-1}(\xi) - 1]\Ai'(\kk_N^{2/3}\zz),\\
D^3_{n,N} & = r_N \Ai'(\kk_N^{2/3}\zz),& D^4_{n,N} & = r_N
\left(\frac{\sigma_{n,N}}{\kk_N}\right) R_N(\xi)
\frac{\partial}{\partial \xi} \ee_2(\kk_N, \xi).
\end{align*}

\noindent \textit{(i)} Consider $D^1_{n,N}$ first. Direct
computation shows $N^{2/3}(\sigma_{n,N}/\kk_N) \to
2(1+1/\sqrt{\gamma})^{1/3}(1+\sqrt{\gamma})(1+\gamma)^{-1}$. Hence
when $N\geq N_0(\gamma)$, we have the bound
\begin{equation*}
    N^{2/3}\left(\frac{\sigma_{n,N}}{\kk_N}\right) \leq
    C\left(1+\frac{1}{\sqrt{\gamma}}\right)^{1/3}\frac{1+\sqrt{\gamma}}{1+\gamma}.
\end{equation*}
Moreover, by the bound \eqref{eq:R-N'-bd} for $R_N'$ and recalling
that $\gamma \geq 1$, we know that when $N\geq N_0(s_0, \gamma)$,
\begin{equation}
\label{eq:D-n-N-1-1}
  N^{2/3}\left(\frac{\sigma_{n,N}}{\kk_N}\right) \labs R'_N(\xi)\rabs
  \leq C\left(1+\frac{1}{\sqrt{\gamma}}\right)^{1/3}\frac{\sqrt{\gamma}}{1+\sqrt{\gamma}}\leq
  C,\quad \text{for all $s\in I_{1,N}$}.
\end{equation}
On the other hand, by \eqref{eq:EMN} and \eqref{eq:kk-N-zz-minus-s},
we obtain
\begin{equation*}
    \labs \Ai(\kk_N^{2/3}\zz) + \ee_2(\kk_N, \xi) \rabs \leq C\,
    \MM(\kk_N^{2/3}\zz)\EE^{-1}(\kk_N^{2/3}\zz)\leq C\,
    \EE^{-1}(\kk_N^{2/3}\zz).
\end{equation*}
When $s\geq 0$, we know from \eqref{eq:kk-N-zz-minus-s} and the
monotonicity of $\EE$ that when $N\geq N_0(s_0,\gamma)$,
$\kk_N^{2/3}\zz \geq s/2$ holds, and hence by
\eqref{eq:EMN},
\begin{equation*}
  \labs \Ai(\kk_N^{2/3}\zz) + \ee_2(\kk_N, \xi) \rabs \leq C\,
  \EE^{-1}(s/2)\leq
  C\,\exp\left(-\frac{1}{3\sqrt{2}}s^{3/2}\right)\leq C\,\exp(-s).
\end{equation*}
If $s_0\leq 0$, for all $s\in [s_0, 0]$, we obtain from
\eqref{eq:kk-N-zz-minus-s} that when $N\geq N_0(s_0,\gamma)$,
$\kk_N^{2/3}\zz \in [3s_0/2, 1]$. Hence, we have
\begin{equation*}
  \EE^{-1}(\kk_N^{2/3}\zz)\exp(s) \leq C(s_0) \equiv \max_{s\in
    [3s_0/2, 1]} e\,\EE^{-1}(s),
\end{equation*}
the right hand side of which is, by its definition, continuous and
non-increasing. Therefore, we conclude that when $N\geq
N_0(s_0,\gamma)$,
\begin{equation}
\label{eq:D-n-N-1-2}
  \labs \Ai(\kk_N^{2/3}\zz) + \ee_2(\kk_N, \xi) \rabs \leq
  C(s_0)\exp(-s), \quad \text{for all $s\in I_{1,N}$.}
\end{equation}

Finally, putting the bounds \eqref{eq:D-n-N-1-1} and
\eqref{eq:D-n-N-1-2} together and recalling that $|r_N|$ could be
bounded by $2$, we obtain that under the condition of Lemma
\ref{lemma:laguerre}, when $N\geq N_0(s_0,\gamma)$, on $I_{1,N}$,
\begin{equation*}
    \labs D^1_{n,N} \rabs\leq C(s_0)N^{-2/3}\exp(-s).
\end{equation*}

\noindent \textit{(ii)} For $D^2_{n,N}$, we first split and control
$|r_NR_N^{-1}(\xi) - 1|$ as
\begin{equation*}
    \labs r_NR_N^{-1}(\xi)-1 \rabs \leq r_N\labs R_N^{-1}(\xi)-1
    \rabs + |r_N-1| = r_N |R_N(\xi)|^{-1}|R_N(\xi) - 1|+ |r_N - 1|.
\end{equation*}
By \eqref{eq:r-N} and \eqref{eq:R-N-minus-1}, when $N\geq
N_0(s_0,\gamma)$, we have $|r_N|\leq 2$, $|R_N(\xi)|^{-1}\leq 2$ and
hence
\begin{equation}
\label{eq:D-n-N-2-1}
  \labs r_NR_N^{-1}(\xi)-1 \rabs \leq C\,N^{-2/3}s + CN^{-1}\leq
  CN^{-2/3}s,\quad \text{for all } s\in I_{1,N}.
\end{equation}
On the other hand, by \eqref{eq:airy-bd}, we obtain
\begin{equation*}
    \labs \Ai'(\kk_N^{2/3}\zz)\rabs \leq
    \NN(\kk_N^{2/3}\zz)\EE^{-1}(\kk_N^{2/3}\zz).
\end{equation*}

When $s\geq 0$, we have from \eqref{eq:kk-N-zz-minus-s} that
$\kk_N^{2/3}\zz \in [s/2, 3s/2]$, and using \eqref{eq:EMN}, we obtain
\begin{equation*}
  \NN(\kk_N^{2/3}\zz)\EE^{-1}(\kk_N^{2/3}\zz) \leq C\,
  (\kk_N^{2/3}\zz)^{1/4}\exp\left(-\frac{1}{3\sqrt{2}}s^{3/2}\right)
  \leq C\, s^{1/4}\exp\left(-\frac{1}{3\sqrt{2}}s^{3/2}\right)\leq
  C\,\exp(-3s/2).
\end{equation*}

If $s_0\leq 0$, we know that when $N\geq N_0(s_0,\gamma)$,
$\kk_N^{2/3}\zz\in [3s_0/2, 1]$ for all $s\in [s_0, 0]$. We then have
\begin{equation*}
  \NN(\kk_N^{2/3}\zz)\EE(\kk_N^{2/3}\zz)\exp(3s/2) \leq C(s_0)\equiv
  \max_{s\in [3s_0/2, 1]} e^{3/2}\NN(s)\EE(s),
\end{equation*}
the right hand side of which is again continuous and non-increasing in
$s_0$. As before, this enables us to conclude that when $N\geq
N_0(s_0,\gamma)$,
\begin{equation}
\label{eq:D-n-N-2-2}
  \labs \Ai'(\kk_N^{2/3}\zz)\rabs C(s_0)\exp(-3s/2),\quad \text{for
    all $s\in I_{1,N}$.}
\end{equation}
Assembling \eqref{eq:D-n-N-2-1} and \eqref{eq:D-n-N-2-2}, we obtain
that when $N\geq N_0(s_0,\gamma)$,
\begin{equation*}
  \labs D^2_{n,N}\rabs \leq C(s_0)N^{-2/3}|s|\exp(-3s/2)\leq
  C(s_0)N^{-2/3}\exp(-s),\quad \text{for all $s\in I_{1,N}$.}
\end{equation*}

\noindent \textit{(iii)} For $D^3_{n,N}$, recalling that $r_N = 1 +
\Oh{N^{-1}}$ and we obtain the following bound under the condition of
Lemma \ref{lemma:laguerre} by using the previously derived bound on
$\Ai'(\kk_N^{2/3}\zz)$:
\begin{equation*}
  \labs D^3_{n,N} \rabs \leq C(s_0)\exp(-s).
\end{equation*}

\noindent \textit{(iv)} For $D^4_{n,N}$, by the definition of $R_N$
and $\dot{\zz}_N$ as well as the bound for $\partial_\xi\ee_2(\kk_N,
\xi)$, we have
\begin{equation*}
  \begin{split}
    \labs D^4_{n,N} \rabs & = \labs
    \left(\frac{\sigma_{n,N}}{\kk_N}\right)
    r_NR_N(\xi)\partial_\xi\ee_2(\kk_N, \xi)\rabs \\
    & \leq
    C\,N^{-2/3}\sigma_{n,N}\kk_N^{-1/3}r_N
    R_N(\xi)\dot{\zz}(\xi)\NN(\kk_N^{2/3}\zz)\EE^{-1}(\kk_N^{2/3}\zz)\\
    & =
    C\,N^{-2/3}r_NR_N^{-1}(\xi)\NN(\kk_N^{2/3}\zz)\EE^{-1}(\kk_N^{2/3}\zz).
  \end{split}
\end{equation*}
All the terms involved in the last bound have been well studied during
our analysis of $D^2_{n,N}$, and applying various results established
there, we obtain that when $N\geq N_0(s_0,\gamma)$,
\begin{equation*}
  \labs D^4_{n,N}\rabs \leq C(s_0)N^{-2/3}\exp(-s),\quad \text{for all
  $s\in I_{1,N}$}.
\end{equation*}

Combining the bounds for the four terms, we obtain from a simple
triangular inequality that when $N\geq N_0(s_0,\gamma)$,
\begin{equation*}
  T_{N,1} \leq C(s_0)\exp(-s), \quad \text{for $s\in
    I_{1,N}$.}
\end{equation*}

We remark that, here and after, we derive a more stringent bound with
the rate term $N^{-2/3}$ whenever possible. Although it is not
necessary here, those bounds with this rate term will become useful in
the later study of $|\psi'_\tau(s)-\Ai'(s)|$.

\emph{Case $s\in I_{2,N}$.} In this case, we define
$\tilde{D}^1_{n,N} = D^1_{n,N}$ and $\tilde{D}^2_{n,N} = D^2_{n,N} +
D^3_{n,N} + D^4_{n,N}$.

\noindent \textit{(i)} To analyze the $\tilde{D}^1_{n,N}$ term, we
first introduce a useful lemma:
\begin{lemma}
  \label{lemma:sqrt-f}
Let $r>0$ be fixed. For $x = x_N(s) = \xnexpr{s}$ and $\xi = x/\kk_N$,
when $s\geq r^2$, we have
\begin{equation*}
  \sigma_{n,N}\sqrt{f(\xi)} \geq r\xi_+/\xi =
  r\mu_{n,N}/(\xnexpr{s}).
\end{equation*}
\end{lemma}
For $\tilde{D}^1_{n,N}$, we could bound it for large $N$ as
\begin{equation*}
  \labs \tilde{D}^1_{n,N} \rabs\leq C\,r_N
  \left(\frac{\sigma_{n,N}}{\kk_N}\right)
  \labs \frac{R'_N(\xi)}{R_N(\xi)} \rabs
  R_N(\xi)\MM(\kk_N^{2/3}\zz)\EE^{-1}(\kk_N^{2/3}\zz).
\end{equation*}
We consider first the $R_N(\xi)\MM(\kk_N^{2/3}\zz)$ term. Recall that
$R_N(\xi) = \kk_N^{1/6}\sigma_{n,N}^{-1/2}\hat{f}^{-1/4}(\xi)$ and
that $|\MM(\kk_N^{2/3}\zz)|\leq C\kk_N^{-1/6}\zz^{-1/4}$ when $N$ is
large. Applying Lemma \ref{lemma:sqrt-f}, we obtain that when $N\geq
N_0(\gamma)$,
\begin{equation*}
  \begin{split}
    R_N(\xi)\MM(\kk_N^{2/3}\zz) & \leq
    C\zz^{-1/4}\sigma_{n,N}^{-1/2}\hat{f}^{-1/4}(\xi) =
    Cf^{-1/4}(\xi)\sigma_{n,N}^{-1/2}  \\
    & \leq
    Cr^{-1/2}\left(\frac{\mu_{n,N}}{\xnexpr{s}}\right)^{-1/2}\leq Cs, \quad
    \text{for all } s\in I_{2,N}.
  \end{split}
\end{equation*}
We remark that our choice of $s_1$ ensures that $s_1\geq r^2$
with $r = 1$.

Switching to the term $|R'_N(\xi)/R_N(\xi)|$, from the definition, we
have
\begin{equation*}
  \frac{R'_N(\xi)}{R_N(\xi)} =
  -\frac{\ddot{\zz}(\xi)}{2\dot{\zz}(\xi)}, \quad \text{and}\quad
  \frac{\ddot{\zz}(\xi)}{\dot{\zz}(\xi)} =
  \frac{f'(\xi)}{2f(\xi)}-\frac{\sqrt{f(\xi)}}{3I(\sqrt{f})}
\end{equation*}
where $I(\sqrt{f}) = \int_{\xi_+}^{\xi}\sqrt{f}$. Simple triangular
inequality gives a direct bound as
\begin{equation*}
  \labs \frac{R'_N(\xi)}{R_N(\xi)} \rabs \leq \frac{1}{4}\labs
  \frac{f'(\xi)}{f(\xi)} \rabs + \frac{1}{6}\frac{\sqrt{f(\xi)}}{I(\sqrt{f})}.
\end{equation*}
For the first term on the right hand side, simple manipulation gives
us
\begin{equation*}
  \labs \frac{f'(\xi)}{f(\xi)} \rabs = \labs \frac{1}{\xi-\xi_+} +
  \frac{1}{\xi-\xi_-} - \frac{2}{\xi} \rabs \leq \frac{4}{\xi - \xi_+}
  = \frac{4\kk_N}{s\sigma_{n,N}}\leq C\frac{\kk_N}{\sigma_{n,N}}.
\end{equation*}
Moreover, we could bound $(\xi-\xi_+)[\sqrt{f}/I(\sqrt{f})]$ as
\begin{equation*}
  \begin{split}
    \frac{(\xi-\xi_+)\sqrt{f(\xi)}}{I(\sqrt{f})} & \leq
    \frac{(\xi-\xi_+)^{3/2}(\xi-\xi_-)^{1/2}}{2\xi\int_{\xi_+}^{\xi}\frac{\sqrt{t-\xi_+}(t-\xi_-)}{2t}dt}\leq
    \frac{(\xi-\xi_+)^{3/2}}{(1-\xi_-/\xi_+)\int_{\xi_+}^{\xi}\sqrt{t-\xi_+}dt}\\
    & = \frac{3}{2(1-\xi_-/\xi_+)} \leq \frac{6}{\xi_+ - \xi_-} \leq
    \frac{3}{4}\left(1+\frac{2n}{N}\right)
  \end{split}
\end{equation*}
Hence, when $N\geq N_0(\gamma)$, we obtain the bound for
$\sqrt{f}/I(\sqrt{f})$ as
\begin{equation*}
  \frac{\sqrt{f(\xi)}}{I(\sqrt{f})} \leq
  \frac{3}{4}\left(\frac{1+2n/N}{\xi - \xi_+}\right) \leq
  \frac{3}{4}\left(1+\frac{2n}{N}\right)\frac{\kk_N}{\sigma_{n,N}}s_1^{-1}N^{-1/6}\leq
  C\frac{\kk_N}{\sigma_{n,N}}.
\end{equation*}
This implies that $|R'_N(\xi)/R_N(\xi)|$ is bounded by
$C\kk_N/\sigma_{n,N}$ which further ensures
\begin{equation*}
  r_N\left(\frac{\sigma_{n,N}}{\kk_N}\right)\labs
  \frac{R'_N(\xi)}{R_N(\xi)} \rabs \leq C.
\end{equation*}
Finally, using \eqref{eq:kk-N-zz-s-1} and the fact that $s_1$ is a
fixed constant, we obtain that when $N\geq N_0(\gamma)$,
\begin{equation*}
  \labs \tilde{D}^1_{n,N} \rabs \leq Cs\exp(-3s/2)\leq Cs^{-4}\exp(-s)
  \leq CN^{-2/3}\exp(-s),\quad \text{for all $s\in I_{2,N}$.}
\end{equation*}

\noindent \textit{(ii)} For $\tilde{D}^2_{n,N}$, we first recall its
definition as
\begin{equation*}
  \tilde{D}^2_{n,N} = r_N R_N(\xi)\left[ \Ai'(\kk_N^{2/3}\zz)R_N^{-2}(\xi) + \left(\frac{\sigma_{n,N}}
    {\kk_N}\right) \frac{\partial}{\partial \xi} \ee_2(\kk_N, \xi)
    \right].
\end{equation*}
By definition of $R_N$ and the large $N$ bounds on $r_N$,
$\partial_\xi\ee_2(\kk_N,\xi)$ and $\Ai'$, we have
\begin{equation*}
  \labs \tilde{D}^2_{n,N} \rabs \leq
  C\,R_N^{-1}(\xi)\NN(\kk_N^{2/3}\zz)\EE^{-1}(\kk_N^{2/3}\zz).
\end{equation*}
The asymptotics of the phase function $\NN$ suggest that
\begin{equation*}
  R_N^{-1}(\xi)\NN(\kk_N^{2/3}\zz)\leq
  CR_N^{-1}(\xi)\kk_N^{1/6}\zz^{1/4} =
  Cf^{1/4}(\xi)\sigma_{n,N}^{1/2}.
\end{equation*}
For $\sigma_{n,N}\sqrt{f(\xi)}$, we could simply bound it as
\begin{equation*}
  \sigma_{n,N}\sqrt{f(\xi)} =
  \frac{\sigma_{n,N}\sqrt{(\xi-\xi_+)(\xi-\xi_-)}}{2\xi}\leq
  \frac{\sigma_{n,N}}{2}.
\end{equation*}
Observing that for $s\in I_{2,N}$, $\sigma_{n,N}\leq
C(\gamma)N^{1/3}\leq Cs^4$, we obtain
\begin{equation*}
  R_N^{-1}(\xi)\NN(\kk_N^{2/3}\zz)\leq
  C\sigma_{n,N}^{1/2}f^{1/4}(\xi)\leq C\sigma_{n,N}^{1/2}\leq Cs^2.
\end{equation*}
Once more, by \eqref{eq:kk-N-zz-s-1} and our choice of $s_1$ [see
\ref{subsec:a-s-1}], we obtain
\begin{equation*}
  \labs \tilde{D}^2_{n,N}\rabs \leq Cs^2\exp(-3s/2) \leq
  Cs^{-4}\exp(-s)\leq C\,N^{-2/3}\exp(-s).
\end{equation*}
This finally gives a bound of the form $C\,N^{-2/3}\exp(-s)$ for
$T_{N,1}$ on $I_{2,N}$.

By a simple triangular inequality, we combine our bounds on
$T_{N,1}$ and $T_{N,2}$ on both $I_{1,N}$ and $I_{2,N}$ together and
obtain that under the condition of Lemma \ref{lemma:laguerre}, when
$N\geq N_0(s_0,\gamma)$,
\begin{equation*}
  \labs \partial_s\theta_{n,N}(x_{n,N}(s)) \rabs \leq
  C(s_0)\exp(-s),\quad \text{for all $s\geq s_0$.}
\end{equation*}

\paragraph{Bound for $|\psi'_\tau(s)|$.}
We have pointed out that $\psi_\tau$ is of the form
$\frac{1}{\sqrt{2}}\rho_{N}\theta_{n-1,N-1}(x_{n-1,N-1}(s))$ with
$\rho_N = 1 + \Oh{N^{-1}}$. Hence, we have $\psi'_\tau(s)$ as
\begin{equation*}
  \psi'_\tau(s) = \frac{1}{\sqrt{2}}\rho_N
  \partial_s\theta_{n-1,N-1}(x_{n-1,N-1}(s)),
\end{equation*}
for which our bound on $\sigma_{n,N}\partial_s\theta_{n,N}(s)$ apply
directly and we obtain that under the condition of Lemma
\ref{lemma:laguerre}, when $N\geq N_0(s_0,\gamma)$,
\begin{equation*}
  \labs \psi'_\tau(s)  \rabs\leq C(s_0)\exp(-s),\quad \text{for all $s\geq
  s_0$.}
\end{equation*}

\paragraph{Bound for $|\psi'_\tau(s) - G'(s)|$.}
By the expression of $\psi'_\tau$, we could split
$|\psi'_\tau(s)-G'(s)|$ as
\begin{equation}
\label{eq:psi'-tau-G'}
\begin{split}
  \labs \psi'_\tau(s)-G'(s) \rabs \leq &
  \frac{1}{\sqrt{2}}|\rho_N-1|\labs \partial_s\theta_{n-1,N-1}(x_{n-1,N-1}(s))
  \rabs \\ & + \frac{1}{\sqrt{2}}\labs \partial_s\theta_{n-1,N-1}(x_{n-1,N-1}(s))-\Ai'(s) \rabs.
\end{split}
\end{equation}
By our bound on $|\partial_s\theta_{n-1,N-1}(x_{n-1,N-1}(s))|$ and
recalling that $\rho_N = 1+\Oh{N^{-1}}$, the first term is then
bounded by $C(s_0)N^{-1}\exp(-s)$. We focus on the quantity
$|\partial_s\theta_{n,N}(x_{n,N}(s))-\Ai'(s)|$ to bound the second
term.

We split the quantity of interest into two parts as the following:
\begin{equation*}
  \begin{split}
    \labs \partial_s\theta_{n,N}(x_{n,N}(s))-\Ai'(s) \rabs & \leq \labs
    \sigma_{n,N}F'_{n,N}(x_{n,N}(s))\frac{\mu_{n,N}}{x_{n,N}(s)} -
    \Ai'(s)\rabs + \labs
    \sigma_{n,N}F_{n,N}(x_{n,N}(s))\frac{\mu_{n,N}}{x_{n,N}^2(s)} \rabs
    \\
    & = \mathcal{T}_{N,1}(s) + \mathcal{T}_{N,2}(s).
  \end{split}
\end{equation*}
The $\mathcal{T}_{N,2}(s)$ term is exactly the same as $T_{N,2}(s)$
defined in the previous study of
$\partial_s\theta_{n,N}(x_{n,N}(s))$ and hence we quote the bound
derived there directly as
\begin{equation*}
  \mathcal{T}_{N,2}(s) \leq C(s_0)N^{-2/3}\exp(-s),\quad \text{for all }s\geq
  s_0.
\end{equation*}

Switching to the $\mathcal{T}_{N,1}(s)$ term, we divide the whole region into
the two disjoint intervals $I_{1,N} = [s_0, s_1N^{1/6})$ and $I_{2,N}
= [s_1N^{1/6},\infty)$ again.

\emph{Case $s\in I_{1,N}$.} Exploiting a similar strategy in
splitting $\sigma_{n,N}F'_{n,N}(x)$, on $I_{1,N}$, we decompose
$\mathcal{T}_{N,1}(s)$ as $ \mathcal{T}_{N,1}(s) =
\sum_{i=1}^5\mathcal{D}^i_{n,N}$, with $\mathcal{D}^i_{n,N} =
D^i_{n,N}\mu_{n,N}/x_{n,N}(s)$ for $i = 1, 2$ and $4$,
\begin{equation*}
  \mathcal{D}^3_{n,N} =
  r_N\frac{\mu_{n,N}}{x_{n,N}(s)}\left[\Ai'(\kk_N^{2/3}\zz) -
    \Ai'(s)\right],\quad \text{and} \quad \mathcal{D}^5_{n,N} =
  \left[r_N\frac{\mu_{n,N}}{x_{n,N}(s)}-1\right]\Ai'(s).
\end{equation*}

For $i=1,2$ and $4$, using our previous bounds on $D^i_{n,N}$ and
noting that $|\mu_{n,N}/x_{n,N}(s)|$ could be bounded by $2$ on
$I_{1,N}$, we obtain directly that, when $N\geq N_0(s_0,\gamma)$,
\begin{equation*}
  \labs \mathcal{D}^i_{n,N}\rabs \leq C(s_0)N^{-2/3}\exp(-s),\quad \text{for $i = 1,2$
    and $4$, and all $s\in I_{1,N}$.}
\end{equation*}

For $\mathcal{D}^3_{n,N}$, by a first order Taylor expansion and the
identity $\Ai"(s) = s\Ai(s)$ for all $s$, we have that, for some
$s^*$ in the middle of $\kk_N^{2/3}\zz$ and $s$,
\begin{equation*}
    \labs\mathcal{D}^3_{n,N}\rabs  =
    r_N\labs\frac{\mu_{n,N}}{x_{n,N}(s)}\rabs\labs s^*\Ai(s^*)\rabs\labs\kk_N^{2/3}\zz
    - s\rabs
    \leq  C\,N^{-2/3}s^2\labs s^*\Ai(s^*)\rabs,
\end{equation*}
where the inequality holds when $N\geq N_0(s_0,\gamma)$ and comes
from \eqref{eq:kk-N-zz-minus-s} and the large $N$ bounds for $r_N$
and $\mu_{n,N}/x_{n,N}(s)$.

When $s\geq 0$, we know from the definition of $\zz$ that
$\kk_N^{2/3}\zz \geq 0$ and hence $s^*\geq 0$. Moreover,
\eqref{eq:kk-N-zz-minus-s} implies that when $N$ is large,
$\kk_N^{2/3}\zz$ and hence $s^*$ will be greater than $s/2$. Thus,
by \eqref{eq:EMN} and the monotonicity of $\EE$, we obtain
\begin{equation*}
    \labs s^*\Ai(s^*)\rabs \leq C\,s\EE^{-1}(s/2)\leq
    C\,s\exp\left( -\frac{1}{3\sqrt{2}}s^{3/2} \right)\leq
    C\,\exp(-3s/2).
\end{equation*}

If $s_0\leq 0$, as before, we consider all $s\in [s_0,0]$. Once
again, we obtain from \eqref{eq:kk-N-zz-minus-s} that for large $N$,
$\kk_N^{2/3}\zz\in [3s_0/2,1]$ and hence $s^*\in [3s_0/2, 1]$. Then
for all $s\in [s_0,0]$, when $N\geq N_0(s_0,\gamma)$,
\begin{equation*}
\labs s^*\Ai(s^*)\rabs \exp(3s/2) \leq C(s_0) \equiv \max_{s\in
[3s_0/2,1]} e^{3/2}\labs s\Ai(s) \rabs.
\end{equation*}
This $C(s_0)$ is continuous and non-increasing in $s_0$.

Thus, we could conclude that when $N\geq N_0(s_0,\gamma)$, for all
$s\in I_{1,N}$,
\begin{equation*}
    \labs\mathcal{D}^3_{n,N}\rabs  \leq  C\,N^{-2/3}s^2\labs
    s^*\Ai(s^*)\rabs \leq C(s_0)N^{-2/3}s^2\exp(-3s/2)\leq
    C(s_0)N^{-2/3}\exp(-s).
\end{equation*}

In $\mathcal{D}^5_{n,N}$, recalling
$\sigma_{n,N}/\mu_{n,N}=\Oh{N^{-2/3}}$ and $r_N = 1 + \Oh{N^{-1}}$,
we have that, when $N\geq N_0(s_0,\gamma)$, for all $s\in I_{1,N}$,
$\labs s_0 + \mu_{n,N}/\sigma_{n,N} \rabs \geq
\frac{1}{2}(\mu_{n,N}/\sigma_{n,N})$ and hence
\begin{equation*}
  \begin{split}
    \labs r_N\frac{\mu_{n,N}}{x_{n,N}(s)}-1\rabs & \leq r_N\labs
    \frac{\mu_{n,N}}{x_{n,N}(s)}-1\rabs + \labs r_N - 1 \rabs
    = r_N|s|\labs s+ \frac{\mu_{n,N}}{\sigma_{n,N}} \rabs^{-1} +
    \labs r_N-1 \rabs \\ & \leq r_N|s|\labs s_0 + \frac{\mu_{n,N}}{\sigma_{n,N}} \rabs^{-1} +
    \labs r_N-1 \rabs \leq C\,N^{-2/3}|s|+C\,N^{-1}.
  \end{split}
\end{equation*}
For $\Ai'(s)$, by \eqref{eq:EMN} and \eqref{eq:airy-bd}, we obtain
directly that
\begin{equation*}
  \labs \Ai'(s) \rabs \leq
  C(s_0)|s|^{1/4}\exp\left(-\frac{2}{3}s^{3/2}\right),
\end{equation*}
where $C(s_0)$ could be chosen as
\begin{equation*}
    \max_{s\in [s_0,\infty)} \labs \Ai'(s) \rabs \left(1 +
    |s|^{1/4}\right)^{-1}\exp\left(\frac{2}{3}s^{3/2}\right),
\end{equation*}
which is continuous and non-increasing.

Putting two parts together, we obtain that for all $s\in I_{1,N}$,
\begin{equation*}
  \labs \mathcal{D}^5_{n,N}\rabs \leq C(s_0)N^{-2/3}\left( |s| + CN^{-1/3}
  \right)|s|^{1/4}\exp\left(-\frac{2}{3}s^{3/2}\right) \leq
  C(s_0)N^{-2/3}\exp(-s).
\end{equation*}
We could then assemble all the bounds on $\mathcal{D}^i_{n,N}$ using
the triangular inequality and conclude that under the condition of
Lemma \ref{lemma:laguerre}, when $N\geq N_0(s_0,\gamma)$,
\begin{equation*}
  \mathcal{T}_{N,1}(s)\leq C(s_0)N^{-2/3}\exp(-s),\quad \text{for all }s\in
  I_{1,N}.
\end{equation*}

\emph{Case $s\in I_{2,N}$.} In this case, we could act more
heavy-handedly. In particular, by the asymptotics of $T_{N,1}(s)$ on
$I_{2,N}$ and the asymptotics of $\Ai'$, we have
\begin{equation*}
  \begin{split}
    \mathcal{T}_{N,1}(s) & \leq \labs
    \sigma_{n,N}F'_{n,N}(x_{n,N}(s))\frac{\mu_{n,N}}{x_{n,N}(s)} \rabs
    + \labs \Ai'(s)\rabs
    \leq CN^{-2/3}\exp(-s) + Cs^{1/4}\exp(-3s/2) \\
    & \leq CN^{-2/3}\exp(-s) + CN^{-2/3}s^{4+1/4}\exp(-3s/2)
    \leq CN^{-2/3}\exp(-s).
  \end{split}
\end{equation*}

We then obtain the bound $C(s_0)N^{-2/3}\exp(-s)$ for
$\mathcal{T}_{N,1}(s)$ and hence also for
$|\partial_s\theta_{n,N}(x_{n,N}(s))-\Ai'(s)|$ for all $s\in
[s_0,\infty)$. Applying the bound to the second term in
\eqref{eq:psi'-tau-G'}, we obtain that under the condition of Lemma
\ref{lemma:laguerre}, when $N\geq N_0(s_0,\gamma)$,
\begin{equation*}
  \labs \psi'_\tau(s)-G'(s)\rabs \leq C(s_0)N^{-2/3}\exp(-s),\quad
  \text{for all $s\in [s_0,\infty)$.}
\end{equation*}

\paragraph{Improved bound for $|\psi_\tau - G|$.}
  The above bound on $|\psi'_\tau(s)-G'(s)|$ could be used to derive a
  more stringent bound for $|\psi_\tau(s)-G(s)|$ as the following:
  \begin{equation}\label{eq:psitau-g-improved}
    \begin{split}
      \labs \psi_\tau(s) -G(s)\rabs & = \labs
      \int_s^{2s}[\psi'_\tau(t)-G'(t)]dt -\left[ \psi_\tau(2s) - G(2s)
      \right] \rabs \\
      & \leq \int_s^{2s}\labs \psi'_\tau(t) - G'(t) \rabs dt + \labs
      \psi_\tau(2s)-G(2s) \rabs \\
      & \leq \int_s^{2s}C(s_0)N^{-2/3}e^{-t}dt +
      C(s_0)N^{-2/3}\exp(-s)
      \leq C(s_0)N^{-2/3}\exp(-s).
    \end{split}
  \end{equation}
This is exactly the bound that we have claimed in Lemma
\ref{lemma:laguerre}.

\subsubsection{Bounds for quantities related to $\phi_\tau(s)$}
\label{subsec:phi-bd}
In this part, we employ a trick that was first used in
\citet[p.320]{johnstone01} to derive bounds for quantities related
to $\phi_\tau$ from those for quantities related to $\psi_\tau$.

Recall that $\phi_\tau$ could be expressed as
\begin{equation*}
  \phi_\tau(s) =
  \frac{1}{\sqrt{2}}\tilde{\rho}_NF_{n-2,N}(x_{n-1,N-1}(s))\frac{\mu_{n-2,N}}
  {x_{n-1,N-1}(s)},
\end{equation*}
where $\tilde{\rho}_N = 1+\Oh{N^{-1}}$ [see \ref{subsec:a-constants}
for its proof]. The problem of $\phi_\tau$ is that the centering and
scaling constants $(\mu_{n-1,N-1},\sigma_{n-1,N-1})$ in the
transformation $x_{n-1,N-1}(s)$ does not agree with the ``optimal''
constants $(\mu_{n-2,N}, \sigma_{n-2,N})$ for the related function
$F_{n-2,N}$. To circumvent this problem, we introduce a new
independent variable $s'$ as the following [one should not confuse
it with the $s'$ previously appeared in Section \ref{sec:proof}]:
\begin{equation}
  \label{eq:s-prime-def}
  \mu_{n-1,N-1}+s\sigma_{n-1,N-1} = \mu_{n-2,N}+s'\sigma_{n-2,N}.
\end{equation}
Then $s' =
(\mu_{n-1,N-1}-\mu_{n-2,N})/\sigma_{n-2,N}+s\sigma_{n-1,N-1}/\sigma_{n-2,N}$. By
defining
\begin{equation}
\label{eq:Delta-N}
  \Delta_N = \frac{\mu_{n-1,N-1}-\mu_{n-2,N}}{\sigma_{n-2,N}},
\end{equation}
we have $s'-s = \Delta_N + [\sigma_{n-1,N-1}\sigma_{n-2,N}^{-1}]s$ and
$\phi_\tau(s)$ could be rewritten as
\begin{equation*}
  \phi_\tau(s)=\frac{1}{\sqrt{2}}\tilde{\rho}_NF_{n-2,N}(x_{n-2,N}(s'))
  \frac{\mu_{n-2,N}}{x_{n-2,N}(s')}.
\end{equation*}
Before we proceed, we list two important properties as the following
[with proof given in \ref{subsec:a-constants}]:
\begin{equation}
\label{eq:Delta-N-bd}
  \Delta_N = \Oh{N^{-1/3}}\quad \text{and} \quad 1\leq
  \frac{\sigma_{n-1,N-1}}{\sigma_{n-2,N}}=1+\Oh{N^{-1}}.
\end{equation}

\paragraph{Bounds for $|\phi_\tau(s)|$ and $|\phi'_\tau(s)|$. }
Applying our previous bounds for $|\theta_{n,N}(x_{n,N}(s))|$ and
$|\partial_s\theta_{n,N}(x_{n,N}(s))|$, and using
\eqref{eq:Delta-N-bd}, we obtain that under the condition of Lemma
\ref{lemma:laguerre}, for all $s\in [s_0,\infty)$
\begin{equation*}
  \labs \phi_\tau(s)\rabs = \frac{1}{\sqrt{2}}\tilde{\rho}_N\labs
  \theta_{n-2,N}(x_{n-2,N}(s')) \rabs \leq C(s_0)\exp(-s')\leq
  C(s_0)\exp(-s);
\end{equation*}
\begin{equation*}
  \begin{split}
    \labs \phi'_\tau(s)\rabs & = \frac{1}{\sqrt{2}}\tilde{\rho}_N\labs
    \partial_s\theta_{n-2,N}(x_{n-2,N}(s'))\rabs =
    \frac{1}{\sqrt{2}}\tilde{\rho}_N\labs \partial_{s'}\theta_{n-2,N}(x_{n-2,N}(s'))
    \rabs \frac{ds'}{ds} \\
    & \leq C(s_0)\exp(-s')\frac{\sigma_{n-1,N-1}}{\sigma_{n-2,N}}\leq
    C(s_0)\exp(-s).
  \end{split}
\end{equation*}

\paragraph{Bounds for $|\phi_\tau(s)-G(s)-\Delta_NG'(s)|$ and
  $|\phi'_\tau(s)-G'(s)-\Delta_NG"(s)|$.}

We consider $|\phi_\tau(s)-G(s)-\Delta_NG'(s)|$ in detail and the
derivation for $|\phi'_\tau(s)-G'(s)-\Delta_NG"(s)|$ is essentially
the same.

By our definition of $s'$ and recalling that $\Ai"(s) = s\Ai(s)$, we
have the Taylor expansion of $G(s')$ as
\begin{equation*}
  \begin{split}
    G(s') & = G(s) + (s'-s)G'(s) + \frac{1}{2}(s'-s)^2G"(s^*) \\
    & = G(s) + \Delta_NG'(s) +
    \frac{1}{\sqrt{2}}\left[\frac{\sigma_{n-1,N-1}}{\sigma_{n-2,N}}-1
    \right]s\Ai'(s) + \frac{1}{2\sqrt{2}}(s'-s)^2s^*\Ai(s^*),
  \end{split}
\end{equation*}
with $s^*$ lies at somewhere between $s$ and $s'$.

Hence, by our bound on $|\psi_\tau(s)-G(s)|$, we obtain
\begin{equation}
\label{eq:phi-tau-G}
  \labs \phi_\tau(s) - G(s) - \Delta_NG'(s) \rabs \leq
  C(s_0)N^{-2/3}\exp(-s')+CN^{-1}\labs s\Ai'(s) \rabs + C(s'-s)^2\labs
  s^*\Ai(s^*)\rabs.
\end{equation}
On $[s_0,\infty)$, we have $C(s_0)N^{-2/3}\exp(-s')\leq
C(s_0)N^{-2/3}\exp(-s)$ for the first term in the above bound.
Moreover, by \eqref{eq:EMN} and \eqref{eq:airy-bd}, the second term
satisfies
\begin{equation*}
  CN^{-1}\labs s\Ai'(s)\rabs \leq
  C(s_0)N^{-1}|s|^{1+1/4}\exp\left(-\frac{2}{3}s^{3/2}\right)\leq
  C(s_0)N^{-1}\exp(-s),\quad \text{for all $s\geq s_0$.}
\end{equation*}
For the last term, we split $[s_0,\infty)$ into $I_{1,N}\cup
I_{2,N}$ as usual. For $s\in I_{1,N}$, when $N\geq N_0(s_0,\gamma)$,
\begin{equation*}
    (s-s')^2 = \left[ \Delta_N + \left(\frac{\sigma_{n-1,N-1}}
    {\sigma_{n-2,N}} - 1\right)s \right]^2 \leq \left[ CN^{-1/3} + CN^{-1}s
    \right]^2 \leq \left(CN^{-2/3}\right)\wedge 1.
\end{equation*}
We obtain from the above bound that $|s^*-s|\leq 1$ and hence by
\eqref{eq:EMN} and \eqref{eq:airy-bd},
\begin{equation*}
  C(s-s')^2\labs s^*\Ai(s^*) \rabs \leq
  C\,N^{-2/3}(|s|+1)\EE^{-1}(s-1)\leq
  C(s_0)N^{-2/3}\exp(-s),
\end{equation*}
where $C(s_0)$ could be chosen as $\max_{s\in [s_0,\infty)}
Ce^s\EE^{-1}(s-1)$.

On $I_{2,N}$, we have $s'\geq s/2$ from \eqref{eq:s-prime-s/2} and
hence $s^*\geq s/2$. By \eqref{eq:EMN} and \eqref{eq:airy-bd}, we
obtain that when $N\geq N_0(s_0,\gamma)$,
\begin{equation*}
  C(s'-s)^2\labs s^*\Ai(s^*)\rabs \leq
  Cs^3\exp\left(-\frac{1}{3\sqrt{2}}s^{3/2}\right)\leq
  Cs^{-4}\exp(-s)\leq CN^{-2/3}\exp(-s).
\end{equation*}

Therefore, we have shown that, for all $N\geq N_0(s_0,\gamma)$ and
$s\geq s_0$, the right hand side of \eqref{eq:phi-tau-G} is further
controlled by $C(s_0)N^{-2/3}\exp(-s)$, which is exactly the desired
bound for $\labs \phi_\tau(s) - G(s) - \Delta_NG'(s) \rabs$.

\appendix

\section{Technical Details}
\label{appendix}

\subsection{Properties of $\beta_N, \rho_N, \tilde{\rho}_N, \Delta_N$ and $\sigma_{n-1,N-1}/\sigma_{n-2,N}$}
\label{subsec:a-constants}

\subsubsection{Property of $\beta_N$}
\label{subsubsec:beta-N}

We are to show that
\begin{equation*}
  \beta_N = \frac{1}{\sqrt{2}} + \Oh{N^{-1}}.
\end{equation*}

First of all, we recall that $\phi_\tau(s)$ is defined to be $0$
when $\tau(s) = \tauexpr{s} < 0$, i.e., when $s\in (-\infty,
-\tmu_{n,N}/\tsigma_{n,N})$. Hence, we have
\begin{equation*}
  \begin{split}
    \beta_N & = \frac{1}{2}\int_{-\infty}^\infty \phi_\tau(s)ds =
    \frac{1}{2}\int_0^\infty \phi(x;\aa_N - 1)dx \\
    & =
    \frac{N^{1/4}(n-1)^{1/4}\Gamma^{1/2}(N+1)}{2\sqrt{2}\Gamma^{1/2}(n)}\int_0^\infty
    x^{\aa_N/2-1}e^{-x/2}L_N^{\aa_N-1}(x)dx \\
    & =
    \frac{2^{-(\aa_N-1)/2}N^{1/4}(n-1)^{1/4}\Gamma^{1/2}(n)\Gamma\left(\frac{N+3}{2}\right)}
    {(N+1)\Gamma^{1/2}(N+1)\Gamma\left(\frac{n+1}{2}\right)}.
  \end{split}
\end{equation*}
Applying Sterling's formula
\begin{equation*}
  \Gamma(z) =
  \Bigg(\frac{2\pi}{z}\Bigg)^{1/2}\Bigg(\frac{z}{e}\Bigg)^{z}\left[1+\Oh{\frac{1}{z}}\right],
\end{equation*}
we obtain that
\begin{equation*}
  \begin{split}
    \beta_N & =
    \frac{2^{-(\aa_N-1)/2}N^{1/4}(n-1)^{1/4}}{N+1}\frac{\Bigl(\frac{2\pi}{n}\Bigr)^{1/4}\Bigl(\frac{n}{e}\Bigr)^{n/2}
      \Bigl(\frac{4\pi}{N+3}\Bigr)^{1/2}\Bigl(\frac{N+3}{2e}\Bigr)^{(N+3)/2}}{\Bigl(\frac{2\pi}{N+1}\Bigr)^{1/4}
      \Bigl(\frac{N+1}{e}\Bigr)^{(N+1)/2}\Bigl(\frac{4\pi}{n+1}\Bigr)^{1/2}\Bigl(\frac{n+1}{2e}\Bigr)^{(n+1)/2}}
    \left(1 + \Oh{N^{-1}}\right) \\
    & = \frac{1}{\sqrt{2e}}\left(1 - \frac{1}{n+1}\right)^{n/2}\left(1
      + \frac{2}{N+1}\right)^{(N+1)/2 + 3/4} \left(1 +
      \Oh{N^{-1}}\right) = \frac{1}{\sqrt{2}} \left(1 +
      \Oh{N^{-1}}\right).
  \end{split}
\end{equation*}
The last equality is exactly the asymptotics that we need for $\beta_N$.

\subsubsection{Asymptotics of $\rho_N$ and $\tilde{\rho}_N$}
\label{subsubsec:rho-N}

In this part, we show that the asymptotics of $\rho_N$ and
$\tilde{\rho}_N$ satisfy
\begin{equation*}
  \rho_N, \tilde{\rho}_N = 1 + \Oh{N^{-1}}.
\end{equation*}
We consider $\rho_N$ first. By definition, we have
\begin{equation*}
  \rho_N =
  \frac{N^{1/4}(n-1)^{1/4}\sigma_{n-1,N-1}^{1/2}\tsigma_{n,N}}{\tmu_{n,N}}
  = \frac{N^{1/4}(n-1)^{1/4}\sigma_{n-1,N-1}^{3/2}}{\mu_{n-1,N-1}}.
\end{equation*}
Plugging in the definition of $\sigma_{n-1,N-1}$ and $\mu_{n-1,N-1}$,
we obtain that
\begin{equation*}
  \begin{split}
    \rho_N & =
    N^{1/4}(n-1)^{1/4}\left(\sqrt{N-\tfrac{1}{2}}+\sqrt{n-\tfrac{1}{2}}\right)^{-1/2}
    \left(\frac{1}{\sqrt{N-\frac{1}{2}}}+\frac{1}{\sqrt{n-\frac{1}{2}}}\right)^{1/2}
    \\
    & =
    \left(\frac{N}{N-\frac{1}{2}}\right)^{1/4}\left(\frac{n-1}{n-\frac{1}{2}}\right)^{1/4}
    = 1 + \Oh{N^{-1}}.
  \end{split}
\end{equation*}

For $\tilde{\rho}_N$, we have from its definition that
\begin{equation*}
  \begin{split}
    \tilde{\rho}_N & =
    \frac{N^{1/4}(n-1)^{1/4}\sigma_{n-2,N}^{1/2}\tsigma_{n,N}}{\mu_{n-2,N}}
    =
    \frac{\sigma_{n-1,N-1}}{\sigma_{n-2,N}}\frac{N^{1/4}(n-1)^{1/4}\sigma_{n-2,N}^{3/2}}{\mu_{n-2,N}}
    \\
    & = \frac{\sigma_{n-1,N-1}}{\sigma_{n-2,N}}N^{1/4}(n-1)^{1/4}\left(\sqrt{N+\tfrac{1}{2}}+\sqrt{n-\tfrac{3}{2}}\right)^{-1/2}
    \left(\frac{1}{\sqrt{N+\frac{1}{2}}}+\frac{1}{\sqrt{n-\frac{3}{2}}}\right)^{1/2}\\
    & =
    \frac{\sigma_{n-1,N-1}}{\sigma_{n-2,N}}\left(\frac{N}{N+\frac{1}{2}}\right)^{1/4}\left(\frac{n-1}{n-\frac{3}{2}}\right)^{1/4}
    = 1 + \Oh{N^{-1}}.
  \end{split}
\end{equation*}
The last equality holds since $\sigma_{n-1,N-1}/\sigma_{n-2,N} = 1 +
\Oh{N^{-1}}$ as claimed in \eqref{eq:Delta-N}, which is to be shown
below in \ref{subsubsec:Delta-N}.

\subsubsection{Properties of $\Delta_N$ and $\sigma_{n-1,N-1}/\sigma_{n-2,N}$}
\label{subsubsec:Delta-N}

We focus on $\Delta_N$ first. As a reminder, we recall its definition as
\begin{equation*}
  \Delta_N = \frac{\mu_{n-1,N-1}-\mu_{n-2,N}}{\sigma_{n-2,N}}.
\end{equation*}
By \citet[A.1.2]{nek06}, we have for the numerator that $\mu_{n-1,N-1}
- \mu_{n-2,N} = \Oh{1}$.  For the denominator, if we let denote
$\left(n-\tfrac{3}{2}\right)/\left(N+\tfrac{1}{2}\right)$ by
$\gamma_{n,N}$, we then have
\begin{equation*}
  \begin{split}
    \frac{1}{\sigma_{n-2,N}} & =
    \left(\sqrt{N+\tfrac{1}{2}}+\sqrt{n-\tfrac{3}{2}}\right)^{-1}
    \left(\frac{1}{\sqrt{N+\frac{1}{2}}}+\frac{1}{\sqrt{n-\frac{3}{2}}}\right)^{-1/3}
    \\
    & = \frac{1}{1 + \sqrt{\gamma_{n,N}}}\left(1 +
      \frac{1}{\sqrt{\gamma_{n,N}}}\right)\left(N +
      \frac{1}{2}\right)^{-1/3} = \Oh{N^{-1/3}}.
  \end{split}
\end{equation*}
The last equality holds since $\gamma_{n,N}$ is bounded below for
all $n > N$. Combining the two estimates, we establish that
\begin{equation*}
  \Delta_{N} = \Oh{N^{-1/3}}.
\end{equation*}

We now switch to prove that
\begin{equation*}
  1\leq \sigma_{n-1,N-1}/\sigma_{n-2,N} = 1 + \Oh{N^{-1}}.
\end{equation*}
The fact that $\sigma_{n-1,N-1}/\sigma_{n-2,N} = 1 + \Oh{N^{-1}}$
has been proved in \citet[A.1.3]{nek06}. On the other hand, we have
from the second last display of \citet[A.1.3]{nek06} that
\begin{equation*}
  \left( \frac{\sigma_{n-1,N-1}}{\sigma_{n-2,N}} \right)^3 = \left[ 1
    + \frac{\sqrt{n/N} - \sqrt{N/n}}{n+N} + \Oh{n^{-2}}\right]\left[1
    + \frac{1}{2}\left(\frac{1}{n}+\frac{1}{N}\right) +
    \Oh{n^{-2}}\right].
\end{equation*}
Both terms become greater than $1$ when $N\geq N_0(\gamma)$ and hence $\sigma_{n-1,N-1}\sigma_{n-2,N}^{-1} \geq 1$
for large $N$. Actually, the inequality holds for any $n>N\geq
2$. However, what we have proved here is sufficient for our argument
in Section \ref{subsec:phi-bd}.

\subsection{Evaluation of the entries of $K_\tau$}
\label{subsec:a-k-tau}

In this part, we work out the explicit expressions for the entries of
$K_\tau$ given in \eqref{eq:k-tau-entries}. To this end, we proceed
term by term.

\paragraph{$K_{\tau,11}$ term.} For $K_{\tau,11}$, we have from
its definition that
\begin{equation*}
  \begin{split}
    K_{\tau,11}(s,t) & = \sigma_{n,N}S_{N,1}(\tauexpr{s}, \tauexpr{t})
    \\
    & = \sigma_{n,N}\left[ S_{N,2}(\tau(s),\tau(t)) +
      \psi(\tau(s))(\ee\phi)(\tau(t)) \right] \\
    & = S_\tau(s,t) + \sigma_{n,N}\psi(\tau(s))(\ee\phi)(\tau(t)).
  \end{split}
\end{equation*}
For the second term in the last expression, we have
$\sigma_{n,N}\phi(\tau(s))=\psi_\tau(s)$ and
\begin{equation*}
  \int_y^\infty \phi(z)dz = \int_t^\infty \phi(\tau(u))\tau'(u)du =
  \tsigma_{n,N}\int_t^\infty \phi(\tau(u))du = \int_t^\infty \phi_\tau(u)du.
\end{equation*}
Hence, the second term equals $\psi_\tau(s)(\ee\phi_\tau)(t)$ and we
obtain
\begin{equation*}
  K_{\tau,11}(s,t) = S_\tau(s,t) + \psi_\tau(s)(\ee\phi_\tau)(t).
\end{equation*}

\paragraph{$K_{\tau,12}$ term.} We first recall the definition of
$K_{\tau,12}$ as
\begin{equation*}
  K_{\tau,12}(s,t) = -\tsigma_{n,N}\sqrt{\tau'(s)\tau'(t)}\partial_2
  S_{N,1}(\tau(s),\tau(t)).
\end{equation*}
For the involved partial derivative, we have
\begin{equation*}
    \partial_2 S_{N,1}(\tau(s),\tau(t)) =
    \frac{1}{\tau'(t)}\frac{\partial}{\partial t}
    \frac{K_{\tau,11}(s,t)}{\sqrt{\tau'(s)\tau'(t)}}
    = \frac{1}{\tsigma_{n,N}}\frac{\partial}{\partial t}
    \frac{S_\tau(s,t)+\psi_\tau(s)(\ee\phi_\tau)(t)}{\tsigma_{n,N}}
    = \frac{1}{\tsigma_{n,N}^2} \partial_t S^R_\tau(s,t),
\end{equation*}
with $S^R_\tau$ defined as in \eqref{eq:S-tau-R}. Observing that
$\tau'(s) = \tau'(t) = \tsigma_{n,N}$, we obtain
\begin{equation*}
  K_{\tau,12}(s,t) = -\partial_t S^R_\tau(s,t).
\end{equation*}

\paragraph{$K_{\tau,21}$ term.} By its definition, we have
\begin{equation*}
  K_{\tau,21}(s,t) =
  \frac{\sqrt{\tau'(s)\tau'(t)}}{\tsigma_{n,N}}\left[\ee
    S_{N,1}(\tau(s),\tau(t)) - \ee(\tau(s)-\tau(t))\right].
\end{equation*}
Observing that $\tau$ is a monotone transformation, we obtain
\begin{equation*}
  \ee(\tau(s) - \tau(t)) = \ee(s-t).
\end{equation*}
For the quantity $\ee S_{N,1}(\tau(s),\tau(t))$, by using the above identity, we have
\begin{equation*}
    \ee S_{N,1}(\tau(s),\tau(t)) = \int
    S_{N,1}(\tau(u),\tau(t))\ee(\tau(s)-\tau(u))\tau'(u)du
    = \int S^R_\tau(u,t)\ee(s-u)du = (\ee S^R_\tau)(s,t).
\end{equation*}
Plugging all these identities back into the definition of
$K_{\tau,21}$, we obtain the expression
\begin{equation*}
  K_{\tau,21}(s,t) = (\ee S^R_\tau)(s,t) - \ee(s-t).
\end{equation*}

\paragraph{$K_{\tau,22}$ term.} The formula for $K_{\tau,22}$ is
obtained directly from that of $K_{\tau,11}$ by switch $s$ and $t$.

\subsection{Behavior of $R_N(\xi), R_N'(\xi)$ and $\kk_N^{2/3}\zz$}
\label{subsec:a-auxfunc}

In this part, we investigate the behavior of $R_N(\xi)$, $R_N'(\xi)$
and $\kk_N^{2/3}\zz$ which is essential in deriving the Laguerre
asymptotics. Before we start, we remark that throughout our
discussion, we consider only the case where $s\in I_{1,N} = [s_0, s_1N^{1/6})$.

\subsubsection{Properties of $R_N(\xi)$ and $R_N'(\xi)$}
\label{subsubsec:R-N}

Recall the definition $R_N(\xi) =
(\dot{\zz}(\xi)/\dot{\zz}(\xi_+))^{1/2}$, we obtain that
\begin{equation}
  \label{eq:R-N'}
  R_N(\xi_+) = 1,\quad \text{and}\quad R_N'(\xi) =
  -\frac{1}{2}\dot{\zz}_N^{1/2}\dot{\zz}(\xi)^{-3/2}\ddot{\zz}(\xi).
\end{equation}

By our derivation in the Liouville-Green approximation, we know that
$\xi = \xi_+ +s\sigma_{n,N}/\kk_N$ and as has been shown before, when
$N\geq N_0(\gamma)$,
\begin{equation*}
  N^{2/3}\left(\frac{\sigma_{n,N}}{\kk_N}\right)\leq
  4\left(1+1/\sqrt{\gamma}\right)^{1/3}\left(1+\sqrt{\gamma}\right)(1+\gamma)^{-1}
  \leq C.
\end{equation*}
As $N\to \infty$, we have
\begin{equation}
  \label{eq:xi-conv}
  \sup_{s\in I_{1,N}}\labs \xi-\xi_+ \rabs = \Oh{s_1N^{-1/2}}\to 0.
\end{equation}
We then have the following first order Taylor expansion
\begin{equation*}
  R_N(\xi) = R_N(\xi_+) + R_N'(\xi^*)(\xi-\xi_+),\quad \text{for some
    $\xi^*\in [\xi\wedge \xi_+, \xi\vee \xi_+]$.}
\end{equation*}
Hence, when $N\geq N_0(s_0, \gamma)$, we have
\begin{equation}
  \label{eq:R-N-minus-1}
  \labs R_N(\xi)-1\rabs \leq \labs R_N'(\xi^*) \rabs
  4\left(1+1/\sqrt{\gamma}\right)^{1/3}\left(1+\sqrt{\gamma}\right)(1+\gamma)^{-1}
  N^{-2/3}|s|.
\end{equation}

In order to bound $|R_N(\xi)-1|$ uniformly on $I_{1,N}$ and also of
its own interest, we are to derive a bound for $|R_N'(\xi)|$ by some
constant that does not depend on $N$ and is uniform for $s\in
I_{1,N}$. By the definition of $R_N'(\xi)$ in \eqref{eq:R-N'}, this
relies on the understanding of the quantities $\dot{\zz}_N,
\dot\zz(\xi)$ and $\ddot\zz(\xi)$.

First, we consider the asymptotics of $\dot\zz_N$. Using the
notation $m_\pm = m\pm 1/2$, we obtain from simple calculation that
as $N\to \infty$,
\begin{equation}
  \label{eq:dot-zz-N-asym}
  \dot\zz_N =
  \frac{n_+^{1/6}N_+^{1/6}\left(n_++N_+\right)^{1/3}}{2^{1/3}\left(\sqrt{n_+}+\sqrt{N_+}\right)^{4/3}}
  \longrightarrow
  \frac{\gamma^{1/6}(1+\gamma)^{1/3}}{2^{1/3}\left(1+\sqrt{\gamma}\right)^{4/3}}.
\end{equation}

Second, we check the behavior of $\dot\zz(\xi)$. For simplicity, we let
$\xi_\pm^\infty = \lim_{N\to \infty}\xi_\pm$ and simple manipulation
gives us
\begin{equation*}
  \xi_+^\infty = 2\left( 1+\sqrt{\gamma} \right)^2/(1+\gamma),\quad
  \text{and} \quad \xi_+^\infty - \xi_-^\infty =
  8\sqrt{\gamma}/(1+\gamma).
\end{equation*}

We assume first that $s_0\geq 0$. By the definition of $\dot\zz(\xi)$
for $\xi\geq \xi_+$, we recognize it as
\begin{equation*}
\dot\zz(\xi) = \left[
  \frac{3}{2}\int_{\xi_+}^{\xi}\left(\frac{z-\xi_-}{\xi_+-\xi_-}\right)^{1/2}\frac{\xi\sqrt{z-\xi_+}dz}{z}
\frac{\sqrt{\xi_+-\xi_-}}{2\xi}\right]^{-1/3} \frac{\sqrt{(\xi-\xi_+)(\xi-\xi_-)}}{2\xi}.
\end{equation*}
When $s\in I_{1,N}$ with $s_0\geq 0$, we always have the bounds
\begin{equation}
  \label{eq:ratio-bd-1}
  1\leq \left(\frac{z-\xi_-}{\xi_+-\xi_-}\right)^{1/2}\frac{\xi}{z}
  \leq \left(\frac{\xi-\xi_-}{\xi_+-\xi_-}\right)^{1/2}\frac{\xi}{\xi_+}.
\end{equation}
Plugging these bounds into our modification of $\dot\zz(\xi)$, we
obtain the lower and upper bounds for $\dot\zz(\xi)$ as
\begin{equation*}
  \frac{\xi_+^{1/3}(\xi-\xi_-)^{1/3}}{2^{2/3}\xi} \leq \dot\zz(\xi)
  \leq \frac{(\xi-\xi_-)^{1/2}}{2^{2/3}\xi^{2/3}(\xi_+-\xi_-)^{1/6}},
\end{equation*}
where as $N\to \infty$, both bounds converge to the same limit:
\begin{equation*}
  \frac{(\xi_+^\infty -
    \xi_-^\infty)^{1/3}}{2^{2/3}(\xi^\infty_+)^{2/3}} = \lim_{N\to
    \infty} \dot{\zz}_N.
\end{equation*}
We remark that because of \eqref{eq:xi-conv}, the convergence is
uniform on $I_{1,N}$, which is crucial for deriving finite $N$ bounds
from the limit.

If $s_0<0$, we only need to consider the case where $s\in [s_0, 0]$,
for the case where $s\geq 0$ has essentially been considered in the
above derivation. When $s\in [s_0, 0]$, the definition of
$\dot\zz(\xi)$ is changed to
\begin{equation*}
  \begin{split}
    \dot\zz(\xi) & =
    \left(\frac{3}{2}\int_{\xi}^{\xi_+}\frac{\sqrt{(\xi_+-z)(z-\xi_-)}}{2z}dz\right)^{-1/3}\frac{\sqrt{(\xi_+-\xi)(\xi-\xi_-)}}{2\xi}\\
    & =
    \left[\frac{3}{2}\int_{\xi}^{\xi_+}\left(\frac{z-\xi_-}{\xi_+-\xi_-}\right)^{1/2}\frac{\xi\sqrt{\xi_+- z}dz}{z} \frac{\sqrt{\xi_+-\xi_-}}{2\xi}\right]^{-1/3}
    \frac{\sqrt{(\xi_+-\xi)(\xi-\xi_-)}}{2\xi}
  \end{split}
\end{equation*}
In this case, we have for all $s\in [s_0, 0]$,
\begin{equation}
  \label{eq:ratio-bd-2}
  \left(\frac{\xi-\xi_-}{\xi_+-\xi_-}\right)^{1/2}\frac{\xi}{\xi_+}\leq
  \left(\frac{z-\xi_-}{\xi_+-\xi_-}\right)^{1/2}\frac{z}{\xi_+}\leq 1.
\end{equation}
We notice that all the bounds tend to $1$ when $N\to \infty$. Hence,
plugging these bounds to our modification of $\dot\zz(\xi)$, we
obtain the lower and upper bounds for it that tend to the same limit
as when $s_0\geq 0$. Thus, we conclude for $\dot\zz(\xi)$ that when
$N\geq N_0(s_0,\gamma)$,
\begin{equation}
  \label{eq:dot-zz-xi}
  C_1 \lim_{N\to \infty}\dot\zz_N \leq \dot{\zz}(\xi) \leq
  C_2 \lim_{N\to \infty}\dot\zz_N, \quad \text{for all $s\in
    I_{1,N}$.}
\end{equation}
Such a derivation is valid, since the convergence to the limit is
uniform for $s\in I_{1,N}$.

Finally, we study the behavior of $\ddot\zz(\xi)$. To this end, we
first derive a convenient representation for it. By the definition of
$\zz$, we have $(\dot\zz)^2 = f\zz^{-1}$.  We then take derivative
with respect to $\xi$ on both sides and collect to get
\begin{equation*}
  \ddot\zz = \frac{f'\zz - f\dot\zz}{2\dot\zz\zz^2}.
\end{equation*}
Furthermore, we plug in $\zz = f/\dot\zz^2$ and obtain the final
representation as
\begin{equation*}
  \ddot\zz = \frac{f'\dot\zz-\dot\zz^4}{2f}.
\end{equation*}
Noticing the definition of $f$, we could regard the above
representation as the product of three factors: $\dot\zz(\xi)$,
$(f'(\xi)-\dot\zz(\xi)^3)/(\xi-\xi_+)$ and $2\xi^2/(\xi-\xi_-)$. The
first factor $\dot\zz$ has already been studied. We first investigate
the second factor: $(f'(\xi)-\dot\zz(\xi)^3)/(\xi-\xi_+)$.

As before, we start with the assumption that $s_0\geq 0$. By the
definition of $f$, we have
\begin{equation*}
  f'(\xi) = \frac{\xi-\xi_-}{4\xi^2} +
  \frac{\xi-\xi_+}{4\xi^2}-\frac{(\xi-\xi_+)(\xi-\xi_-)}{2\xi^3}.
\end{equation*}
For $f'(\xi)-\dot\zz(\xi)^3$, we consider first the quantity $\II(\xi)
= (\xi-\xi_-)/{4\xi^2}-\dot\zz(\xi)^3$. By \eqref{eq:ratio-bd-1} and
straightforward calculation, we obtain
\begin{equation*}
  \left[1 -
    \left(\frac{\xi-\xi_-}{\xi_+-\xi_-}\right)^{1/2}\right]\frac{\xi-\xi_-}{4\xi^2}
  \leq \II(\xi)
  \leq \left(1-\frac{\xi_+}{\xi}\right)\frac{\xi-\xi_-}{4\xi^2}.
\end{equation*}
Hence, we obtain that when $N\geq N_0(s_0,\gamma)$, for all $s\in I_{1,N}$
\begin{equation*}
  \labs \frac{\II(\xi)}{\xi-\xi_+}\rabs \leq
  \frac{\xi-\xi_-}{\xi^3},\quad \text{and hence}\quad
  \labs \frac{f'(\xi)-\dot\zz(\xi)^3}{\xi-\xi_+} \rabs \leq
  \frac{1}{4\xi^2}+\frac{2(\xi-\xi_-)}{\xi^3}\leq \frac{9}{4\xi^2}\leq
  \frac{C}{(\xi_+^\infty)^{2}}.
\end{equation*}
Moreover, when $N\geq N_0(s_0,\gamma)$, we could also have
\begin{equation*}
  \labs \dot\zz(\xi) \rabs \leq \left(\xi_+^\infty -
    \xi_-^\infty\right)^{1/3}(\xi_+^\infty)^{-2/3}, \quad \text{and}
  \quad \labs \frac{2\xi^2}{\xi_+-\xi_-} \rabs \leq
  \frac{4(\xi_+^\infty)^2}{\xi_+^\infty - \xi_-^\infty}.
\end{equation*}
Multiplying the three bounds, we finally obtain that when $N\geq
N_0(s_0,\gamma)$,
\begin{equation}
  \label{eq:ddot-zz}
  \labs \ddot\zz(\xi) \rabs \leq C(\xi_+^\infty)^{-2/3}(\xi_+^\infty -
  \xi_-^\infty)^{-2/3} =
  C\gamma^{-1/3}(1+\sqrt{\gamma})^{-4/3}(1+\gamma)^{4/3}, \quad
  \text{for all $s\in I_{1,N}$.}
\end{equation}

We remark that when $s_0<0$, we just focus on $s\in [s_0, 0]$. In this
case, the quantity $\II(\xi)$ becomes
\begin{equation*}
  \II(\xi) = \frac{\xi-\xi_-}{4\xi^2} - \left(
    \frac{\sqrt{(\xi_+-\xi)(\xi-\xi_-)}}{2\xi}
  \right)^3\left[\frac{3}{2}\int_{\xi}^{\xi_+}
    \left(\frac{z-\xi_-}{\xi_+-\xi_-}\right)^{1/2}
    \frac{\xi\sqrt{\xi_+-z}dz}{z}\frac{\sqrt{\xi_+-\xi_-}}{2\xi}\right]^{-1}
\end{equation*}
with \eqref{eq:ratio-bd-2} holds. Everything else follows just as in
the study of $\dot\zz(\xi)$. In particular, \eqref{eq:ddot-zz} still
holds.

Finally, by the definition of $R_N'(\xi)$ in \eqref{eq:R-N'} and our
analysis of $\dot\zz_N$, $\dot\zz(\xi)$ and $\ddot\zz(\xi)$, we have
that when $N\geq N_0(s_0,\gamma)$,
\begin{equation}
  \label{eq:R-N'-bd}
  \labs R_N'(\xi) \rabs \leq C\gamma^{-1/2}(1+\gamma),\quad \text{for
    all $s\in I_{1,N}$.}
\end{equation}
This bound, together with \eqref{eq:R-N-minus-1}, gives
\begin{equation}
  \label{eq:R-N-minus-1-bd}
  \labs R_N(\xi) - 1\rabs \leq
  C\sqrt{\gamma}(1+\sqrt{\gamma})(1+1/\sqrt{\gamma})^{1/3}N^{-2/3}s
  \leq CN^{-2/3}|s|,\quad \text{for all $s\in I_{1,N}$.}
\end{equation}

\subsubsection{Behavior of $\kk_N^{2/3}\zz$}
\label{subsubsec:kk-N-zz}

Exploiting a simple Taylor expansion at $\xi_+$ to the second order,
we obtain that
\begin{equation*}
  \begin{split}
    \kk_N^{2/3}\zz(\xi) & = \kk_N^{2/3}\zz(\mu_{n,N}/\kk_N +
    s\sigma_{n,N}/\kk_N) \\
    & = \kk_N^{2/3}\zz(\xi_+) + \kk_N^{-1/3}\sigma_{n,N}\dot\zz_Ns +
    \frac{1}{2}\kk_N^{-4/3}\sigma_{n,N}^2\ddot\zz(\xi^*)s^2.
  \end{split}
\end{equation*}
Recalling that $\zz(\xi_+) = 0$ and that
$\sigma_{n,N}\kk_N^{-1/3}\dot\zz_N = 1$, we obtain that
\begin{equation*}
  \kk_N^{2/3}\zz(\xi)-s =
  \frac{1}{2}\frac{\sigma_{n,N}}{\kk_N}\frac{\ddot\zz(\xi^*)}{\dot\zz_N}s^2.
\end{equation*}
According to our previous discussion, we have
\begin{equation*}
  \frac{\sigma_{n,N}}{\kk_N} = \Oh{N^{-2/3}}\quad \text{and} \quad
  \frac{\ddot{\zz}(\xi^*)}{\dot\zz_N} = \Oh{1}, \quad \text{for all
    $s\in I_{1,N}$.}
\end{equation*}

Hence for all $s\in I_{1,N}$, when $N\geq N_0(s_0,\gamma)$, we have
\begin{equation*}
  \labs \kk_N^{2/3}\zz(\xi) -s \rabs \leq CN^{-2/3}s^2.
\end{equation*}
Note that on $I_{1,N}$, $|s|\leq s_1N^{1/6}$ and hence we could modify
the above bound to be
\begin{equation}
  \label{eq:kk-N-zz-minus-s}
  \labs \kk_N^{2/3}\zz(\xi) -s \rabs \leq \left(CN^{-2/3}s^2\right) \wedge
  \frac{|s|}{2}
  \wedge 1,\quad \text{for all $s\in I_{1,N}$.}
\end{equation}

\subsection{Choice of $s_1$ and its consequences}
\label{subsec:a-s-1}

The key point in our choice of $s_1$ is to ensure that when $s\geq
s_1$, we have
\begin{equation}
  \label{eq:kk-N-zz-s-1}
  \frac{2}{3}\kk_N\zz^{3/2}\geq \frac{3}{2}s.
\end{equation}

To this end, recall that in \citet[A.8]{johnstone01}, one could choose
$\tilde{s}_1(\gamma) = C(\gamma)(1+\delta)$ with some $\delta>0$, such
that when $s\geq \tilde{s}_1(\gamma)$, we have $\sqrt{f(\xi)}\geq
2/\sigma_{n,N}$ and hence if $s\geq 4\tilde{s}_1(\gamma)$,
\begin{equation*}
  \frac{2}{3}\kk_N\zz^{3/2} = \kk_N\int_{\xi_+}^\xi\sqrt{f(z)}dz \geq
  \kk_N\frac{2}{\sigma_{n,N}}(s-\tilde{s}_1(\gamma))\frac{\sigma_{n,N}}{\kk_N}
  = 2(s-\tilde{s}_1(\gamma))\geq \frac{3}{2}s.
\end{equation*}
Moreover, by the analysis in \citet[A.6.4]{nek06},
$\tilde{s}_1(\gamma)$ could be chosen independently of $\gamma$ and
hence we could define our $s_1$ to be
\begin{equation*}
  s_1 = 4\tilde{s}_1
\end{equation*}
which is independent of $\gamma$ and such that \eqref{eq:kk-N-zz-s-1}
holds. Moreover, for our convenience of arguments, we could also
impose the constraint that $s_1 \geq 1$.

After specifying our choice of $s_1$, we spell out two of its
consequences. The first of them is that when $s\geq
s_1\geq 1$,
\begin{equation}
  \label{eq:EE-bd}
  \EE^{-1}(\kk_N^{2/3}\zz) \leq C\exp(-3s/2) \leq C\exp(-s).
\end{equation}
This is from the observation that $\EE(x)\geq
C\exp(\frac{2}{3}x^{3/2})$ and hence
\begin{equation*}
  \EE^{-1}(\kk_N^{2/3}\zz)\leq
  C\exp\left(-\frac{2}{3}\kk_N\zz^{3/2}\right)\leq C\exp(-3s/2).
\end{equation*}

The other consequence is about the behavior of $s'$ defined in
\eqref{eq:s-prime-def} when $s\geq s_1$. Remembering that $s_1\geq 1$,
we then have that when $s\geq s_1$ and $N\geq N_0(\gamma)$,
\begin{equation}
  \label{eq:s-prime-s/2}
  s'-\frac{s}{2} = \Delta_N + \left(\frac{\sigma_{n-1,N-1}}{\sigma_{n-2,N}} -
    \frac{1}{2}\right)s\geq \Delta_N + \frac{s_1}{2} \geq \Delta_N + \frac{1}{2} \geq 0.
\end{equation}
The last inequality holds when $N\geq N_0(\gamma)$ for $\Delta_N =
\Oh{N^{-1/3}}$.

\subsection{Proofs of Proposition \ref{prop:det-bd} and Lemma \ref{lemma:sqrt-f}}
\label{subsec:a-proof}

\subsubsection{Proof of Proposition \ref{prop:det-bd}}
By the definition \eqref{eq:det} of the determinant for operators in
class $\mathcal{A}$, we have a first decomposition as
\begin{equation*}
\begin{split}
\labs \det(I-A) - \det(I-B) \rabs \ =\ & \labs
\text{det}_2(I-A)\exp\left(-\trace{A}\right) -
\text{det}_2(I-B)\exp\left(-\trace{B}\right)
\rabs\\
\ \leq\  & \labs \text{det}_2(I-A) - \text{det}_2(I-B) \rabs
\left[\labs \exp\left(-\trace{A}\right) -
\exp\left(-\trace{B}\right) \rabs +
\exp\left(-\trace{B}\right)\right] \\
&\quad +\ \labs \text{det}_2(I-B) \rabs \labs
\exp\left(-\trace{A}\right) - \exp\left(-\trace{B}\right) \rabs.
\end{split}
\end{equation*}

According to \citet[p.69, Theorem 7.4]{gohberg00}, we have the bound
for the $2$-determinant as
\begin{equation*}
    \labs  \text{det}_2(I-A) - \text{det}_2(I-B) \rabs \leq
    \|A-B\|_2\exp\left[ \frac{1}{2}\left(1 + \|A\|_2 + \|B\|_2\right)^2 \right].
\end{equation*}
Moreover, for any $A, B\in \mathcal{A}$, the Hilbert-Schmidt norm
satisfies
\begin{equation*}
    \|A-B\|_2 \leq \sum_i\|A_{ii}-B_{ii}\|_2 + \sum_{i\neq
    j}\|A_{ij}-B_{ij}\|_2 \leq \sum_i\|A_{ii}-B_{ii}\|_1 + \sum_{i\neq
    j}\|A_{ij}-B_{ij}\|_2.
\end{equation*}
Recalling that for any trace class operator $A$, $\trace(A)\leq
\|A\|_1$, we obtain
\begin{equation*}
    \exp\left(-\trace{B}\right) \leq \exp{\labs \trace{B} \rabs}\leq \exp\left(\labs \trace{B_{11}} \rabs + \labs
    \trace{B_{22}}
    \rabs\right) \leq \exp\left(\|B_{11}\|_1 + \|B_{22}\|_1\right).
\end{equation*}

Observing that for $|x|\leq 1/2$, $\labs e^x - 1 \rabs\leq 2|x|$, we
obtain that, when $\sum_{i}\|A_{ii} - B_{ii}\|_1 + \sum_{i\neq
j}\|A_{ij} - B_{ij}\|_2\leq 1/2$,
\begin{equation*}
\begin{split}
    \labs \exp\left(-\trace{A}\right) - \exp\left(-\trace{B}\right) \rabs & \leq
    2\exp\left(-\trace{B}\right)\labs \trace{A} - \trace{B} \rabs \leq
    2\sum_i\|A_{ii}-B_{ii}\|_1 e^{\|B_{11}\|_1 + \|B_{22}\|_1} \\
    & \leq 2\left(\sum_{i}\|A_{ii} - B_{ii}\|_1 + \sum_{i\neq
j}\|A_{ij} - B_{ij}\|_2\right)\exp\left(\|B_{11}\|_1 +
\|B_{22}\|_1\right).
\end{split}
\end{equation*}

Plugging all these bounds into our first decomposition, we obtain an
intermediate bound as $M(A, B)\left( \sum_i\|A_{ii}-B_{ii}\|_1 +
\sum_{i\neq j}\|A_{ij} - B_{ij}\|_2 \right)$, where
\begin{equation*}
    M(A,B) = 2\,\labs \det(I-B) \rabs + 2\exp\left[\frac{1}{2}\left(1 + \|A\|_2 + \|B\|_2\right)^2 + \sum_i\|B_{ii}\|_1
    \right].
\end{equation*}
Under the given condition,
\begin{equation*}
\begin{split}
    1 + \|A\|_2 + \|B\|_2 & \leq 1 + 2\|B\|_2 + \|A-B\|_2 \\
    & \leq 1 + 2\|B\|_2 + \sum_i\|A_{ii}-B_{ii}\|_1 + \sum_{i\neq j}\|A_{ij} -
    B_{ij}\|_2 \leq 2 + 2\|B\|_2,
\end{split}
\end{equation*}
which reduce $M(A,B)$ to the constant $M(B)$ claimed.

\subsubsection{Proof of Lemma \ref{lemma:sqrt-f}}

By definition, we have
\begin{equation*}
    \dot{\zz}_N^3 = \frac{\xi_+ - \xi_-}{4\xi_+^2} =
    \frac{\kk_N}{\sigma_{n,N}^3}.
\end{equation*}
Thus, we obtain from direct calculation that
\begin{equation*}
    \sqrt{f(\xi)} = \frac{(\xi-\xi_+)(\xi-\xi_-)}{2\xi} \geq
    r\sqrt{\frac{\sigma_{n,N}}{\kk_N}}\frac{\sqrt{\xi_+-\xi_-}}{2\xi_+}\frac{\xi_+}{\xi}
    =
    r\sqrt{\frac{\sigma_{n,N}}{\kk_N}}\sqrt{\frac{\kk_N}{\sigma_{n,N}^3}}\frac{\xi_+}{\xi}
    = \frac{r\xi_+}{\sigma_{n,N}\xi}.
\end{equation*}

\section{Logarithmic Transformation and the Smallest Eigenvalue}
\label{appendix-b}

In this part, we give a brief account of how one could derive the
similar second order accuracy results claimed in \eqref{eq:res-log}
and \eqref{eq:res-log-smallest} with a logarithmic transformation.
In many aspects, the derivation here for Laguerre orthogonal
ensembles [as based on \citet[Proposition 4.2]{adler00}] is parallel
to what \citet{johnstone07} did for Jacobi orthogonal ensembles.

\subsection{Logarithmic transformation for the largest eigenvalue}
\label{sec:log-trans}

For the largest eigenvalue, we assume the same setting as that in
the beginning of Section \ref{sec:determ-form-orth}. With $\phi_k$
defined in \eqref{eq:phi-k}, let
\begin{equation}\label{eq:phi-bar-k}
    \bar{\phi}_k(x;\tilde{\aa}) =
    (-1)^j\phi_k(x;\tilde{\aa})/\sqrt{x}.
\end{equation}
Then setting $a_N = \sqrt{N(N+\aa_N-1)}$, we have the following
alternative way of expressing $S_{N,1}$ in term of $S_{k,2}$, the
correlation kernel occurring in LUE$(k,\tilde{\aa})$ model:
\begin{equation}\label{eq:S-N-1-alt}
S_{N,1}(x,y;\aa_N - 1) = \sqrt{\frac{y}{x}}S_{N-1,2}(x,y;\aa_N) +
\sqrt{\frac{N-1}{N}}\frac{a_N}{2}\,\bar{\phi}_{N-1}(x;\aa_N)(\ee\bar{\phi}_{N-2})(y;\aa_N).
\end{equation}
As a comparison, the central formula \eqref{eq:central} could be
rewritten as
\begin{equation*}
    S_{N,1}(x,y;\aa_N - 1) =
S_{N,2}(x,y;\aa_N - 1) +
\frac{a_N}{2}\,\bar{\phi}_{N-1}(x;\aa_N)(\ee\bar{\phi}_N)(y;\aa_N-2).
\end{equation*}
The equivalence of the above two representations is given in
\citet[Appendix]{adler00} and hence omitted here.

We make use of the representation \eqref{eq:S-N-1-alt} to give an
alternative second order accuracy argument with a logarithmic
transformation. Recalling $\aa_N = n - N$, we define
\begin{equation*}
    \mu_k = \left(\sqrt{k+\tfrac{1}{2}} + \sqrt{k + \aa_N +
    \tfrac{1}{2}}\right)^2, \sigma_k = \left(\sqrt{k+\tfrac{1}{2}} + \sqrt{k + \aa_N +
    \tfrac{1}{2}}\right)\left(\frac{1}{\sqrt{k+\tfrac{1}{2}}} + \frac{1}{\sqrt{k + \aa_N +
    \tfrac{1}{2}}}\right)^{1/3}.
\end{equation*}
Then we let
\begin{equation*}
    \check{\phi}_k(x) =
    (-1)^k\frac{(N-1)^{1/4}(N-1+\aa_N)^{1/4}}{\sqrt{2}}x^{1/2}\phi_k(x;\aa_N).
\end{equation*}
For $\hat{S}_{N-1,2}(u,v;\aa_N) = S_{N-1,2}(e^u,
e^v;\aa_N)\,e^{u/2}e^{v/2}$, we could represent it as
\begin{equation*}
    \hat{S}_{N-1,2}(u,v; \aa_N) = \int_0^\infty\left[
    \check{\phi}_{N-2}(e^{u+w})\check{\phi}_{N-1}(e^{v+w})
    + \check{\phi}_{N-1}(e^{u+w})\check{\phi}_{N-2}(e^{v+w})
    \right] dw.
\end{equation*}

We then define
\begin{equation*}
    \nu_{n,N} = \log\mu_{N-1},\quad \tau_{n,N}
    =\sigma_{N-1}/\mu_{N-1},\quad\text{and}\quad \tau(s) = \exp(\nu_{n,N} +
    s\tau_{n,N}).
\end{equation*}
The $\tau$-transformation induces the following transformed Laguerre
polynomials:
\begin{equation*}
    \psi_\tau(s) = \tau_{n,N} \check{\phi}_{N-1}(\tau(s)),\quad
    \phi_\tau(s) = \tau_{n,N} \check{\phi}_{N-2}(\tau(s)).
\end{equation*}

Define $S_\tau(s,t) =
\sqrt{\tau'(s)\tau'(t)}\,S_{N-1,2}(\tau(s),\tau(t);\aa_N)$, we have
the following integral representation from the expression for the
$\hat{S}_{N-1,2}$ kernel:
\begin{equation}\label{eq:s-tau-int}
    S_\tau(s,t) = \int_0^\infty
    [\phi_\tau(s+z)\psi_\tau(t+z)+\phi_\tau(t+z)\psi_\tau(s+z)]dz.
\end{equation}
Moreover, if we define the following quantities [fix $s_0\in
\mathbb{R}$, with $s,t\geq s_0$]
\begin{equation*}
    q_{N}(s) = \sqrt{\tau'(s_0)/\tau'(s)},\quad \text{and} \quad
    S^R_\tau(s,t) = S_\tau(s,t) + \psi_\tau(s)(\ee \phi_\tau)(t),
\end{equation*}
we have
\begin{equation*}
    F_{N,1}(s') = P(x_1\leq \tau(s')) = P((\log x_1-\nu_{n,N})/\tau_{n,N}
    \leq s') = \sqrt{\det(I-K_\tau)},
\end{equation*}
where the new operator $K_\tau$ has a $2\times 2$ matrix kernel with
entries given by
\begin{equation}
\label{eq:k-tau-entries-log}
\begin{split}
    K_{\tau,11}(s,t) =
    q_{N}(s)q_{N}^{-1}(t)S^R_\tau(s,t); \qquad &
    K_{\tau,12}(s,t) = -q_{N}(s)q_{N}(t)\partial_t
    S^R_\tau(s,t); \\
    K_{\tau,21}(s,t) = q_{N}^{-1}(s)q_{N}^{-1}(t)[\ee_1
    S^R_\tau(s,t)-\ee(s-t)]; \qquad &
    K_{\tau,22}(s,t)  = K_{\tau,11}(t,s).
\end{split}
\end{equation}

By Proposition \ref{prop:det-bd}, we need to obtain entrywise bound
for $K_\tau - K_{GOE}$ here. To this end, a convenient
representation of the kernel difference as in Section
\ref{sec:repr-kern-diff} is most helpful.

For the transformed Laguerre polynomials $\phi_\tau$ and
$\psi_\tau$, we have
\begin{equation*}
    \int_{-\infty}^\infty \psi_\tau  = 0,\quad \text{and}\quad
    \int_{-\infty}^\infty \phi_\tau  = \frac{(N-1)^{1/4}(n-1)^{1/4}\Gamma^{1/2}(n-1)
    \Gamma\left(\frac{N+1}{2}\right)}{2^{\aa_N -
    2}(N-1)\Gamma^{1/2}(N-1)\Gamma\left(\frac{n}{2}\right)}.
\end{equation*}
For notational convenience, let $\tilde\beta_N =
\frac{1}{2}\int_{-\infty}^\infty \phi_\tau = \frac{1}{\sqrt{2}} +
\Oh{N^{-1}}$.

With the replacement of $\ee$ by $\tee$ in \eqref{eq:tee} and the
matrices $\tilde{L}, L_1$ and $L_2$ introduced in Section
\ref{sec:repr-kern-diff}, we obtain that
\begin{equation*}
\begin{split}
    K_\tau  = Q_N(s)\left[ K^R_\tau + K^F_{\tau,1} + K^F_{\tau,2} + K^\ee
    \right]Q_N^{-1}(t).
\end{split}
\end{equation*}
with the unspecified components given by
\begin{align*}
    K^R_{\tau} &= \tilde{L}[S_\tau - \psi_\tau\otimes
    \tee\phi_\tau], &\quad & K^F_{\tau,1} = \tilde\beta_N
    L_1[\psi_\tau(s)],& \quad & K^F_{\tau,2} = \tilde\beta_N
    L_2[\psi_\tau(t)],
\end{align*}
where $Q_N(s) = \text{diag}(q_N(s), q_N^{-1}(s))$ and as before
$G(s) = \Ai(s)/\sqrt{2}$.

For
\begin{equation*}
    \Delta_N = \frac{\log\mu_{N-1}
    -\log\mu_{N-2}}{\sigma_{N-2}/\mu_{N-2}},
\end{equation*}
set $G_N = G + \Delta_N G'$, we have
\begin{equation*}
    \begin{split}
    K^R_\tau - K_R & = \tilde{L}\left[ S_\tau - S_A - \psi_\tau\otimes\tee\phi_\tau
    + G\otimes \tee(G_N - \Delta_N G')\right] \\
    & = \tilde{L}\left[S_\tau - S_A + \Delta_N G\otimes G\right] -
    \tilde{L}\left[\psi_\tau\otimes \tee\phi_\tau - G\otimes \tee
    G_N\right] \\
    & = \delta^{R,I} + \delta^F_0.
    \end{split}
\end{equation*}

If we write $S_{A_N} = G\diamond G_N + G_N\diamond G$, we have
$\delta^{R,I} = \tilde{L}[S_\tau - S_{A_N}]$.

Finally, we organize $K_\tau - K_{GOE}$ as
\begin{equation*}
    K_\tau - K_{GOE} = \delta^{R,D} + \delta^{R,I} + \delta^F_0 +
    \delta^F_1 + \delta^F_2 + \delta^\ee,
\end{equation*}
where the unspecified terms are defined as the following:
\begin{equation*}
    \begin{split}
    \delta^{R,D} & = Q_N(s)K^R_\tau Q_N^{-1}(t) - K^R_\tau,\\
    \delta^F_i & = Q_N(s)K^F_{\tau,i}Q_N^{-1}(t) - K^F_i,\quad \text{for $i = 1,2$,
    and}\\
    \delta^\ee & = Q_N(s)K^\ee Q_N^{-1}(t) - K^\ee.
    \end{split}
\end{equation*}

With the above representation of the kernel difference, we could
apply the machineries in \citet{johnstone07} to obtain the desired
second order accuracy of the Tracy-Widom approximation to the
distribution of $(\log x_1 - \nu_{n,N})/\tau_{n,N}$. After
establishing the result in RMT notation, we replace $N$ by $p$ and
hence obtain the bound in \eqref{eq:res-log}.

\subsection{The smallest eigenvalue}
\label{sec:smallest}

We first restate the claim in \eqref{eq:res-log-smallest} in a more
friendly way. Let $\nu_{n,N}^-$ and $\tau_{n,N}^-$ be the centering
and scaling constants defined in \eqref{eq:small-center-scale}, with
$p$ replaced by $N$. Then for $x_N$ the smallest eigenvalue in the
model \eqref{eq:loe-pdf} [with $\tilde{\aa} = \aa_N - 1$], there
exists a continuous and nonincreasing function $C(\cdot)$, such that
for all real $s_0$, there is an integer $N_0(s_0,\gamma)$ for which
we have that for any $s\geq s_0$ and $N\geq N_0(s_0,\gamma)$,
\begin{equation}
\label{eq:res-log-smallest-1}
    \labs P\{\log x_N > \nu_{n,N}^- - s\tau_{n,N}^-\} - F_1(s) \rabs
    \leq C(s_0)N^{-2/3}\exp(-s/2).
\end{equation}

Fix $x_0 \geq 0$ and consider any $x' \in [0,x_0]$. To prove
\eqref{eq:res-log-smallest-1}, we first observe that for $x_N$ in
model \eqref{eq:loe-pdf}, when $N$ is even, choosing $\chi =
I_{0\leq x\leq x'}$, we have
\begin{equation*}
    P\{x_N > x'\} = \sqrt{\det(I-K_N\chi)},
\end{equation*}
where $K_N$ is the same operator as for $x_1$, which has the kernel
\eqref{eq:K-N-kernel}. If we think of $K_N$ as Hilbert-Schmidt
operator on $L^2([0,x'];\rho)\oplus L^2([0,x'];\rho^{-1})$ with
$\rho$ any weight function chosen from some proper class, then the
above formula changes to
\begin{equation*}
    P\{x_N > x'\} = \sqrt{\det(I-K_N)}.
\end{equation*}

Introduce the transformation
\begin{equation*}
    \tau(s) = \exp(\nu_{n,N}^- - s\tau_{n,N}^-),
\end{equation*}
and let $s_0 = \tau^{-1}(x_0)$ and $s_0\leq s' = \tau^{-1}(x')$, we
have $\tau^{-1}([0,x']) = [s',\infty)$. By defining $\phi_\tau =
-\tau_{n,N}^- \check{\phi}_{N-1}(\tau(s))$ and $\psi_\tau =
-\tau_{n,N}^- \check{\phi}_{N-2}(\tau(s))$ and using the alternative
representation \eqref{eq:S-N-1-alt}, the formal derivation for the
largest eigenvalue in \ref{sec:log-trans} could be carried out
analogously for the smallest eigenvalue. In particular, we have the
integral representation \eqref{eq:s-tau-int} for $S_\tau(s,t) =
\sqrt{\tau'(s)\tau'(t)} S_{N-1,2}(\tau(s),\tau(t);\aa_N)$ and
\begin{equation*}
    P(x_N > \tau(s')) = P((\log x_N-\nu_{n,N}^-)/\tau_{n,N}^- >
    -s') = \sqrt{\det(I-K_\tau)},
\end{equation*}
with $K_\tau$ thought of as Hilbert-Schmidt operator on
$L^2([s',\infty);\rho\circ\tau)\oplus L^2([s',\infty);\rho^{-1}\circ
\tau)$ with entries given by \eqref{eq:k-tau-entries-log}. We remark
that the actual definition of $\phi_\tau$ and $\psi_\tau$ used in
these formulas have changed, albeit the formal representations
remain the same.

The rest of the proof for the smallest eigenvalue becomes the
routine procedure of a) finding a representation for the kernel
difference $K_\tau - K_{GOE}$ and b) studying the asymptotic
behavior of the transformed Laguerre polynomials $\phi_\tau$ and
$\psi_\tau$. The former is very similar to the largest eigenvalue
case while the latter could be obtained by applying the
Liouville-Green approach to analyze the behavior of the solution to
the differential equation \eqref{eq:diff-eqn} around the lower
turning point $\xi_-$.

\section*{Acknowledgment}
The author is most grateful to Professor Iain Johnstone for his
indispensable advice during the development of this project. The
author would also like to thank Debashis Paul for sharing a draft of
his paper on the smallest eigenvalue. This work is supported in part
by grants NSF DMS 0505303 and NIH EB R01 EB001988.

\bibliographystyle{plainnat}

\bibliography{accuracybib}

\begin{thebibliography}{33}
\providecommand{\natexlab}[1]{#1}
\providecommand{\url}[1]{\texttt{#1}}
\expandafter\ifx\csname urlstyle\endcsname\relax
  \providecommand{\doi}[1]{doi: #1}\else
  \providecommand{\doi}{doi: \begingroup \urlstyle{rm}\Url}\fi

\bibitem[Adler et~al.(2000)Adler, Forrester, Nagao, and van Moerbeke]{adler00}
M.~Adler, P.~J. Forrester, T.~Nagao, and P.~van Moerbeke.
\newblock Classical skew orthogonal polynomials and random matrices.
\newblock \emph{J. Statist. Phys.}, 99:\penalty0 141--170, 2000.

\bibitem[Anderson(2003)]{anderson}
T.~W. Anderson.
\newblock \emph{An Introduction to Multivariate Statistical Analysis}.
\newblock John Wiley and Sons, 3rd edition, 2003.

\bibitem[Candes and Tao(2006)]{candestao06}
E.~Candes and T.~Tao.
\newblock Near optimal signal recovery from random projections: Universal
  encoding strategies?
\newblock \emph{IEEE Trans. Inform. Theory}, 52:\penalty0 5406--5425, 2006.

\bibitem[de~Bruijn(1955)]{db55}
N.~G. de~Bruijn.
\newblock On some multiple integrals involving determinants.
\newblock \emph{J. Indian Math. Soc.}, 19:\penalty0 133--151, 1955.

\bibitem[Donoho(2004)]{donoho04}
D.~L. Donoho.
\newblock For most large underdetermined systems of linear equations the
  minimal $\ell^1$-norm solution is also the sparsest solution.
\newblock 2004.

\bibitem[Dumitriu and Edelman(2002)]{de02}
I.~Dumitriu and A.~Edelman.
\newblock Matrix models for beta ensembles.
\newblock \emph{J. Math. Phys.}, 43\penalty0 (11):\penalty0 5830--5847, 2002.

\bibitem[Edelman and Persson(2002)]{persson05}
A.~Edelman and P.-O. Persson.
\newblock Numerical methods for eigenvalue distributions of random matrices.
\newblock Technical report, Massachusetts Institute of Technology, 2002.

\bibitem[El~Karoui(2006{\natexlab{a}})]{nek-inf}
N.~El~Karoui.
\newblock On the largest eigenvalue of {W}ishart matrices with identity
  covariance when $n, p$ and $p/n\to \infty$.
\newblock \emph{arXiv:math/03093355v1}, 2006{\natexlab{a}}.

\bibitem[El~Karoui(2006{\natexlab{b}})]{nek06}
N.~El~Karoui.
\newblock A rate of convergence result for the largest eigenvalue of complex
  white {W}ishart matrices.
\newblock \emph{Ann. Probab.}, 34:\penalty0 2077--2117, 2006{\natexlab{b}}.

\bibitem[Forrester(2004)]{pjfbook}
P.~J. Forrester.
\newblock Log-gases and random matrices.
\newblock Book manuscript, 2004.

\bibitem[Gohberg et~al.(2000)Gohberg, Goldberg, and Krupnik]{gohberg00}
I.~Gohberg, S.~Goldberg, and N.~Krupnik.
\newblock \emph{Traces and Determinants of Linear Operators}, volume 116 of
  \emph{Operator Theory, Advances and Applications}.
\newblock Birkh\"{a}user Verlag, Basel, 2000.

\bibitem[Golub and van Loan(1996)]{golub}
G.~H. Golub and C.~F. van Loan.
\newblock \emph{Matrix Computations}.
\newblock The Johns Hopkins University Press, 3rd ed. edition, 1996.

\bibitem[Johansson(2000)]{johansson00}
K.~Johansson.
\newblock Shape fluctuations and random matrices.
\newblock \emph{Comm. Math. Phys.}, 209:\penalty0 437--476, 2000.

\bibitem[Johnstone(2001)]{johnstone01}
I.~Johnstone.
\newblock On the distribution of the largest eigenvalue in principal component
  analysis.
\newblock \emph{Ann. Statist.}, 29:\penalty0 295--327, 2001.

\bibitem[Johnstone(2006)]{johnstone06}
I.~Johnstone.
\newblock High dimensional statistical inference and random matrices.
\newblock \emph{arXiv:math/0611589}, 2006.

\bibitem[Johnstone(2007)]{johnstone07}
I.~Johnstone.
\newblock Canonical correlation analysis and {J}acobi ensembles: {T}racy
  {W}idom limits and rates of convergence.
\newblock Unpublished manuscript, 2007.

\bibitem[Koev and Edelman(2006)]{koev}
P.~Koev and A.~Edelman.
\newblock The efficient evaluation of the hypergeometric function of a matrix
  argument.
\newblock \emph{Math. Comp.}, 75:\penalty0 833--846, 2006.

\bibitem[Luca et~al.(2008)Luca, Ringquist, Klei, Lee, Gieger, Wichmann,
  Schreiber, Krawczak, Lu, Styche, Devlin, Roeder, and Trucco]{luca}
D.~Luca, S.~Ringquist, L.~Klei, A.~B. Lee, C.~Gieger, H.~E. Wichmann,
  S.~Schreiber, M.~Krawczak, Y.~Lu, A.~Styche, B.~Devlin, K.~Roeder, and
  M.~Trucco.
\newblock On the use of general control samples for genome-wide association
  studies: genetic matching highlights causal variants.
\newblock \emph{Am. J. Hum. Genet.}, 82:\penalty0 453--463, 2008.

\bibitem[Mardia et~al.(1979)Mardia, Kent, and Bibby]{mkb}
K.~V. Mardia, J.~T. Kent, and J.~M. Bibby.
\newblock \emph{Multivariate Analysis}.
\newblock Academic Press, 1979.

\bibitem[Muirhead(1982)]{muirhead}
R.~J. Muirhead.
\newblock \emph{Aspects of Multivariate Statistical Theory}.
\newblock John Wiley and Sons, 1982.

\bibitem[Nagao and Forrester(1995)]{nagao95}
T.~Nagao and P.~J. Forrester.
\newblock Asymptotic correlations at the spectrum edge of random matrices.
\newblock \emph{Nucl. Phys. B}, 435:\penalty0 401--420, 1995.

\bibitem[Olver(1974)]{olver74}
F.~W.~J. Olver.
\newblock \emph{Asymptotics and Special Functions}.
\newblock Academic Press, 1974.

\bibitem[Patterson et~al.(2006)Patterson, Price, and Reich]{patterson}
N.~Patterson, A.~L. Price, and D.~Reich.
\newblock Population structure and eigenanalysis.
\newblock \emph{PLoS Genet.}, 2:\penalty0 e190, 2006.
\newblock \doi{10.1371/journal.pgen.0020190}.

\bibitem[Paul(2006)]{paul06}
D.~Paul.
\newblock Distribution of the smallest eigenvalue of {W}ishart$({N},n)$ when
  ${N}/n\to 0$.
\newblock Unpublished manuscript, 2006.

\bibitem[Price et~al.(2006)Price, Patterson, Plenge, Weinblatt, Shadick, and
  Reich]{price}
A.~L. Price, N.~J. Patterson, R.~M. Plenge, M.~E. Weinblatt, N.~A. Shadick, and
  D.~Reich.
\newblock Principal components analysis corrects for stratification in
  genome-wide association studies.
\newblock \emph{Nat. Genet.}, 38:\penalty0 904--909, 2006.

\bibitem[Reed and Simon(1980)]{reed&simon80}
M.~Reed and B.~Simon.
\newblock \emph{Methods of Modern Mathematical Physics. Vol. I: Functional
  Analysis}.
\newblock Academic Press, 1980.

\bibitem[Roy(1953)]{roy53}
S.~N. Roy.
\newblock On a heuristic method of test construction and its use in
  multivariate analysis.
\newblock \emph{Ann. Math. Stat.}, 24:\penalty0 220--238, 1953.

\bibitem[Szeg\"{o}(1975)]{szego}
G.~Szeg\"{o}.
\newblock \emph{Orthogonal Polynomials}.
\newblock Amer. Math. Soc., Providence, RI., 4th edition, 1975.

\bibitem[Tracy and Widom(2005)]{tw05}
C.~A. Tracy and H.~Widom.
\newblock Matrix kernels for the {G}aussian orthogonal and symplectic
  ensembles.
\newblock \emph{Ann. Institut. Fourier, Grenoble}, 55:\penalty0 2197--2207,
  2005.

\bibitem[Tracy and Widom(1994)]{tw94}
C.~A. Tracy and H.~Widom.
\newblock Level-spacing distributions and the {A}iry kernel.
\newblock \emph{Commun. Math. Phys.}, 159:\penalty0 151--174, 1994.

\bibitem[Tracy and Widom(1996)]{tw96}
C.~A. Tracy and H.~Widom.
\newblock On orthogonal and symplectic matrix ensembles.
\newblock \emph{Commun. Math. Phys.}, 177:\penalty0 727--754, 1996.

\bibitem[Tracy and Widom(1998)]{tw98}
C.~A. Tracy and H.~Widom.
\newblock Correlation functions, cluster functions, and spacing distributions
  for random matrices.
\newblock \emph{J. Statist. Phys.}, 92:\penalty0 809--835, 1998.

\bibitem[Widom(1999)]{widom99}
H.~Widom.
\newblock On the relation between orthogonal, symplectic and unitary matrix
  ensembles.
\newblock \emph{J. Statist. Phys.}, 94:\penalty0 347--364, 1999.

\end{thebibliography}


\end{document}